\documentclass[]{amsart}
\usepackage{amsfonts,amssymb,epsfig}
\usepackage[latin1]{inputenc}
\usepackage[all]{xy}
\usepackage{amsmath}
\usepackage{delarray}
\usepackage{amssymb}
\usepackage{amscd}
\usepackage[mathscr]{eucal}
\usepackage{diagrams}
\usepackage{amsxtra}

\newtheorem{thm}{Theorem}[subsubsection]

\theoremstyle{definition}

\newtheorem{rem}[thm]{Remark}

\numberwithin{equation}{subsection}

\newcommand{\cD}{\mathcal D}

\newcommand{\cT}{\mathcal T}

\newcommand{\cS}{\mathcal S}

\newcommand{\bbR}{\mathbb R}
\newcommand{\bbT}{\mathbb T}
\newcommand{\bbC}{\mathbb C}
\newcommand{\bbN}{\mathbb N}

\newcommand{\bbZ}{\mathbb Z}

\newcommand{\bbD}{\mathbb D}
\newcommand{\bbI}{\mathbb I}

\newcommand{\St}{$\check{\textrm{S}}$tefan}
\newcommand{\sD}{\mathsf D}
\newcommand{\sT}{\mathsf T}
\newcommand{\sS}{\mathsf S}
\newcommand{\sF}{\mathsf F}
\newcommand{\sE}{\mathsf E}

\newcommand{\gP}{\mathfrak{P}}

\newcommand{\sov}[1]{\ls\negthickspace\raisebox{.5ex}{${\un{\ssst{#1}}{\bigtriangledown}}$}}
\newcommand{\ptop}[1]{\pi^{\text{top}}_{#1}}
\newcommand{\itop}[1]{\iota^{\text{top}}_{#1}}

\newcommand{\gr}{\mathcal}

\newcommand{\sst}{\scriptstyle}
\newcommand{\ssst}{\scriptscriptstyle}

\newcommand{\bdi}{\begin{diagram}}
\newcommand{\edi}{\end{diagram}}

\newcommand{\ls}{\lesssim}

\def\ttt{\stackrel{\circlearrowleft}{\longleftrightarrow}}

\newarrow{TO}----{curlyvee}

\newarrow{TeXto}----{->}
\newarrow{ClopTo}{sqhook}---{vee}
\newarrow{OpTo}{hooka}---{vee}
\newarrow{iTo}{triangle}---{curlyvee}

\newarrow{ITo}{blacktriangle}---{curlyvee}
\newarrow{sTo}----{triangle}
\newarrow{STo}----{blacktriangle}
\newarrow{Line}-----%

\newarrow{iLine}{triangle}----
\newarrow{ILine}{blacktriangle}----
\def\rH{\Rightarrow}
\def\lH{\Leftarrow}
\def\uH{\Uparrow}
\def\dH{\Downarrow}
\newarrowhead{=>}{\rH}{\lH}{\dH}{\uH}
\newarrow{DoubleLine}=====%
\newarrow{DoubleTo}===={=>}%
\newarrow{DoubleTO}===={curlyvee}%
\newarrow{Equal}=====%
\newarrow{DoublesTo}===={triangle}
\newarrow{DoubleSTo}===={blacktriangle}
\newarrow{DoubleiTo}{triangle}==={=>}
\newarrow{DoubleITo}{blacktriangle}==={=>}
\newarrow{DotTo}....{->}
\newarrow{DotsTo}....{triangle}
\newarrow{DotSTo}....{blacktriangle}
\newarrow{DotiTo}{triangle}...{->}
\newarrow{DotITo}{blacktriangle}...{->}
\newarrow{DashTO}{}{dash}{}{dash}{curlyvee}
\newarrow{DashsTo}{dash}{dash}{dash}{dash}{triangle}
\newarrow{DashSTo}{dash}{dash}{dash}{dash}{blacktriangle}
\newarrow{DashiTo}{triangle}{dash}{}{dash}{>}
\newarrow{DashITo}{blacktriangle}{dash}{dash}{dash}{>}
\newarrow{DashLine}{dash}{dash}{dash}{dash}{dash}
\newarrow{DashiLine}{triangle}{dash}{dash}{dash}{dash}
\newarrow{DashILine}{blacktriangle}{dash}{dash}{dash}{dash}
\newarrow{DotLine}.....

\def\rsum{\raisebox{2pt}{$\:\scriptscriptstyle+\:$}}
\def\lsum{\raisebox{2pt}{$\:\scriptscriptstyle+\:$}}
\def\dsum{\vphantom b\scriptscriptstyle+}
\def\usum{\vphantom b\scriptscriptstyle+}

\newarrowfiller{+}\rsum\lsum\dsum\usum
\newarrow{SumSTo}{+}{+}{+}{+}{blacktriangle}
\newarrow{SumSTo}{+}{+}{+}{+}{blacktriangle}
\newarrow{SumLine}{+}{+}{+}{+}{+}
\newarrow{SumTo}{+}{+}{+}{+}{->}
\newarrow{SumiTo}{triangle}{+}{+}{+}{->}

\def\rrnd{\raisebox{1.5pt}{$\:\scriptscriptstyle{\circ}\:$}}
\def\lrnd{\raisebox{1.5pt}{$\:\scriptscriptstyle{\circ}\:$}}
\def\drnd{\vphantom b\scriptscriptstyle{\circ}}
\def\urnd{\vphantom b\scriptscriptstyle{\circ}}

\newarrowfiller{rnd}{\rrnd}{\lrnd}{\drnd}{\urnd}
\newarrow{RndSTo}{rnd}{rnd}{rnd}{rnd}{blacktriangle}
\newarrow{RndLine}{rnd}{rnd}{rnd}{rnd}{rnd}

\def\rbul{\raisebox{1.5pt}{$\:\scriptscriptstyle{\bullet}\:$}}
\def\lbul{\raisebox{1.5pt}{$\:\scriptscriptstyle{\bullet}\:$}}
\def\dbul{\vphantom b\scriptscriptstyle{\bullet}}
\def\ubul{\vphantom b\scriptscriptstyle{\bullet}}

\newarrowfiller{bul}{\rbul}{\lbul}{\dbul}{\ubul}
\newarrow{BulSTo}{bul}{bul}{bul}{bul}{blacktriangle}
\newarrow{BulLine}{bul}{bul}{bul}{bul}{bul}

\def\imm{\text{\scriptsize{\textcircled{i}}}}
\def\rimm{\imm}
\def\limm{\imm}
\def\uimm{\imm}
\def\dimm{\imm}
\newarrowmiddle{imm}{\rimm}{\limm}{\uimm}{\dimm}
\newarrow{imTo}--{imm}-{->}
\newarrow{plfTo}{LaTeX}-{imm}-{->}
\newarrow{fidTo}{curlyvee}-{imm}-{->}
\newarrow{istTo}{boldlittlevee}-{imm}-{->}

\newarrow{qimTo}{o}---{->}
\def\eq{\text{\textcircled{e}}}
\def\etl{\text{\scriptsize{\textcircled{e}}}}
\def\retl{\etl}
\def\letl{\etl}
\def\uetl{\etl}
\def\detl{\etl}
\newarrowmiddle{etl}{\retl}{\letl}{\uetl}{\detl}
\newarrow{etaTo}--{etl}-{->}
\newarrow{etasTo}--{etl}-{triangle}
\newarrow{ploTo}{triangle}-{etl}-{->}
\newarrow{pofTo}{littleblack}-{etl}-{->}

\def\sub{\text{\scriptsize{\textcircled{s}}}}
\def\rsub{\sub}
\def\lsub{\sub}
\def\usub{\sub}
\def\dsub{\sub}
\newarrowmiddle{sub}{\rsub}{\lsub}{\usub}{\dsub}
\newarrow{subTo}--{sub}-{->}

\def\rph{\phantom{-}}
\def\lph{\phantom{-}}
\def\uph{\phantom{-}}
\def\dph{\phantom{-}}

\newarrowfiller{PH}{\rph}{\lph}{\uph}{\dph}

\newarrow{PHL}{PH}{PH}{PH}{PH}{PH}

\def\sq{\square\ }

\def\spb{\hat{\square}}
\def\gpb{\ov{\diamond}{\square}}
\def\ipb{\check{\square}}

\def\pb{\text{\raisebox{1.6pt}{{\fboxsep0.7pt\framebox{\tiny{\sffamily pb}}}}}}

\def\ovra{\overrightarrow}
\def\ov{\overset}
\def\un{\underset}
\def\lv{\textit{\raisebox{.3ex}{{\scriptsize{v\,}}}}}
\def\li{\textit{\raisebox{.3ex}{{\scriptsize{i\,}}}}}
\def\ls{\textit{\raisebox{.3ex}{{\scriptsize{s\,}}}}}

\def\vQD{\ls\negthickspace\raisebox{.5ex}{${\un{Q}{\bigtriangledown}}$}\sD}

\def\emphc{\mathcal}
\def\emphf{\mathfrak}

\def\h{\quad}
\def\H{\qquad}
\def\HH{\H\H}

\begin{document}

\title[groupoid geometries]{In Ehresmann's footsteps:\\
from group geometries\\to groupoid geometries}

\author{Jean PRADINES}
\address{26, rue Alexandre Ducos, F31500 Toulouse, France}
\email{jpradines@wanadoo.fr} 
\keywords{Principal bundles, Lie group and Lie groupoid actions, \textsc{Ehresmann}}

\subjclass{58H05}

\thanks{The author expresses his grateful thanks to the head and the staff of Banach Centre who welcomed the meeting in such a prestigious and comfortable setting. He also thanks warmly all the Polish colleagues who spared no effort in contributing to the superb organization of this Conference, dedicated to the Mathematical Legacy of Charles \textsc{Ehresmann}, on the occasion of the hundredth anniversary of his birthday }

\thanks{The author is indebted to Paul Taylor, whose package \emph{diagrams} is used throughout this text for drawing a variety of diagrams and arrows.}

\date{28/09/2006}

\dedicatory{%
\begin{flushright}
Dedicated to the memory of Charles \textsc{Ehresmann}\\
on the occasion of the hundredth anniversary of his birth\\
as a token of respect and admiration.
\end{flushright}
}%

\begin{abstract}

The \emph{geometric understanding} of \emph{Cartan connections} led Charles \textsc{Ehres\-mann} from the Erlangen program of (abstract) \emph{transformation groups} to the enlarged program of \emph{Lie groupoid actions}, via the basic concept of \emph{structural groupoid} acting through the fibres of a (smooth) \emph{principal fibre bundle} or of its \emph{associated} bundles, and the basic examples stemming from the \emph{manifold of jets} (fibred by its source or target projections).

We show that the remarkable relation arising between the actions of the \emph{structural group} and the \emph{structural groupoid} (which are mutually determined by one another and commuting) may be viewed as a very special (\emph{unsymmetrical}!) instance of a general \emph{fully symmetric} notion of ``\emph{conjugation between principal actions}'' and between ``\emph{associated actions}'', encapsulated in a nice ``\emph{butterfly diagram}''. In this prospect, the role of the \emph{local triviality} looks more incidental, and may be withdrawn, allowing to encompass and bring together much more general situations.

We describe various examples illustrating the ubiquity of this concept in Differential Geometry, and the way it unifies miscellaneous aspects of fibre bundles and foliations.

We also suggest some tracks (to be developed more extensively elsewhere) for a more efficient implementation of the basic principle presently known as ``\emph{internalization}'', pioneered by \textsc{Ehresmann} in his very general theory of ``\emph{structured'' categories and functors}, towards the more special but very rich and far-reaching study of the above-mentioned \emph{Lie groupoid actions}.

Still now, due to \emph{misleading and conflicting terminologies}, the latter concept seems too often neglected (and sometimes misunderstood) by too many geometers, and has long been generally ignored or despised by most ``pure categorists'', though it will be presented here as one of the gems of \textsc{Ehresmann}'s legacy.

\end{abstract}
\maketitle

\setcounter{tocdepth}{2}
\tableofcontents

\section{General introduction.}
\label{genintro}

\subsection{Structural groupoid.}

In this lecture we intend to focus on one of the turning points in Charles \textsc{Ehres\-mann}'s work : we mean, when dealing with the theory of principal bundles, the introduction, parallel to that of the \emph{structural group}, of the \emph{structural groupoid}\footnote%
{Note that nowadays these are often known too as \emph{gauge group} and \emph{gauge groupoid} \cite{MK}, rather conveniently though ambiguously, since these terms possess also, according to the authors, quite different (and conflicting) meanings (especially in Physics, in gauge theories \cite{BGT, NGTP}), which are related to the (local or global) \emph{automorphisms} of the structure, or to its local \emph{trivializations}. For this reason we shall avoid these too ambiguous terms, keeping here \textsc{Ehresmann}'s original terminology.}. %
The actions, on the bundle, of the group (global action, transitive along the fibres) and the groupoid (partial action from fibre to fibre, transitive) are commuting.
Together with the study of the category of jets (and more specifically the groupoid of invertible jets), this led him to the basic idea (rediscovered much later by various authors) of considering (\emph{small}) \emph{internal groupoids} (more generally categories) in various (\emph{large}) ``\emph{structuring}'' \emph{categories}, originally and basically the category of [morphisms between] manifolds. In particular this yields an inductive definition of multiple (smooth) groupoids, which play an important, though too often hidden, role. 

This unifying concept involves a very far-reaching \emph{intertwining between algebra and geometry}.

\textsc{Ehresmann} also stressed the interest of choosing the source and target maps and\,/\,or the so-called ``\emph{anchor map}'' \cite{MK} (i.e. target-source map, which we prefer to call, more suggestively, ``\emph{transitor}'') in suitable \emph{subclasses} of the ``structuring'' category. We shall show it is convenient to impose to these subclasses some suitable \emph{stability properties} of a very general nature, and propose various choices unifying various theories. We point out that the stability properties we need are very easily satisfied in the topos setting, but \emph{one gets much wider ranging theories} (especially when aiming at applications to Differential Geometry) \emph{when not demanding the structuring category to be a topos}, which would require much more special properties, \emph{not satisfied} in \textbf{Dif} for instance.

\subsection{Interpretations of the unsymmetrical situation.}\label{inter}

We note that the concepts of structural group and groupoid (denoted just below by $H$ and $G$) have a strongly intuitive (geometrical or ``physical'') interpretation. When interpreting the elements of the \emph{principal}\footnote
{Actually these interpretations will become perhaps clearer and more flexible when dealing with \emph{associated} bundles, to be described below (\ref{ass}).}
 \emph{bundle} $P$ (regarded as frames) as ``measures of some events'' or ``observations'', the fibres, which are also the elements of the base $B$ or the objects of the groupoid $G$, may be viewed as ``observers'' (or disjoint classes of observations), while the base $B$ itself makes up a class of interconnected observers, defining a ``\emph{coherent point of view}''. The structure of the principal bundle $P$ is then interpreted as allowing a ``\emph{change in point of view}''.

In that respect, the unique object of the group would play the role of an abstract (privileged!) ``\emph{universal observer}'', and the group $H$ itself of an ``\emph{absolute}'' gauge reference. In other contexts, such as Galois theory or renormalization theories, this group is often interpreted as an ``\emph{ambiguity}'' group.

Alternatively, from this absolute viewpoint, the base space $B$ may be interpreted as an \emph{abstract} absolute universe (arising from ``experimental'' coherence of observations, and picturing some emerging ``objective reality''), liable to evolve globally (or alternatively to be observed in various states) under the global action of $H$.

On the other hand the groupoid $G$ (with base $B$) allows \emph{direct comparisons} (here uniquely determined by a pair of observations) between the fibres of the various (\emph{relative}) observers, without using the medium of the \emph{absolute} observer. However the latter (relative) comparisons actually involve the \emph{isotropy groups}\footnote%
{This is \textsc{Ehresmann}'s terminology, which is perfectly convenient. They are called ``coherence groups'' in \textsc{van Est}'s papers \cite{vE}; they make up what is called rather ambiguously and improperly the inner groupoid in \cite{MK}.}.

A (global) \emph{trivialization}, when it happens to exist, of $G$, or equivalently a (global) \emph{section} of the principal fibration $p:P\rightarrow B$, yields an absolute comparison, hence a ``classical'' universe, for which one would have $P=B\times H$, while in general the ``objective universe'' $B$ is constructed as just a \emph{quotient} of the ``space of observations'' $P$. However it should be noted that the \emph{choice} of such a section ( even when it does \emph{exist} globally) is \emph{not} determined by the sole data of the principal fibration, and constitutes an extra datum: it is parametrized by the group of maps from $B$ to $H$, which is \emph{not} in general a finite dimensional Lie group, unlike $H$ itself.

An \emph{infinitesimal connection} may be described in geometrical terms as a (generally non-integrable) field of \emph{infinitesimal jets of trivializations}, or equivalently of local \emph{sections}\footnote{
\textsc{Ehresmann} would favour the complementary interpretation as jets of \emph{retractions} on target fibres ($\beta$-fibres). However, favouring the sectional interpretation, we showed in \cite{P67, P68} that the basic geometric features of infinitesimal, and more specifically linear, connections, such as curvature, torsion, covariant differentiation, parallel transport, exponential map, are not only preserved and generalized, but better understood, when \emph{replacing the banal groupoid} $B\times B$ \emph{of pairs of points of the base} $B$ (\emph{ubiquitous though hidden}!) by a (\emph{visible}, and \emph{possibly non-locally trivial}) Lie groupoid $G$, and the transitor or anchor map $\tau_G:G\rightarrow B\times B$ (a surmersion iff $G$ is locally trivial and transitive) by a base-preserving \emph{surmersive} morphism $f: G'\rightarrow G$, then \emph{specializing} to the case when this morphism is the canonical projection $s_G:SG\rightarrow G$ of \textsc{Ehresmann}'s infinitesimal prolongation $SG$ (jets of bisections) (see Rem.~ \ref{prol} below). See Rem.~ \ref{loctriv} below for examples of non locally trivial groupoids in which connections in this sense might be useful.}
 of the ``transitor'' $\tau_{G}=(\beta_G,\alpha_G):G\rightarrow B\times B$, or of $p:P\rightarrow B$. This yields absolute comparisons only along some given path, depending upon the path.

Now the transition between the two previous ``points of view'' (absolute and relative) is realized by means of the projections of the bundle $P$ onto the bases of the group $H$ and of the groupoid $G$ (the latter being $B$, and the former a \emph{singleton}!), and by the basic fact that the group and the groupoid induce on $P$ (by pulling back along these projections) the ``same'' (here trivial and trivialized) groupoid $K=H\times(P\times P)$ with base $P$ (more precisely isomorphic groupoids). 

From a purely algebraic point of view, this describes a \emph{Morita equivalence} between the structural group $H$ and the structural groupoid $G$. But in \textsc{Ehresmann}'s internal setting, this may acquire a lot of different more precise ``\emph{smooth}'' meanings, depending on the above-mentioned various choices. It is also a special case (isomorphism) of the ``\emph{generalized morphisms}'' in the sense of \textsc{Skandalis} and
 \textsc{Haefliger}\footnote
{At least in the purely topological setting, i.e. when forgetting the smoothness. It should be noted that in this framework the local triviality conditions seem to bring about some difficulties in the general definition of the composition of these generalized morphisms; these difficulties are more easily mastered in the smooth setting, via the diagrammatic description we propose here, which relies on Godement theorem.}
 (cf.~ \cite{H}), which indeed may be interpreted as a very geometric description of the arrows of the \emph{category of fractions} \cite{GZ}\footnote
{It is worth noting that such a geometric description is specific of the present situation, and quite different from the so-called ``calculus of fractions'', available under general conditions \cite{GZ} which are \emph{not} fulfilled here \cite{P89}.}
 in which the \emph{smooth equivalences} (between smooth groupoids) are formally inverted, i.e. turned into actual \emph{isomorphisms} (cf. \cite{P89, P04}).

This (trivial) groupoid $K$ with base $P$ (which is the core of the structure, and might be called the ``\emph{conjugation groupoid}'') inherits from the two projections a very rich extra structure, which may be recognized as a particular instance of a very special and interesting structure of \emph{smooth} double groupoid (indeed a ``double equivalence relation''). In the \emph{purely algebraic} context, this structure has been described in the literature under various names and equivalent ways: rule of three, affinoid (\textsc{A.~ Weinstein}), pregroupoid (A.~ \textsc{Kock}). See also \cite{P85}.

\subsection{Restoring symmetry}

Now it turns out that all the essential features of this description remain still valid (and indeed are better understood) when replacing the starting structural \emph{group} $H$ by a (\emph{possibly non-transitive}) \emph{groupoid} (with base $C$, no longer a singleton!), and then the situation becomes \emph{perfectly symmetrical}, or reversible (up to isomorphisms)\footnote
{This concept of a fully symmetrical presentation of ``conjugate principal actions'' had been sketched a very long time ago in an unpublished communication \cite{P77}, and presented in some other oral communications. We seize here the opportunity of giving further developments and explanations.}(section \ref{sec:pact}).
Note that now the core $K$ is itself generally \emph{no longer trivial}. In the purely algebraic (hence also topos) setting, this symmetric situation has been described in several equivalent ways by \textsc{A. Kock}\footnote
{See his contribution to the present Conference.}
, in the language of torsors (which we call here \emph{principal actors}) and bitorsors.

The symmetry of the situation will be fully apparent when drawing what we call the ``\emph{butterfly diagram}'' (see below \ref{cdpm}).

In the ``physical'' interpretation, we can say that we have \emph{two classes of observers} (or two correlated ``coherent points of view'' or ``descriptions''), \emph{now placed on equal footings} in the absence of a privileged ``absolute'' point of view, and a \emph{transition between two ``conjugate points of view''}, determined by the two projections: 
	\[B\stackrel{p}{\lTO}P\stackrel{q}{\rTO}C\ .
\]
The \emph{orbits} of the action of one of the two conjugate groupoids on $P$ become the \emph{objects} of the \emph{other} one.

In spite of the ``coherence'' of observations, encapsulated in the groupoid $K$ (induced on $P$ by both $G$ and $H$ along $p$ and $q$), it is no longer possible to construct a privileged ``objective reality'' in the restricted sense described above: the framework of ``orbital structures'' \cite{P89, P04} is then better adapted.

We also note that it is no longer demanded the groupoids to be transitive even when they are (topologically!) connected (this means that possibly certain pairs of observers cannot compare their observations, though belonging to the same class), nor also the rank of the anchor map to be constant (then the isotropy groups are allowed to vary, and this might be thought as modeling \emph{symmetry breaks} or \emph{phase changes} \cite{NGTP}).

Indeed, though, from a purely \emph{algebraic} point of view, any groupoid is the disjoint sum of its transitive components or orbits (each one being individually equivalent to a group\footnote
{There lies the origin of the still widespread disinterest for groupoids among many ``categorists'', who are forgetting than (even from a purely algebraic point of view), classifying objects does not imply classifying morphisms!}),
 the \emph{topology} and the \emph{smoothness} involve a very subtle and deep interconnection between these components, which in our opinion, would deserve a better understanding and exploration from a geometrical (and physical) point of view (see \cite{P03}).

For instance one of the simplest invariants arising from this interconnection is the ``\emph{fundamental group}'' of the orbit space, constructed along different and independent ways by \textsc{Haefliger} \cite{H} and \textsc{van Est} \cite{vE} (for a universal characterization of this group, see \cite{PA} and \cite{P89}). In this prospect the classical fundamental group would just arise from the interconnection of the topology of a space with the banal groupoid structure on its set of pairs); this has been advocated by \textsc{R.~ Brown} for a very long time (see \cite{BRO}).

We shall give in \ref{exactiv} various geometric examples of this unifying (purely diagrammatic) situation (and more general ones), which encompasses (among others) the construction of associated bundles, the realization of a non abelian cocycle, including Haefliger cocycles, the construction of the holonomy groupoids of foliations, and the Palais globalization of a local action law \cite{Pal}.

Due to the necessity of making choices, we shall not describe here the corresponding infinitesimal situation, which of course was \textsc{Ehresmann}'s main motivation for introducing the structural groupoid, in order to understand the meaning of ``infinitesimally connecting'' the fibres. This will be tackled by other lecturers at the present Conference (notably \cite{MA} and \cite{Lib}).

\subsection{Prospects.}

Our conclusion will be that the transition, from the Kleinian conception of Geometry as the study of transformation groups or \emph{group actions}, to the Ehresmannian enlarged point of view, consisting in considering \emph{Lie groupoid actions}, involves a conceptual revolution which parallels (and possibly might model) the physical revolution from the classical concept of an absolute universe to the modern visions about our physical space. Better qualified observations may be found in \textsc{R. Hermann}'s comments \cite[(VII)]{ECO}.

It might well be that this revolution arising from Ehresmann's prophetic visions is not enough well understood presently, and is liable to come out into unexpected developments.

\section{From discrete group actions to structured groupoid actions\\or from \textsc{Erlangen} program to \textsc{Ehresmann}'s program}
\subsection{From transformation groups to Lie group actions and transformation pseudogroups}
\label{tg}

\subsubsection{Lie group actions}
\label{compatibility}
As is well known, the modern conception of Geometry was formulated in a very \emph{global} and essentially \emph{algebraic} way by \textsc{F.~ Klein} in his Erlangen Program (1872) as the study of invariants of some specific \emph{transformation group} acting ( often transitively) on an \emph{abstract set}.

Then, with the emergence of the general concepts of \emph{structures} and \emph{morphisms}, it was progressively realized that it is convenient to endow the set and the group itself with some specific extra structures of the same species, which are intimately related by the action. The latter should be ``\emph{compatible}'' with these structures, in a sense to be made precise below.

{\footnotesize
\begin{rem}
\label{remextra}
As we are aiming basically at differentiable structures, it will be convenient for our purpose to consider this structure on the group as an \emph{extra datum}, and we shall not argue here about the existence of those structures on such a group (such as compact-open topologies for instance) which may arise in an intrinsic way from the very structure of the space on which it acts. That sort of discussions involves \emph{specific} properties of the structuring category (such as the property of being or not ``cartesian closed'' \cite{ML}) and leads to important but technical and specific problems which are quite different in the topological and in the smooth cases for instance, and should not be mingled with the very \emph{general} and \emph{simple} notion of (``structured'') action we are considering and developing here. This attitude allows to completely short-circuit the difficulties and controversies which arose at the very birth of fibre bundles about the role of the structural group (see \textsc{M.~ Zisman}'s comments on that point in \cite[~(II)]{ECO}).
\end{rem}
}%

More precisely \emph{the action has to be described in terms of morphisms of the ``structuring category''} (which is often \textbf{Dif} or \textbf{Top}), \emph{as well as the group law}. Such a formulation remains \emph{algebraic}, in spite of the added extra structures, \emph{via the categorical language}. This is one of the simplest instances of ``\emph{internalization'' of an algebraic structure}, which \textsc{Ehresmann} extended later to groupoids, categories and more general action laws.

Basic examples of such added structures are topological or smooth (i.e. differentiable manifold) structures. Then morphisms are continuous or smooth maps, and this means that one is considering \emph{topological} or \emph{Lie groups} (in the modern global sense), \emph{acting} \emph{continuously} or \emph{smoothly} on a topological space or a manifold.

\subsubsection{Pseudogroups vs groupoids.}
\label{pseudogroup}
Such a conception is highly \emph{global}, hence rigid.

As a differential geometer, and with various motivating examples in mind, Ch. \textsc{Ehres\-mann} realized the need for what he called \emph{local} structures, for which he was led to introduce very general definitions. Then the \emph{localization} of the concept of transformation group brought him to produce a precise definition of the general notion of transformation \emph{pseudogroup}. The pseudogroups play a basic role in the general fundamental concept of \emph{atlas}.

We shall not dwell at length upon what will be considered here as an \emph{intermediate step}\footnote%
{In our opinion this is not the most original part of \textsc{Ehres\-mann}'s work, since it might be rewritten essentially in the better known language of (pre)sheaves (see for instance \cite{GOD, T}), which was developed later in parallel ways.},
 for which we refer to the contribution of Andrée \textsc{Charles Ehres\-mann} addressed for the present proceedings, and also to \textsc{J.~ Ben\-abou}'s comments in \cite[~(X)]{ECO}. We just mention (as an anticipation of what follows) that the so-called \emph{complete pseudogroups} may be considered as special cases of \emph{étale} topological (or smooth) \emph{groupoids}; this means that they can be endowed with a topology (or a manifold structure) such that the source (and target) projections are étale (i.e. are local homeo\-morphisms or diffeo\-morphisms). 

We remind also than any pseudogroup may be ``\emph{completed}'' by the general process of ``sheafification'' or ``germification'' giving rise to a sheaf, or equivalently an étale space, starting from a presheaf. This allows to deal only with ``complete'' pseudogroups (in \textsc{Ehresmann}'s sense).

We point out too that a (complete) transformation pseudogroup has to be considered as a localization of a \emph{discrete} (or abstract) transformation group rather than of a general \emph{Lie} group action, and is in this respect \emph{much more special} than the general concept of \emph{structured} groupoid action which is stressed in the sequel\footnote%
{In ``modern'' category theory (see \textsc{A.~ Kock}'s comments \cite[(VII)]{ECO}), this concept might be (in principle!) described in terms of the so-called ``\emph{internal presheaves}'' and ``\emph{discrete (op)fibrations}''. In our opinion these \emph{highly unfortunate terminologies}, though well established in Category Theory and in Algebraic Geometry and Topology (see for instance \cite{GIR, Ja, Jb}), are completely \emph{misleading} as soon as one is interested in \emph{topological} categories or groupoids (introduced much earlier by \textsc{Ehresmann}), and cannot be used by differential geometers (who anyway generally ignore it!), since the main interest of this internalization for the study of fibre bundles is to allow \emph{non-discrete} fibres and \emph{non-étale} groupoids! So that such a terminology makes it \emph{virtually impossible} to understand the drastic change from \emph{pseudogroups} to general \emph{groupoids}. Unfortunately again, as we shall see below, \textsc{Ehres\-mann}'s terminology is equally misleading for different reasons!},
 and \emph{which goes far beyond}.

It is worth noting too that any (\emph{possibly non-complete}) transformation pseudogroup might be viewed also as a groupoid \emph{in a quite different way}, taking as \emph{objects} the \emph{domains} of the transformations (which are generally open sets), and using the \emph{strict composition law} for local transformations. But this latter groupoid does not admit a finite dimensional manifold structure; it will not be considered here.

{\footnotesize 
\begin{rem}
\label{prol}
Though this is out of the scope of the present lecture, let us mention, again to avoid \emph{frequent and persistent confusions between general groupoids and pseudogroups} (see Rem.~ \ref{Lie pseudogroup} below about \emph{Lie pseudogroups}), that (in the opposite direction) the localization of a (topological) \emph{groupoid} leads to the \emph{pseudogroup} of what we think convenient to call the ``\emph{local translations}'' associated to the so-called \emph{local bisections}, i.e. local sections of both source and target maps. Of course, in the case of a \emph{group}, one recovers the (automatically global!) translations, whence our terminology.
(This relation between groupoids and pseudogroups might be stated more precisely in terms of adjointness.)

This pseudogroup was introduced and used by \textsc{Ehresmann} as the \emph{local prolongation} of the groupoid (together with the \emph{infinitesimal} prolongations, when infinitesimal jets are replacing germs, which he calls ``\emph{local jets}'', giving rise to \emph{non-étale} topologies and groupoids), and (more or less partially) rediscovered much later by many authors under various names (for instance \textsc{C.~ Albert}'s ``\emph{glisse\-ments}''; see also \cite{MK}).
\end{rem}
}

\subsection{From pseudogroups and Lie group actions to Lie groupoid actions: a turning point}
\label{lga}
\subsubsection{Fibre bundles.}
Though implicitly present in previous works of \textsc{Ehres\-mann}, the general concept of \emph{fibre bundle} (\emph{espace fibré}) was introduced explicitly by \textsc{Ehres\-mann} in 1941 \cite [/14/]{OC} in collaboration with \textsc{J.~ Feld\-bau} (and independently from \textsc{N.~ Steen\-rod}).

This concept evolved throughout \textsc{Ehres\-mann}'s work, and one can guess it played the role of a mainspring for this work. For more detailed historical comments about this evolution, we refer to the precious comments by \textsc{Ehres\-mann} himself \cite [~/136/]{OC}, and by \textsc{M.~ Zis\-man} and \textsc{A.~ Haef\-liger} in \cite[(II), (IX)]{ECO}, while the applications of this concept in Differential Geometry are described in \textsc{P.~ Liber\-mann}'s comments \cite[~(IV)]{ECO}. This reflects the richness of the concept. For the modern ``classical'' presentation, we refer to \textsc{N.~ Bour\-baki} \cite{VAR}.

In \cite [/14/]{OC} there is no topology on the structural group, while it becomes a topological group in \cite [/15/ and /16/]{OC} (see comments in Rem.~ \ref{remextra} about this point), and the problem of the reduction of the structural group arises in \cite [~/16/]{OC}, while in \cite[~/18/]{OC} the emphasis is put only on the local triviality. The definition by means of a pseudogroup and an atlas is emphasized in \cite[~/20/]{OC}.

We note that in these previous works the \emph{local triviality} conditions are automatically involved in the very definitions and are considered as a basic ingredient of the structure.

\subsubsection{Structural groupoid}
A \emph{turning point} occurs in \cite [~/28/]{OC} with the appearance of the \emph{structural groupoid} suggestively denoted by $PP^{-1}$ (here $P$ is interpreted as a bundle of frames, and the corresponding co-frames belong to $P^{-1}$), viewed as a \emph{quotient} of $P\times P$ (by the diagonal action of the structural group), with its \emph{manifold structure}, and its (partial) transitive \emph{action} on the bundle. 

\emph{The arrows of this structural \emph{ (small!) groupoid} are isomorphisms between the fibres of the bundle, commuting with the ``vertical'' \emph{(global) non-transitive} action of the structural \emph{group} along the fibres}. 

This concept makes explicit the leading geometrical idea underlying the early work of 1938 \cite[~/10/]{OC}, and will lead to the geometric formulation of higher order \emph{connections} in \cite[~/46/]{OC}. Since we are focusing here on the \emph{global} aspects rather than the \emph{infinitesimal} ones, we refer to \textsc{C.-M.~ Marle}'s and \textsc{P. Libermann}'s lectures \cite{MA,Lib} to explain how this concept allowed \textsc{Ehres\-mann} to formalize the geometric ideas implicit in \textsc{E.~ Car\-tan}'s works on infinitesimal connections.

\subsubsection{Generalized fibrations}
\label{terminology}
From now on \textsc{Ehres\-mann} will consider that the essence of a (\emph{generalized}) \emph{fibration} consists in the datum of a \emph{small} groupoid (more generally a category) acting\footnote%
{As emphasized below (\ref{spemor}), this action may be described by a groupoid morphism (a functor). But, as above-mentioned, categorists would call this functor a ``\emph{discrete} (\emph{op})\emph{fibration}''. \emph{This creates an inextricable mess}! 
More generally the more and more widespread intrusion in Category Theory of (pre-existing and well established!) \emph{topological} terms (such as connected, discrete, undiscrete, compact, closed, fibred, foliated, and so on) for naming purely \emph{algebraic} properties, has disastrous outcomes when dealing with \textsc{Ehres\-mann}'s \emph{topological} (and smooth) \emph{categories}, since such terms acquire in this latter framework \emph{two quite independent and conflicting meanings}.
We of course definitely reject here such misuses.}
 on a space, with some extra structures (notably topological or smooth), assumed to be ``\emph{compatible}'', in a suitable sense, with such an algebraic datum. 

(This conception will be discussed in more details in the sequel.)

{\footnotesize 
\begin{rem}
\label{universe}
Indeed we have intendedly added ``\emph{small}'' (this just means that the objects as well as the arrows of the groupoid are belonging to sets\footnote
{We recall that, introducing the language of universes $\emphf{U}$, it is always possible to deal only with $\emphf{U}-$small categories, allowing $\emphf{U}$ to vary when needed (this corresponds to a kind of ``change in scale''). When $\emphf{U}$ is a universe, its \emph{elements} are called \emph{small} sets while its \emph{subsets} are \emph{large} sets. This avoids to speak of ``classes'', for which there is no universally accepted detailed axiomatics available. (Actually it seems that \textsc{Ehres\-mann} was thinking in a language of universes though most of its written work uses a language of classes).})
 in the previous statement, in order to stress that a (small) groupoid has to be thought as a fairly mild (though fundamental) generalization of a group rather than a special case of (large) category. On the contrary \textsc{Ehres\-mann}, who was at the same time and on a parallel way deeply involved in formulating a very general theory of ``\emph{structures}'' (in which the properties of \emph{forgetful functors} between ``\emph{large}'' categories play a basic role\footnote
{For instance \textsc{Bourbaki}'s axiom of ``transport de structure'' may be interpreted by saying that the (large) groupoid of bijections between sets acts on the class of sets endowed with a structure.}),
 attached the utmost importance to connecting the formal aspects of the theory of \emph{fibred spaces} and of the emergent theory of ``\emph{species of structures}'', adopting and extending unfortunately the language of the latter to the former\footnote
{Personally we think the contrary would have been much preferable : for us it seems more comfortable from a psychological point of view to regard a \emph{functor} as a special case of the (informal) general concept of \emph{morphism between certain kinds of algebraic structures} (where the term ``morphism'' suggests some suitable kind of ``compatibility'' with these structures), than to consider a morphism of groups (or small groupoids) as a special instance of a functor; but \textsc{Ehres\-mann} felt the opposite (as many pure ``categorists'' do).}.

We are convinced that this had historically disastrous consequences, and largely explains the misunderstanding which occurred with a large part of the mathematical community.

In actual fact we think that, in spite of certain common formal aspects, the most interesting problems \emph{are seldom the same} in the two theories, so that a too strictly parallel and common treatment results in obscuring and weakening rather than enriching both of them. (For instance the search for universal arrows is uninteresting when all arrows are invertible, while, on the other hand, the groupoids have essentially the same specificity among categories as groups among monoids, and also equivalence relations among preorders, and nobody would claim that group theory is nothing else than a specialization of the more general theory of monoids, and that equivalence relations are just a trivial special case of the more general preorder relations.)

Moreover at that age the main stream of ``categorists'' were not prepared to the concept of ``internalization'' (which was rediscovered and popularized much later), being focused on the \textbf{Hom} point of view, while on the other hand many ``working mathematicians'' felt some distrust of what was often considered as ``abstract nonsense''; so that the basic and simple concept of structural groupoid, with its action on the bundle (and its smooth extra structure), disguised as a (more or less degenerate) instance of very abstract and general ``species of structure'', remained largely ignored, in spite of its deep geometrical significance and usefulness, and is still far from being as widely understood and used as it should be in our opinion.
\end{rem}
}%

This point of view is clearly emphasized in \cite[/37/ and /44/]{OC} (and later in \cite[/111/]{OC}), and related to this basic example of fibre bundle consisting of the manifold of \emph{infinitesimal jets} $J^r(V_n,V_m)$, fibred by the source (resp. target) projection on $V_n$ (resp. $V_m)$ and the action of the groupoid of invertible jets on the source (resp. target) side.

\subsubsection{Ehres\-mann's program}

One can realize that \emph{most of} Ch.~\textsc{Ehres\-mann}'s \emph{subsequent work} was devoted to finding the most general definitions and theorems arising from the previous \emph{basic situation}. This will be referred to in the sequel as ``\emph{Ehres\-mann's program}'', and the above presentation stresses how clearly it parallelizes and enlarges \emph{Erlangen program} in order to include (notably) the understanding of \emph{infinitesimal connections}, following closely the early geometric ideas of \cite[/10/]{OC}.

\subsubsection{Topological and differentiable categories and groupoids}
The precise conditions of ``\emph{compatibility}'' in the above formulation when the extra structure is a topology or a manifold structure were described in 1959 \cite[/50/]{OC}, and the basic differentiable case was emphasized in \cite[/103/]{OC} in 1967.

These will be discussed more thorougly in the sequel. Roughly speaking they can be conveyed by saying that the algebraic groupoid structure as well as its action have to be described by \emph{morphisms} of topological spaces or manifolds. 

In particular this leads to the basic concepts of \emph{topological} and \emph{Lie groupoids}.

{\footnotesize 
\begin{rem}\label{Lie pseudogroup}
\textsc{Ehres\-mann}'s \emph{differentiable} (here frequently called also \emph{smooth}) \emph{groupoids} are nowadays known as \emph{Lie groupoids}, following \textsc{A. Weinstein} and \textsc{P. Dazord} \cite{CDW} (1987). A detailed up-to-date presentation by \textsc{K. Mackenzie} (including this new terminology, which we use here throughout) became recently available in textbook form in \cite{MK}. 

However recall that, till 1987, this latter term was kept by various authors for the special case of \emph{locally trivial} (or equivalently \emph{locally transitive}) \emph{transitive} groupoids (which are essentially the structural groupoids of locally trivial principal bundles, as recalled below). This change in terminology seems to be widely accepted nowadays, but of course it may be a source of confusion when reading ancient papers (notably \textsc{Ngo van Que}, \textsc{A. Kumpera}, \emph{et alii}, including the ancient version of \cite{MK}, published precisely in 1987).

The present terminology is now \emph{coherent} with the terminology we had introduced in 1967 for the corresponding infinitesimal objects, in the general (possibly non locally transitive) case, under the name of \emph{Lie algebroids} \cite{P67} (while such infinitesimal objects had been considered previously by \textsc{C.~ Ehres\-mann} and \textsc{P. ~Liber\-mann}, more or less implicitly, only in the original locally trivial case).

A basic example of (non-étale!) such Lie groupoids is the groupoid of r-jets of local diffeomorphisms of a manifold $B$, which \textsc{Ehres\-mann} used to denote by $\Pi^r(B)$. He used them for a \emph{geometrical definition of certain partial differential equations systems}, viewed as smooth subgroupoids $G$ of some $\Pi^r(B)$. Then the (local) \emph{solutions} of such a system make up a complete \emph{pseudogroup}, which \textsc{Ehres\-mann} named a \emph{Lie pseudogroup}. As explained in Rem. \ref{prol}, the \emph{germs} of solutions, with the \emph{étale} topology, may be regarded as making up a (\emph{very special, since étale}!) Lie groupoid $S$.

Unfortunately \textsc{B. Malgrange} in 2001 \cite{Mal} introduced this same generic term of \emph{Lie groupoid} (or $\gr{D}$-groupoid) only with this very restricted and special meaning of \emph{Lie pseudogroup} (germs of solutions), and only in the \emph{analytic case} (more precisely he considers ``analytic manifolds'' which are \emph{possibly non-smooth}, i.e. admitting some \emph{singularities}). Due to the importance of this work, this terminology is becoming widespread in this area.

This may be especially confusing and misleading in the very important and interesting special case of \emph{finite type Lie pseudogroups}, studied by \textsc{Ehres\-mann} in 1958 (cf \cite[/48/]{OC}, where a very general and powerful version of \textsc{Palais} theorem was stated and proved in a very elegant way, as an application of this theory, using \textsc{E. Cartan}'s closed subgroup theorem: see \cite{Lib}): this names the case when the defining system of equations $G$ may be chosen (possibly after ``prolongation'') in such a way that there exists a unique germ of solution through each point of~ $G$. Then $G$ (equations) and $S$ (germs of solutions) have the \emph{same} underlying set (and \emph{algebraic} structure of groupoid), but \emph{two quite different structures of Lie groupoids} in the general sense, the former \emph{non-étale} (which encapsulates a rich information), the latter finer and \emph{étale} (whith the \emph{germ topology} of the Lie \emph{pseudogroup}).
\end{rem}
}%

\begin{rem}\label{loctriv}
It should be noted that the \emph{local triviality} is then reflected by a special property (namely \emph{transitivity} and \emph{local transitivity}) of the structural groupoid and may be easily withdrawn (together with the local description by means of atlases) in order to encompass more general situations in which the main basic facts remain still valid, and which may cover many other interesting applications. \textsc{Ehres\-mann} pointed out explicitly that this \emph{allows the fibres to vary} (we shall come back to this point below), in contrast with the locally trivial case.

It may be proved \cite{P85, P03} that on any Lie groupoid there is defined (\emph{in a universal way}) a finer manifold structure (determining a \emph{\St\ foliation} \cite{ST}) which turns it into a locally trivial Lie groupoid (the latter is then the sum of its locally trivial transitive components)\footnote{
In the case of an \emph{étale} Lie groupoid (in particular a \emph{pseudogroup}) this new structure is always uninteresting since \emph{discrete}!}. 
This allows to say that \textsc{Ehres\-mann}'s \emph{concept of differentiable groupoids gives a precise meaning to the concept of (generalized) principal bundle in which the structural group may ``vary smoothly''}, which, we guess, might and should play an important role in Physics for gauge theories (see \textsc{R.~ Hermann}'s comments \cite[(VII)]{ECO}) in the interpretation of symmetry breaks and phase changes \cite{BGT, NGTP}. 

Note that \emph{the \emph{(``vertical'')} variation of the \emph{isotropy group} dimension is counterbalanced by the opposite \emph{(``horizontal'')} variation of the \emph{orbit} dimension}. This may be viewed as a refined process of smoothly blowing up a singularity, which results in attaching to it an emerging group ``measuring'' it.

\emph{We think that the precise general mechanisms of this remarkable ``compensation'' are not presently well enough understood, and we are convinced that it is of the utmost importance to clear them up} (a special case is known as \textsc{Wein\-stein}'s conjecture). Very interesting and typical \emph{examples} of this situation are known:
\begin{itemize}
	\item A. \textsc{Connes}'s \emph{tangent groupoid} \cite{NCG};
	\item the groupoids associated to \emph{manifolds with edges and corners} (see \textsc{A. ~Wein\-stein}'s survey \cite{W} and \textsc{B.~ Monthubert}'s lecture \cite{MH});
	\item the \emph{holonomy groupoid} of a regular foliation (introduced by \textsc{Ehres\-mann} in the more general context of topological foliations: see our comments in \cite{P03}), which describes the smooth variation of holonomy galoisian coverings of the leaves; a complete characterization of such groupoids under the name of ``\emph{graphoids}'' (because they behave locally like pure graphs, in spite of their non-trivial isotropy) is given in \cite{P84}, among those groupoids for which the transitor (or anchor map) is an immersion;
	\item more generally all the previous examples are special instances of ``convectors'' (or ``\emph{quasi-graph\-oids}'') which I introduced~in 1985 (in collaboration with \textsc{B.~ Bi\-gon\-net} \cite{PB, P85}), when studying the holonomy of certain \St\ foliations, to designate those groupoids for which the pseudogroup of bisections or local translations (Rem. \ref{prol}) is \emph{faithfully} represented (via the anchor map or ``transitor'') by local diffeomorphisms \emph{of the base}, thus partaking properties of both graphs and pseudo\-groups.
	\item \footnotesize{%
A typical (though very elementary) special \emph{example} is the \emph{smooth} (and even complex analytic) non-transitive groupoid $G$ describing the (non-free) action of the multiplicative group $\bbC^\ast$ on $\bbC$, for which $0$ is an isolated fixed point. The isotropy group of $0$ is $\bbC^\ast$. This groupoid $G$ defines an atlas for an \emph{orbital structure} in the precise sense introduced in \cite{P89, P04}. This orbital structure is not trivial: it has just two orbits, but the orbit $\{0\}$ has $\bbC^\ast$ for its isotropy group\footnote{%
Starting from $\bbR$, one can define analogously an orbital structure with three orbits: $0,+,-$.%
}.%

Though, as a smooth \emph{groupoid}, $G$ is far from being trivial nor locally trivial (the \St\ foliation alluded to in Rem. \ref{loctriv} has just two leaves, one open and one closed), its underlying \emph{manifold} structure is just the product $\bbC^{\ast}\times\bbC$, and the source map $\alpha_G$, which is the canonical projection, determines a \emph{trivial} principal fibration with structural group $\bbC^\ast$ (as well as the target map $\beta_G$).

In spite of the lack of local triviality for $G$, its canonical \textsc{Ehresmann}'s prolongation $SG$ (jets of bisections: see footnote in \ref{inter}) is canonically isomorphic to the product of $G$ by a Lie group, and $s_G:SG\rightarrow G$ admits a canonical section, which is a connection in the generalized sense alluded to in this footnote\footnote{%
The connections recently proposed by \textsc{A. Connes} for describing the structure of space-time seem to be of this type, with one more complex dimension.
}.

When dropping the origin, $G$ induces on the punctured complex plane the banal groupoid $\bbC^\ast\times\bbC^\ast$.

Adding a second fixed point $\infty$, this action may be extended to the Riemannian sphere $\dot{\bbC}= \mathbf{P}_1(\bbC)$, with an action groupoid $\dot{G}$. There are now three orbits.

One can also note than by \emph{blowing up} $0\in\bbC$, one gets the Möbius strip $M= \widehat{\bbC}$. The action of $\bbC^\ast$ extends (still non-freely!) on this blow up. The action groupoid $\widehat{G}$ (with base $M$: now the base has been blown up!) defines an atlas for a \emph{different orbital structure}; the isotropy group of the singular orbit is now~ $\bbZ_2$. One might also blow up $\infty\in\dot{\bbC}$, which yields now the \emph{real} projective plane $\mathbf{P}_2(\bbR)$, or both $0$ and $\infty$, and introduce the analogous groupoids.

Such objects are \emph{implicit} in the elementary theory of functions of one complex variable\footnote{
It would be probably enlightening to make explicit all the groupoids which are implicitly used in the (much more sophisticated) constructions described in \cite{ramar}.
}.

One might also resume the previous constructions when restricting the action of the multiplicative group $\bbC^\ast$ to its subgroup $\mathbf{U}(1)$. The action groupoid is now a subgroupoid $H$ of $G$. Its orbits are concentric circles. These are the leaves of a \St\ foliation, and $H$ is its (generalized) holonomy groupoid in the sense introduced in \cite{P85}.

When blowing up anew the origin in the base, one now recovers the (regular) canonical foliation of the Möbius strip $M$ by means of winding circles, and its (classical) \emph{holonomy groupoid} $\widehat{H}$ with base $M=\widehat{\bbC}$. The orbit space may be interpreted now as the ``folded line''.

A lot of other interesting examples arise when considering various subgroups of the complex projective group and various covers of the previous groupoids.
}%
\end{itemize}
\end{rem}

\subsection{From Lie groupoids to structured groupoids}
\label{sg}
\subsubsection{Structured categories}

\textsc{Ehres\-mann} realized that the previous basic examples, as well as other ones he had encountered (ordered categories, double categories) had to be viewed as special examples of what he called ``\emph{structured}'' categories and groupoids \cite[/63/]{OC}, better known nowadays among categorists as \emph{internal} categories and groupoids (see \cite[(XI) and (XII)]{ECO} for more detailed explanations and comments on this terminology).

Essentially this just means that, when describing the purely \emph{algebraic} situation, consisting of the groupoid \emph{composition law} and the groupoid \emph{action} (as well as the corresponding \emph{axioms}), by means of suitable \emph{diagrams} of maps (verifying certain conditions), the compatibility conditions may be expressed by demanding these maps to be \emph{morphisms of the ``structuring category''} (which may be \textbf{Top} or \textbf{Dif}, or other large categories).

\subsubsection{Dealing with partial composition and action laws}
However the point is that, though this might seem at first sight a straightforward generalization of what was said above in \ref{compatibility} about group actions, it turns out that, in the interesting situations, especially in \textbf{Dif}, very serious technical (and conceptual) difficulties arise from the fact that one has to deal now with composition laws (and action laws) which are only \emph{partially defined}, and the diagrammatic expression of the axioms may be also less obvious. More precisely the data and the axioms involve the consideration of certain pullbacks (i.e. \emph{fibred products} instead of products), and these pullbacks may fail to exist in general as actual objects of the considered category; then certain conditions (for instance of locally constant or maximal rank), intended to ensure the needed \emph{existence of these pullbacks}, look at first sight very specific to special categories, and may seem to have no counterpart in more general categorical contexts. 

It seems that essentially most of \textsc{Ehres\-mann}'s subsequent work was devoted to discover definitions and theorems general and precise enough to encompass all the interesting examples he knew, and peculiarly to bypass or circumvent the special difficulties inherent to the basic category \textbf{Dif}.

\subsection{Trails for pursuing Ehres\-mann's program}
\label{ep}
\subsubsection{Necessity for further developments}

In spite of the huge volume of published papers (see \cite{OC}), \textsc{Ehres\-mann}'s program as defined above is clearly uncompleted, and in the last papers the smooth case seems to be progressively neglected or put aside (or at least outside the general theory). 

Moreover a great part of these publications are very hard to read, due to their terminology and very heavy notations deliberately chosen out of the main mathematical stream, as well as the lack of motivating examples, supposed to be obvious for the reader (see also the comments we made in footnotes of \ref{pseudogroup} and \ref{terminology} and in Rem. \ref{universe}). This is extensively explained and commented in \cite{OC} by Andrée \textsc{Charles Ehres\-mann}, who accomplished remarkable efforts to yield translations, keys and examples in order to make this considerable work more widely accessible. 

However a great part of it seems to remain widely unknown, misunderstood, or poorly understood, though we think it might bear oversimplifications and geometric enlightenments in many areas and are open towards very important developments.

\subsubsection{Objectives}
The aim of the present lecture is twofold:
\begin{itemize}
	\item emphasizing some \emph{leading ideas} of \textsc{Ehres\-mann}'s program which in our opinion have not been widely enough understood presently by the main stream of ``working mathematicians'', particularly differential geometers, though we think they should be a basic asset for the mathematical community;
	\item suggesting some \emph{technical improvements} in the presentation, the terminology
\footnote{As was pointed out above, the very simple concept of \emph{Lie groupoid action} (as a mild but necessary and very important generalization of a \emph{Lie group action}) seems presently completely overshadowed by equally unfortunate terminologies introduced independently (and in conflicting ways) by \textsc{Ehres\-mann} and by ``modern'' categorists.}
, the notations and the formalism, trending to a more efficient implementation and to getting more precise if, at least apparently, slightly less general statements.

Many developments will require further research.
\end{itemize}

\subsubsection{``Right'' generality versus ``maximum'' generality \cite{ML}}

As an overall remark, it may be observed that the search for maximum generality may have two contradictory effects, according to the way it is implemented:

\begin{itemize}
	\item It is certainly highly beneficial when the more general theory actually includes the original one as a special case, and is able moreover:
\begin{itemize}
	\item to bear a new insight and a better understanding of the efficiency of the former theory;
	\item to generate new applications in the same field or in new ones;
	\item to unify seemingly very distant various fields, thereby creating fruitful connections and comparisons.
\end{itemize}
\item It may be harmful if the generality is acquired at the expense of too much weakening the original theory and loosing some of its essential results.
\end{itemize}

\subsubsection{}
The \emph{basic leading ideas} we wish to stress, and which emerge from the above presentation (and appeared progressively throughout the whole of \textsc{Ehres\-mann}'s work), are essentially:
\begin{itemize}
	\item the necessity (particularly in Topology and Differential Geometry) of enlarging the \emph{global} concept of \emph{group actions} into the concept of ``\emph{fibred}'' \emph{actions}, by which we mean here \emph{groupoid actions};
	\item the powerful general principle of ``\emph{internalization}'', which consists in describing the previous \emph{purely algebraic} (and \emph{set-theoretic}) situation (namely group\-oid actions) by means of suitable \emph{diagrams} keeping sense in sufficiently general ``\emph{structuring}'' \emph{categories}, including notably \textbf{Top} and \textbf{Dif}. 
	
	This ``\emph{categorical algebraization}'' trends to supply as much as possible unified expressions of the ``\emph{compatibility}'' between the purely algebraic conditions and the various \emph{extra structures} proceeding from the ``structuring category'', such as \textbf{Top} or \textbf{Dif}. It turns out that \emph{judicious choices} of the alluded diagrams (here indeed lies the core of that kind of approach) may lead to very powerful and illuminating unifying methods.
\end{itemize}

\subsubsection{}
As to the suggested \emph{technical improvements}, developed through the next sections, they bear on the following points (mainly inspired by applications to Differential Geometry):
\begin{itemize}
	\item about the definition of \emph{structured groupoids} and particularly \emph{Lie groupoids}, we propose :

\begin{itemize}
	\item a simplified version of the diagrams used for describing the algebraic structure of groupoids (stressing the specificity of groupoids among categories);
	\item an axiomatization of the properties needed for the subclass to which the source projection is demanded to belong, covering a very wide range of categories used by ``working mathematicians'', and ensuring the existence of all the needed pullbacks (though not all pullbacks).
\end{itemize}
\item as to \emph{groupoid actions}:
\begin{itemize}
	\item a simple diagrammatic description, liable of immediate internalization;
	\item a wide-ranging definition of the concept of \emph{principal} action obtained by internalization of the notion of \emph{free} action; this important and fruitful concept was progressively overshadowed in \textsc{Ehres\-mann}'s work by the most general notion of generalized fibration;
	\item a symmetric description of conjugate and associated actions of the structural group and groupoid on the bundle and on its fibres leading to a unified presentation of a very wide range of basic geometric situations and examples.

\end{itemize}

\end{itemize}

\subsubsection{}

We had earlier opportunities of describing a large part of these technical improvements notably in \cite{P89, P03, P04} \emph{to which we refer for more precisions, explanations and developments}. We shall recall here what is necessary for notations and for understanding how they are tied up to what we name \textsc{Ehres\-mann}'s program.

Then we shall insist on what we call the \emph{conjugation of principal or associated actions} (\ref{cpac}), since, though this presentation was introduced a very long time ago \cite{P75}, we expounded it only in oral lectures, and we think it points out a basic property of the original genuine fibre bundles, which was forgotten, or at least hidden, in subsequent \textsc{Ehres\-mann}'s generalized definition. 

It is also intimately related with the subject of \textsc{A. Kock}'s lecture at the present Conference \cite{K}, and offers for it an ``internalized'' point of view.

\section{Is $\mathbf{Dif}$ really a ``bad'' category \\
or rather a fecund model for ``working mathematicians''?}

\subsection{Has \textbf{Dif} to be ``improved''?}
\label{dif}

It is commonly admitted that the category \textbf{Dif}, consisting of [smooth maps between] smooth\footnote
{We use freely and loosely the term ``smooth'' as an abbreviation for ``differentiable'' with an unspecified class of differentiability. As in \cite{VAR} the manifolds are not assumed to be Hausdorff when not expressly mentioned.}
 manifolds, has bad intrinsic properties, as well as the \emph{forgetful functor} towards \textbf{Set}. \textsc{Ehres\-mann} himself shared this opinion, and, as being basically a Differential Geometer, was much worried with that problem. Indeed, though involved in the search of wide-ranging general concepts, he devoted a lot of energy in keeping contact with the basic applications to manifolds and jets.

A widespread attitude among differential geometers has long been to conclude that Category Theory is uninteresting and useless, at least for their main field of interest.

Another widespread standpoint consists in embedding \textbf{Dif} into a ``better'' category, preferably a \emph{topos} \cite{GIR, Ja, Jb}, or at least a cartesian closed category \cite{ML}. This was done in countless ways by various authors, among which (since this is out of the scope of the present lecture) we shall just mention three:\footnote
{We might have mentioned too for instance \textsc{Frölicher}'s approach by means of \emph{curves} (which is close to \textsc{Souriau}'s approach, as well as \textsc{Chen}'s approach by ``smooth spaces''), and \textsc{Sikorski}'s approach by means of ``differential spaces'', which uses the dual point of view of functions, now prevailing in non-commutative geometry \cite{NCG}.}
\begin{itemize}
	\item \textsc{Ehres\-mann}'s very general processes of completions, enlargements, expansions and others, for which we refer to \cite{OC} and to Andrée's comments;
	\item \textsc{J.-M. Souriau}'s approach by means of \emph{diffeologies}; for an application illustrating this approach we refer to the contribution of \textsc{P. Iglesias} \cite{IGL};
	\item the \textsc{Kock-Lawvere} approach called \emph{Synthetic Geometry}, for which we refer to \textsc{A. Kock}'s comments \cite[(VIII)]{ECO}.
\end{itemize}

Though we do not deny the interest of those sorts of approaches when suitably run, we point out as an overall drawback or danger that they may lead to introduce a huge quantity of pathological objects (for instance not only spaces of leaves of regular foliations, but quotients by arbitrary equivalence relations) which are generally out of control and often quite uninteresting from a differentiable standpoint, and we think that this type of method should be completed by developing characterizations of the ``true'' manifolds among these objects, and techniques for getting precise results about these, rather than ``general'' results about ``more general'' objects.

Indeed one may seriously wonder whether it is really desirable for a differential geometer to be able to ``solve'' very general and unnatural problems of (co)limits by introducing artificial and formal ``solutions'', when neither the problem nor the solution own the least kind of compatibility with the smoothness of the initial data. It seems more natural that he deals primarily with those sorts of singular spaces which, though liable to be very wild, take their origin from smooth data (spaces of leaves of foliations are typical examples, as well as manifolds with corners).

\subsection{Dif as an exemplary category?}
\label{ex}
\subsubsection{Good monics and good epis}

This leads us to propose a \emph{quite different approach} (initiated in \cite{P75}, and explained also in \cite{P03, P04}, to which we refer for more details and examples), which will be developed more extensively elsewhere.
 Our opinion is that the reason why \textbf{Dif} (as well as the forgetful functor towards \textbf{Set}) is considered as ``bad'' from a categorical viewpoint derives from the fact that originally the main founders of Category Theory had quite different and distant models in mind, essentially the categories \textbf{Set} and \textbf{Top} and the categories stemming from \emph{Algebra}, \emph{Algebraic Topology}, \emph{Homological Algebra}, and later \emph{Algebraic Geometry}, and the basic concepts (such as general (co)limits, epis and monics) were coined for those models, as well as the axiomatizations which led notably to the concepts of \emph{abelian categories} and \emph{toposes}.

When compared to these sorts of models, \textbf{Dif} looks poor, since for instance few general (co)limits (notably pullbacks and pushouts) do exist and the general epis and monics may be pathological.

However one may observe that $\emphc{D}=\textbf{Dif}\,$ is endowed with \emph{two distinguished subcategories of monics and epis} respectively, namely the subcategory $\emphc{D}_i$ of \emph{embeddings} (in the sense of \cite{VAR})\footnote
{When dealing with \emph{Hausdorff} manifolds, one has to consider \emph{proper} (i.e. closed) embeddings in order to fulfill the alluded stability properties.}
 and the subcategory $\emphc{D}_s$ of \emph{surmersions} (i.e. surjective submersions \cite{VAR}), which play a basic role in the definition of Lie groupoids (as well as throughout Differential Geometry).

These subcategories have some \emph{remarkable stability properties} (essentially stated in \cite{VAR}\footnote
{These formal properties are indeed intimately linked with variants or consequences (including Godement theorem \cite{SLG, VAR}) of the \emph{implicit function theorem}, in other words with the fact that, for the arrows of $\emphc{D}_i$ and $\emphc{D}_s$, the \emph{local} model coincides with the \emph{linear} model. Since this is the essence of differentiability, it is not amazing that, though these properties look very mild, they are able to encapsulate, as extra data, some essential properties of \textbf{Dif} that are not satisfied in the most general abstract categories.

It is highly remarkable that, in many other categories used by ``working mathematicians'', the same simple formal properties are able to encapsulate some other specific basic theorems.}
) which are not all \emph{simultaneously} satisfied (in an unspecified category) when considering \emph{all}\footnote
{For instance the \emph{product} of two epimorphisms is in general no longer an epimorphism, but this stability property is satisfied by (surjective) open maps (though not by general identification maps), or (surjective) proper maps, or surmersions.}
 general monics and epis (nor any of the general variants the categorists were led to introduce, such as split, strict, regular \ldots) in general categories, but turn out to be (non quite obviously) satisfied by all the monics and epis in the (important and wide-ranging) cases of \emph{abelian categories} as well as \emph{toposes} (and of course \textbf{Set}).

\subsubsection{Internalization}
\label{internal}
We also observed that these properties are precisely those which ensure the existence of the pullbacks (and also a few pushouts) needed for the very definition of Lie groupoids, and for all the basic constructions occurring in the theory of Lie groupoid actions, as it will be described below. 

One may even say that there exists a kind of very powerful empirical and heuristic ``\emph{meta-principle of internalization}'' asserting that translating the algebraic constructions of the theory of groupoid actions in terms of \emph{suitable commutative diagrams} drawn in \textbf{Set}, and stressing the injections, surjections, and pullbacks (and also certain pushouts), allows a ``quasi-automatic'' internalization in the categories where these stability properties are fulfilled (hence particularly in any topos, but also in \textbf{Dif}, which is not a topos, and has not to be embedded beforehand, in a more or less artificial way, into a topos, using such an approach).

\subsubsection{Selecting good limits}\label{select}
Moreover these properties remain still valid for numerous \emph{variants} of these distinguished monics or epis.

\emph{One can say that the \emph{extra data} of $\emphc{D}_i$ and $\emphc{D}_s$ operate as a good choice of sub\-objects, quotients, pullbacks and pushouts, focusing on a class of problems considered as interesting from the smooth standpoint, and preventing from getting lost among too general and too wild objects.}

These problems are described by specific equivalence classes of diagrams involving arrows of $\emphc{D}$, $\emphc{D}_i$ and $\emphc{D}_s$ and possibly certain pullbacks (and pushouts). It is not demanded that these diagrams possess an actual limit inside the category $\emphc{D}$. Indeed one might realize that the introduction of such a limit (in the general sense of Category Theory) would often involve a crucial \emph{loss of information}. For instance we think that the space of leaves of a foliation is better viewed (and indeed more precisely defined) as an equivalence class of foliations in the sense introduced by P. \textsc{Molino} \cite{Mol}, or equivalently as a Morita equivalence class of certain smooth groupoids (in a precise smooth sense), rather than as a universal quotient in the most general categorical meaning( cf. \cite{P04}).

For instance such a quotient would exist as an actual object in the category of diffeologies in the sense of J.-M. ~\textsc{Souriau}; these two approaches lead in natural ways to two non-equivalent definitions of the fundamental group of this object, and it turns out that the former notion is strictly better adapted than the latter.

\subsubsection{Categories for working mathematicians}
We also observed that, very remarkably, a very wide variety of categories used by ``working mathematicians'' may analogously be endowed (often in a quite non obvious way, and in multifarious ways) with two such distinguished subcategories.

It is also very important to observe that, when dealing with categories of commutative diagrams (and primarily \emph{commutative squares}\footnote
{Called ``\emph{quatuors}'' (i.e. quartets) in \textsc{Ehres\-mann}'s terminology.}
) of the original category, there arise (a lot of) natural candidates for special monics and epis which inherit (often in a rather non obvious way) the above-mentioned stability properties, thus creating various new examples. 

In the present framework this applies basically to the category of [morphisms between] Lie (more generally structured) goupoids and their actions. This leads also to an inductive approach of higher order structures (multiple structured groupoids).

We think this approach is exactly in \textsc{Ehres\-mann}'s spirit, as a specialization of the general principles he introduced, which is directly adapted to his main fields of interest in Geometry. 

This leads to what we have called ``diptychs'' \cite{P75}, for which we refer to \cite{P04, P03}, but we recall now some basic definitions and notations for reader's convenience (we intend in a near future to weaken slightly our axioms in order to be able to include a lot of important other structures, such as Riemannian or Poisson manifolds).

\section{A brief review of diptych approach}
\label{sec:Diptap}
At first reading, the technical details may be skipped, just retaining roughly that the following machinery (presented in some details in \cite{P04}, and we intend to develop more extensively in future papers), is aiming at allowing to replace (injections)\,/\,(sur\-jections) of \textbf{Set} by ([closed] embeddings)\,/\,(sur\-mersions) in \textbf{Dif[Haus]} and more generally by (``good'' mono\-morphisms)\,/\,(``good'' epi\-morphisms) in nice categories including very wide-ranging examples, for a more efficient implementation of \textsc{Ehresmann}'s program.

\subsection{(Co)(pre)diptych data and axioms, dip-pairs.}
\label{Dipdax}
\subsubsection{Definitions.}
\label{defdip}
A \emph{prediptych}, denoted by $\mathsf{D}=(\cD;\cD_i,\cD_s)$, or sometimes loosely by $\cD$ alone, is a triple in which:
\begin{itemize}
	\item $\cD$ is a category;
	\item $\cD_i$\,/\,$\cD_s$ is a subcategory of $\cD$, the arrows of which are denoted generally by arrows with a triangular tail\,/\,head such as $\riTo$\,/\,$\rsTo$, or $\rITo$\,/\,$\rSTo$, and so on (here / is written loosely for resp.).
\end{itemize}
Moreover it is demanded that:
\begin{enumerate}
\item[(i)] $\cD_i\cap\cD_s=\cD_\ast$,
\end{enumerate}
where $\cD_\ast$ denotes the sub\-groupoid of invertible arrows (called \emph{isomorphisms}).

\begin{itemize}
	\item $\mathsf{D}$ is called a \emph{diptych} if $\cD$ is equipped with the \emph{extra data} of \emph{finite non void products}, and if the following \emph{extra conditions} are satisfied:

\end{itemize}
\begin{enumerate}
\item[(ii)] $\cD_i$ and $\cD_s$ are \emph{stable by products};
\item[(iii)] (a)\,/\,(b) the arrows of $\cD_i$\,/\,$\cD_s$ are \emph{monos\,/\,epis} (called ``\emph{good}'');
	\item[(iv)]
\begin{enumerate}
	\item ($h=gf\in D_i)\Longrightarrow (f\in D_i)$ (``\emph{strong} $\alpha$-stability'');
	\item $\left((h=gf\in D_s)\,\&\,(f\in D_s)\right)\Longrightarrow (g\in D_s)$ (``\emph{weak} $\beta$-stability'');
\end{enumerate}
\item[(v)](``\emph{transversality}'', denoted by $\cD_s\pitchfork\cD_i$):
\begin{enumerate}
\item[(a)](``\emph{parallel transfer}'') :\\given $s\in\cD_s$, and $i\in\cD_i$ with common target $B$, there \emph{exists} 
\begin{diagram}[h=1.5em,w=2em,tight,labelstyle=\scriptstyle,textflow]
A'&\riTo^{i'}&A\\
\dsTo^{s'}&{\pb}&\dsTo_{s}\\
B'&\riTo_{i}&B\\
\end{diagram}
a pullback square with \emph{moreover} $s'\in\cD_s$, $i'\in\cD_i$, as pictured;

\item [(b)](conversely : ``\emph{descent}'', or ``\emph{reverse transfer}'') :\\ \emph{given} a pullback square as above, but \emph{assuming} now $i'\in\cD_i$, and $s,\!s'\in\cD_s$ (with just given $i\in\cD$), one \emph{concludes} $i\in\cD_i$;
\end{enumerate}
\item[(vi)]pullback squares with all four arrows in $\cD_s$ are \emph{pushouts} \cite{ML} too: 
such squares will be called \emph{s-exact squares}\footnote{
Called \emph{perfect squares} in our previous papers.
}.
\end{enumerate}

Since, for $\cD$ fixed, there may exist on it a great variety of interesting (pre)diptych structures, it may be convenient to say that $(\cD',\cD'')$ is a (\emph{pre})\emph{dip-pair} on $\cD$ whenever $(\cD;\cD',\cD'')$ is a (pre)diptych.

The most interesting diptychs will be \emph{Godement diptychs} (and Godement dip-pairs), to be defined below (\ref{god}).

\emph{Any category $\cD$ may always be considered (implicitely), with its trivial \emph{pre}\-diptych structure} (not a diptych in general):
	\[(\cD;\cD_{\ast},\cD_{\ast}).\tag{triv.pred}
\]
corresponding to the \emph{trivial prediptych pair} $(\cD_{\ast},\cD_{\ast})$.

One may also consider the following canonical \emph{pre}\-diptych structures:
	\[[\cD]_i\un{\text{def}}{=}(\cD;\cD,\cD_\ast)
	\H\text{and}\H
	[\cD]_s\un{\text{def}}{=}(\cD;\cD_\ast,\cD)\,.\tag{i-pred, s-pred}
\]

The \emph{product of two \emph{pre}diptychs} $\sD\times \sD'$ may be defined in an obvious way, but may not be a diptych when $\sD$ and $\sD'$ are, though this is always true when $\sD$ and $\sD'$ have terminal objects (this implies that the objects of $\cD\times \cD'$ are the products of their components).

\subsubsection{Standard and split diptychs}
\label{spe}

The basic diptych:
	\[\boxed{\mathsf{Set} =(\textbf{Set;\,\textbf{Inj},\,Surj})}
\]
has very special additional properties, which are not required for general diptychs.
\begin{itemize}
	\item A diptych will be called \emph{standard} when good monos and epis consist of \emph{all} monos and epis. For a given category, the existence of a standard diptych structure on this category is a very special property.
	\item If moreover all monos and epis are split \cite{ML}, the diptych will be called \emph{split}. (Thus $\mathsf{Set}$ is standard and split.)
\end{itemize}
Note that if $(\cD_i,\cD_s)$ is a dip-pair, the same is true for $(\cD_i,\cD_s^{\text{spl}})$, where $\cD_s^{\text{spl}}$ denotes the class of split arrows of $\cD_s$\,.

{\footnotesize%
\subsubsection{Remarks about axioms and definitions}
\label{remax}\hfill
\begin{enumerate}
	\item In future papers, and with many important new examples in mind, the equality in axiom~ (i) will be weakened into an inclusion, thus defining also \emph{good isomorphisms}. We shall also introduce in the data an \emph{involution} in $\cD$, which here would be trivial.
	\item The existence of the \emph{void} product (i.e. of a \emph{terminal object}), is \emph{not required} in general, in view of basic examples ; when it does exist, it will be denoted by a plain dot $\bullet$ (though it has not to be a singleton), and the diptych is then called \emph{augmented}.

Though, in many examples, not only there exists a terminal object, but moreover the canonical arrows $A\rsTo\bullet$ are in $\cD_s$, this property is \emph{not required} in general. Then those objects $A$ owning this latter property will be called \emph{$s$-condensed}. (In the examples below this may mean: non-empty and moreover connected, compact, finite, maximum...)

As a consequence of the axioms, it turns out that whenever $A$ is $s$-condensed, then, for any object $B$, one has:
	\[\left(A\times B\rsTo^{\text{pr}_2\ \ } B\right)\in\cD_s\ .
\]
\emph{More generally}, an object $A$ will be called \emph{s-basic} when the previous condition is satisfied for $B=A$ (see below \ref{ban}). (For instance, in the standard diptych $\mathsf{Set}$, $\emptyset$ is s-basic though not s-condensed!)

A diptych $\sD$ will be called \emph{s-augmented} when it has a terminal object and \emph{moreover} all its objects are s-condensed (see \ref{loc} below).
	\item Axiom (iv) (a) yields for each arrow $A\stackrel{f}{\rightarrow} B$ of $\cD$ a canonical ``\emph{graph factorization}'':
	\[A\riTo^{\ (\text{id}_A,f)\ \ }A\times B\rTO^{\text{pr}_2\ \ }B
	\tag{graph}
\]
through $\cD_i$ (and through $\cD_s$ when $A$ is $s$-condensed).

This applies in particular to $f=1_A$ and sounds like a \emph{separation property} for $A$: see \cite{P04} for weakening symmetrically this axiom and then introducing the general notion of ``\emph{$i-$scattered}'' objects (which however is \emph{not} the dual of $s$-condensed).

\item Property $\cD_s\pitchfork\cD_i$ of axiom (v) means roughly that the arrows of $\cD_i$ can be pulled back and pushed out along the arrows of $\cD_s$; together with the other axioms (and using the previous remark), it implies the properties $\cD_s\pitchfork\cD$ and $\cD_s\pitchfork\cD_s$, which notably yield \emph{sufficient} conditions for the existence of certain pullbacks.

Those special pullbacks constructed by pulling back along an arrow of $\cD_s$ will be called briefly \emph{special pullbacks}.

All such pullbacks satisfy the following condition (keeping the above notations):
	\[\left((s',i'): A'\riTo B'\times A\right)\in\cD_i\tag{ipb}\ .
\]
\item As shown by elementary examples in \textbf{Dif}, other useful pullbacks may exist, for instance those deriving from the standard transversality condition \cite{VAR}, but also many other ones defined by ``nice'' intersections; however ``\emph{pathological}'' ones may also exist, arising from ``bad'' intersections (though they own the standard universal property of general pullbacks \cite{ML}); it turns out that, in our diptych context (and particularly for \textbf{Dif}: see for instance a useful application in \cite{P86}), the efficient notion is what we call ``\emph{good pull back squares}'', by which we shall mean those pullback squares (in the standard categorical sense) \emph{which furthermore satisfy the} (ipb) \emph{condition} just above (see \ref{tbp} below). This notion is stable by finite products. All special pullbacks (4) are good.

\item Axiom (vi) turns out (in presence of the other axioms) to be equivalent (see \cite{P04}) to demanding the arrows in $\cD_s$ to be \emph{strict} epis, i.e. \emph{coequalizers} \cite{ML}.
\item In many applications the category $\cD$ is endowed with a \emph{forgetful functor} $p:\cD\rightarrow\textbf{Set}$ which is \emph{faithful}; then any morphism $f:A\rightarrow B$ is \emph{uniquely} determined by the triple $(p(f),A,B)$ and $A$, $B$ may be interpreted as ``\emph{structures}'' on the \emph{underlying sets} $p(A)$, $p(B)$; one can speak then of a ``\emph{concrete}'' diptych, a familiar and intuitive situation. But in view of interesting examples, we shall avoid as much as possible such an assumption, trending to work ``without points'', though keeping in mind the concrete situation.

\item It should be noted that the axioms look \emph{nearly self-dual}, though they are \emph{not thoroughly}, and this is important for the most interesting applications, which are seldom self dual (see for instance the basic symmetric simplicial category below in \ref{sym}).

The \emph{dual} (or opposite) of a \emph{pre}diptych $\cD$ is defined as:
	\[\mathsf{D}^{\text{op}}=(\cD^{\text{op}};\left(\cD_s\right)^{\text{op}},\left(\cD_i\right)^{\text{op}}),
		\h\text{i.e.}\h
\left(\cD^\text{op}\right)_i=\left(\cD_s\right)^{\text{op}},\ 
\left(\cD^\text{op}\right)_s=\left(\cD_i\right)^{\text{op}};\tag{dualpred}
\]
A prediptych $\sD$ is called a \emph{co-diptych} if its dual prediptych is a \emph{diptych}: when making explicit the data and axioms for a co-diptych, one has to switch product\,/\,co\-product, epi\,/\,mono, pullback\,/\,pushout.
\end{enumerate}
}

\subsection{Prediptych morphisms and (well) exact diptych morphisms}
\label{dipmor}
A (\emph{pre})\-\emph{diptych morphism}, or shortly (pre)\-dip-morphism:
	\[\boxed{f:\mathsf{D'}\rightarrow\mathsf{D}}\,,
\]
with $\mathsf{D'}=(\cD';\cD'_i,\cD'_s),\ \mathsf{D}=(\cD;\cD_i,\cD_s)$, 
is a functor $f:\cD'\rightarrow\cD$ such that one has:
	\[\boxed{f(\cD'_i)\subset\cD_i,\h f(\cD'_s)\subset\cD_s}\,.\tag{predipmor}
\]
(These conditions vanish when $\mathsf{D'}$ is the trivial prediptych structure (\ref{defdip}) on $\cD'$.)

A \emph{diptych} morphism is called \emph{exact} (or dip-exact when there is some risk of confusion with much stronger notions of exactness \cite{ML}) if moreover it \emph{preserves the special pullbacks} (\ref{remax} (4)).

This type of condition might be regarded as a special case of ``\emph{preservation of second order diptych data}''; this would enable to consider more flexible and more various notions of exactness: see \ref{seco} below.

This does \emph{not} imply the preservation of \emph{finite products}, not required! When finite \emph{non-void} products are preserved too, the morphism is called \emph{well exact}. Nothing is required about preservation of \emph{void} products (terminal objects), whether they exist or not in $\cD$ or $\cD'$.

(Well) (exact) (pre)diptych morphisms define large categories denoted by \textbf{Pre\-dip}, \textbf{Dip}, {\textbf{Dip\-Ex}}, \textbf{Dip\-WEx}.

\section{A few examples of diptychs}
\label{Exdip}
Checking all the axioms for the following examples may be rather lengthy, and of unequal difficulty, but most of them rely on more or less classical results of each theory.

Other examples may be found in \cite{P04}, but we slightly change and improve our notations, and we focus on what is necessary for our present purpose.

\subsection{Main example and its variants: Dif}
\label{mainex}
Apart from $\mathsf{Set}$ (\ref{spe}) the basic example for our present purpose will be (using the terminology of \cite{VAR}):
	\[\boxed{\mathsf{Dif}=(\mathbf{Dif}\,;\mathbf{Emb}\,,\mathbf{Surm})},
\]
where $\mathbf{Dif}$ is the category of [morphisms between] differentiable manifolds, $\mathbf{Emb}$ the subcategory of embeddings, $\mathbf{Surm}$ the subcategory of surmersions (= surjective submersions).
It is augmented but not s-augmented: note that ``s-condensed'' (Rem.~ (2) in \ref{remax}) just means here ``non-empty''.

\emph{Various subcategories} of \textbf{Surm} may also be used for \emph{good epis}, demanding from surmersions \emph{extra properties}, such as being:\\
split (= right invertible), étale, ``retro-connected'' (by which we mean that the fibres are connected, hence the inverse image of any connected subset is still connected)\footnote{
More generally, as a consequence of (non obvious) theorems proved by \textsc{G. Meigniez} (see \cite{MG}), it is possible to replace connected by n-connected. (We hope to apply this important fact, which encapsulates some deep properties, in future papers.)
}.%

This actually defines various sub-diptychs, denoted by:
	\[\mathsf{Dif}^{(\mathsf{spl})},\ \mathsf{Dif}^{(\mathsf{\acute{e}t})},\ \mathsf{Dif}^{(\mathsf{r.c})},\ \mathsf{Dif}^{\mathsf{(r.n-c})}.
\]

For the full subcategory of \emph{Hausdorff} manifolds, one has to use \emph{closed} (or \emph{proper}) embeddings (see the comments of \ref{remax} (3)), which yields:
	\[\boxed{\mathsf{DifH}=(\mathbf{DifHaus};\mathbf{PrEmb},\mathbf{Surm})}.
\]

Then one may demand the surmersions to be also proper [and étale] (which, by a classical theorem of \textsc{Ehresmann}, implies they are \emph{locally trivial fibrations} [finite covers]), and this defines: 
	\[\mathsf{DifH}^{(\mathsf{prop})},\ \ [\mathsf{DifH}^{(\mathsf{fin.cov})}].
\]

According to the above various choices for $\cD_s$, the meaning, for a manifold, of being ``$s-$\emph{condensed'}' in the sense of Rem.~ \ref{remax} (2) may be: non-empty and moreover discrete, connected, $n-$connected, compact, finite, etc. One notices that the condition of being non-empty is dropped for \emph{s--basic} objects.

\subsection{Other large diptychs}
\label{large}
We refer to \cite{P04} for examples with topological spaces, Banach spaces, vector bundles... We just stress that, in the category \textbf{Top}, one cannot take as arrows of $\cD_s$ \emph{all} the identification maps (not stable by products!) but one can take the \emph{open} ones (with possibly many variants): the diptych we get in this way will be denoted by $\mathsf{Top}$.

It is remarkable that any \emph{topos} or \emph{abelian category}, may be endowed canonically with a ``\emph{standard}'' \emph{diptych structure} (\ref{spe}), (and also a standard \emph{co-diptych} one), taking for $\cD_i$ and $\cD_s$ \emph{all} the monomorphisms and \emph{all} the epimorphisms (as in \textbf{Set}): the statements and the proofs of most (not all) of the non-obvious axioms may be found for instance in (or deduced from) \cite{Ja} and \cite{ML}. \emph{That is no longer possible for any general category} (notably for \textbf{Dif}).

\emph{We think this explains why it seems easier (from the point of view of a ``pure categorist'') to work within toposes, but it is much more efficient, though no more difficult, for a ``working mathematician'', to use diptych structures, which allow to encapsulate some crucial and sensible extra information, while granting wider-ranging applications.}

We feel that the \emph{existence}, for a given category, of a \emph{diptych} structure on it, is an important \emph{property} of this category, which makes it suitable for ``working mathematicians''.

Note that the (large) category \textbf{Poset} of (\emph{partially})\footnote{%
We recall that in \textsc{Bourbaki}'s terminology \cite{ENS} (pre)\emph{order} means \emph{partial} (pre)\emph{order}.
} %
ordered sets and monotone functions owns a standard diptych structure, denoted by {\sffamily Poset}, and also a co-diptych structure.
\subsection{Small (pre)diptychs}
\label{small}
Small, and even finite (often non-concrete: cf. Rem. \ref{remax}~ (7)) (pre)\-dip\-tychs will be also of interest in the sequel. Here $\cD$ will be a small category $\cT$ (called ``type of diagram'' below).

\subsubsection{Preordered sets, finite ordinals}\label{preoset}
A (partially\footnote{
See previous footnote.%
}) pre\-ordered set $(E,\leq)$ may be viewed as a (very special) category $\cT$, in which all arrows are both monos and epis.

The standard diptych structure (\ref{spe}) cannot exist, but in the degenerate case of an equivalence relation.

For $a,b\in E$, their [co]product means just their g.l.b. [l.u.b.] or meet $a\wedge b$ [join $a\vee b$] \cite{ML}. It is defined for any pair whenever $(E,\leq)$ is a \emph{lattice}, and notably when the pre\-order is \emph{total}.

When such is the case, the \emph{trivial} pre\-diptych $(\cT;\cT_\ast,\cT_\ast)$ is a (non-standard) \emph{diptych}, as well as (notably) $[\cT]_i$ and $[\cT]_s$ (defined in \ref{defdip}).

For any $n\in\bbN^\star$ we denote by 
	\[\boxed{\mathbf{n}=[n-1]=\{0<1<\ldots<n-1\}\ (=[0,n[\,\subset\bbN)}\tag{ord}
\]
the corresponding (non-empty) \emph{ordinal}, possibly viewed as a \emph{trivial} small diptych.

Forgetting the order, one denotes by $\ov{\times}{\mathbf{n}}$ or $\langle n-1\rangle$ the \emph{banal groupoid} of pairs $n\times n$ (corresponding to the trivial total pre\-order). We have the inclusion of categories:
	\[\boxed{[n-1]=\mathbf{n}\subset\ov{\times}{\mathbf{n}}=\langle n-1\rangle}\,.
\]

It may be sometimes suggestive and useful to write:
	\[\boxed{\textbf{1}=[0]=\pmb{\star}\,,\ \textbf{2}=[1]=\,\pmb{\downarrow}\,,\ \ov{\times}{\mathbf{2}}=\bbI}
\]

For instance the ordinal $\mathbf{2}=\{0<1\}$ owns \emph{three canonical diptych structures}, denoted by: 
$\textbf{2}$ (the trivial one), $[\textbf{2}]_i$, and $[\textbf{2}]_s$, 
pictured respectively by:
	\[0\rTO 1\,,\H 0\riTo 1\,,\H0\rsTo 1\,.
\]
Taking products, we get nine diptychs, for which we introduce transparent \emph{geometrical notations}, for instance:

\begin{eqnarray}\label{diagex}
	&\pmb{\square}=\mathbf{2}\times\mathbf{2},&
	&\un{i}{\pmb{\square}}=\left[\mathbf{2}\right]_i\times\mathbf{2},&
	&\un{s}{\pmb{\square}}=\left[\mathbf{2}\right]_s\times\mathbf{2},&
	&\li\pmb{\square}=\mathbf{2}\times\left[\mathbf{2}\right]_i,&
&\ls\un{i}{\pmb{\square}}=\left[\mathbf{2}\right]_i\times\left[\mathbf{2}\right]_s&&\cdots\\ \notag
&\text{pictured}&& \text{by}&& \text{the}&& \text{following}&& \text{diagrams:}&&\\ \notag
&
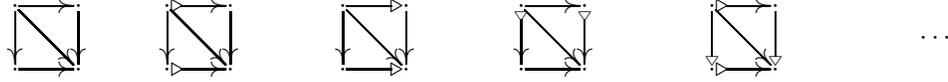
\begin{diagram}[size=1.2em,tight,inline,abut]%
\cdot&\rTO&\cdot\\
\dTO&\rdTO&\dTO\\
\cdot&\rTO&\cdot\\
\end{diagram}&&
\begin{diagram}[size=1.2em,tight,inline,abut]%
\cdot&\riTo&\cdot\\
\dTO&\rdTO&\dTO\\
\cdot&\riTo&\cdot\\
\end{diagram}
&&
\begin{diagram}[size=1.2em,tight,inline,abut]%
\cdot&\rsTo&\cdot\\
\dTO&\rdTO&\dTO\\
\cdot&\rsTo&\cdot\\
\end{diagram}
&&
\begin{diagram}[size=1.2em,tight,inline,abut]%
\cdot&\rTO&\cdot\\
\diTo&\rdTO&\diTo\\
\cdot&\rTO&\cdot\\
\end{diagram}
&&
\begin{diagram}[size=1.2em,tight,inline,abut]%
\cdot&\riTo&\cdot\\
\dsTo&\rdTO&\dsTo\\
\cdot&\riTo&\cdot\\
\end{diagram}
&&\cdots
\notag
\end{eqnarray}

\subsubsection{Groupoid data: $\Upsilon$}\label{grouda}
The pre\-diptych generated by the following diagram (no longer a preorder), denoted by $\mathsf{\Upsilon}$, will be used for describing structured groupoids:
	\[\begin{diagram}[w=3em,tight,labelstyle=\ssst]
2&\rsTo~{\delta}&1&\pile{\rsTo~{\alpha}\\ \liTo~{\omega}\\}&0\\
\end{diagram}\tag{$\Upsilon$}
\]
where $\alpha\,\omega$ is the unit of the object $0$ (see \ref{sym} below for more precision).

\subsubsection{Simplicial category: $\mathbf{\Delta},\pmb{\nabla}$}
\label{simp}
The \emph{non empty finite ordinals} $\mathbf{n}=[n-1]$ of \ref{preoset} are the objects of the simplicial category $\boldsymbol{\Delta}$. Its arrows (\emph{order preserving maps}) may be viewed as exact diptych morphisms. $\mathbf{1}$ is a terminal object.

Though, when including also the empty ordinal \textbf{0}, we get a monoidal category \cite{ML} (here denoted by $\dot{\mathbf{\Delta}}$, and called ``\emph{augmented}''), $\mathbf{\Delta}$ has no products nor co\-products, hence \emph{no diptych structure}, but a canonical (standard) \emph{pre\-diptych} structure, denoted by:
	\[\mathsf{\Delta}=(\mathbf{\Delta};\mathbf{\Delta}_i,\mathbf{\Delta}_s)
\]
 defined by the \emph{injective} and \emph{surjective} morphisms.

Indeed the useful prediptych will be the \emph{dual} one (\ref{remax} (8)): 
	\[\nabla=(\,\pmb{\nabla};\pmb{\nabla}_i,\pmb{\nabla}_s),\text{ where }%
	\ \pmb{\nabla}=\mathbf{\Delta}^{\text{op}},\h
	\pmb{\nabla}_i=(\mathbf{\Delta}_s)^{\text{op}},\h
	\pmb{\nabla}_s=(\mathbf{\Delta}_i)^{\text{op}}\ .
\]
(It may be sometimes useful to \emph{embed} $\mathsf{\Delta}$ in the \emph{diptych of finite \emph{pre}ordered sets}.)
\subsubsection{Symmetric simplicial category: $\pmb{\mathsf{F}},\pmb{\Finv}$}
\label{sym}
More important indeed for structured groupoids\footnote{%
As above-mentioned we focus here on specific properties of structured \emph{groupoids}. It might happen that our approach could be extended to \textsc{Ehresmann}'s structured \emph{categories} (with necessarily much weaker properties), replacing the products data by a monoidal structure, in order to include $\pmb{\nabla}$, but this out of our present scope.
} %
than $\mathbf{\Delta}$ and $\pmb{\nabla}$ will be the category $\pmb{\mathsf{F}}$\footnote{
Often denoted by $\mathbf{\Phi}$ in the literature.
}: 
it has the \emph{same objects} as~ $\mathbf{\Delta}$ (non-empty finite ordinals), but, forgetting the order, its arrows are \emph{all} the set-theoretic maps between the underlying finite cardinals. We call it the \emph{symmetric} simplicial category.

We remind (\ref{preoset}) the embedding:
	\[\textbf{n}\subset\ov{\times}{\mathbf{n}}\,,
\]
where $\ov{\times}{\mathbf{n}}$ denotes
 the \emph{banal} groupoid $\ov{\times}{\mathbf{n}}$ of pairs (such an object will be thought below as modeling all the commutative diagrams in groupoids with $n$ vertices). The arrows of $\,\pmb{\mathsf{F}}$ may be then regarded as (degenerate) groupoid morphisms.

Observe \emph{injections} and \emph{sur\-jections} define on $\pmb{\sF}$ a (\emph{standard}) \emph{diptych structure}:
	\[(\pmb{\mathsf{F}};\pmb{\mathsf{F}}_i,\pmb{\mathsf{F}}_s),
\]
which is \emph{s-augmented} in the sense of Rem. \ref{remax} (2), with terminal object $\mathbf{1}$, but this will not be the most important one for the sequel.

\emph{Much more important} for our purpose will be its \emph{co-diptych structure} (Rem. \ref{remax} (8)), defined by considering \emph{co-products} and \emph{amalgamated sums}, instead of products and fibred products; we denote this \emph{co-diptych} structure by $\mathsf{F}$.

The \emph{dual diptych} (which plays a basic role for groupoids) is denoted by:
	\[\Finv=
	(\pmb{\Finv};\pmb{\Finv}_i,\pmb{\Finv}_s)
	\h\text{with}\h
	\pmb{\Finv}=\pmb{\mathsf{F}}^{\text{op}},\ 
	\pmb{\Finv}_i=(\pmb{\mathsf{F}}_s)^{\text{op}},\ 
	\pmb{\Finv}_s=(\pmb{\mathsf{F}}_i)^{\text{op}}.
\]
It has \emph{no terminal object} (hence cannot be s-augmented), but \textbf{1} is now an \emph{initial} object, and moreover any object \textbf{n} is an \emph{n}-th \emph{iterated product} of \textbf{1}.

Considering the full subcategories of $\pmb{\nabla}$ or $\pmb{\Finv}$ generated by the objects $\mathbf{p}$ or $\ov{\times}{\mathbf{p}}$ for $0\leq p-1\leq n$, we get the \emph{truncated prediptychs} (in which the products are no longer defined): $\nabla^{[n]}$ or $\Finv^{[n]}$.

We have \emph{canonical \emph{pre}diptych embeddings}:
	\[\mathsf{\Delta}\subset\mathsf{F},\ \ \mathsf{\Upsilon}\subset\nabla^{[2]}\subset\nabla\subset\Finv.
\]
With the \emph{notations} used by \cite{ML} for the \emph{canonical generators} of the simplicial category $\mathbf{\Delta}$, the pre\-diptych $\mathsf{\Upsilon}$ (described in \ref{grouda}) is precisely identified with the \emph{dual} of the following sub-pre\-diptych of $\mathsf{F}$ (remind $\textbf{n}=[n-1]$, cf. \ref{preoset}):
\[\begin{diagram}[w=3em,tight,labelstyle=\ssst]
\textbf{3}&\liTo~{\delta^2_0}&\textbf{2}&\pile{\liTo~{\delta_1^1}\\ \rsTo~{\sigma^1_0}\\}&\textbf{1}\\
\end{diagram}\tag{$\pmb{\curlywedge}$}
\]
It should be noted that the pair $(\alpha,\alpha)$ owns a pullback in $\pmb{\Finv}$, but not in $\pmb{\nabla}$.

When \emph{adding the empty ordinal} $\textbf{0}$, we shall speak (in a slightly improper way) of \emph{augmented} simplicial categories and use the notations:
	\[\dot{\mathbf{\Delta}},\ \dot{\pmb{\nabla}},\ \dot{\pmb{\mathsf{F}}}\ 
	\dot{\pmb{{\Finv\ }}}.
\]
{\footnotesize
[The slight ambiguity of this terminology comes from the fact that the \emph{standard} diptych on the category $\dot{\pmb{\mathsf{F}}}$ is \emph{no longer s-augmented}, since now ``\,s-condensed'' (Rem. \ref{remax} (2)) would mean ``\emph{non-empty}''! On the other hand $\dot{{\Finv\ }}$ becomes now an \emph{``\,s-augmented diptych''}, with $\textbf{0}$ as \emph{terminal} object\footnote{
The (slight but basic) dissymmetry of diptych axioms stresses this basic dissymmetry between 0 and 1, but we shall let the philosophers argue about Being and Nothingness!
} (see also \ref{Finv} below), and this is a property of the \emph{co}-diptych structure $\dot{\sF}$ of $\dot{\pmb{\mathsf{F}}}$.]
}%

\subsubsection{Butterfly diagram type}
\label{but}

We introduce now a pre\-diptych which will occur frequently in the sequel. Here $\overline{\textbf{n}}\,$ may be regarded as just a second (distinct) copy of the ordinal \textbf{n} (indeed it should be endowed with the reverse order, but this is immaterial in the present context). We use \emph{notations of} \cite{ML} \emph{for the canonical generators} of the \emph{simplicial category} (canonical injections and surjections), and denote by $\iota_1,\,\iota_2$ the canonical injections of \textbf{2} into $\textbf{4}=\textbf{2}+\textbf{2}$. The pre\-diptych pictured by the right diagram below is defined as the \emph{dual} (= opposite) of the left one; it will be denoted by $\pmb{\bowtie}$ and called the \emph{butterfly diagram type} (see \ref{cancon} below for explaining its role for groupoids and notations).

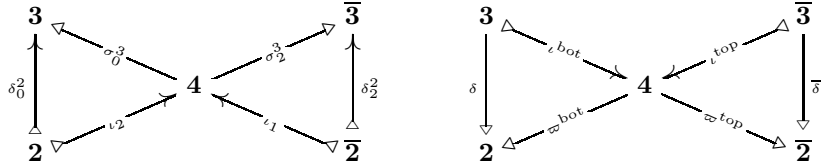
\begin{center}
\begin{diagram}[h=1.3em,w=3em,tight,inline,labelstyle=\ssst]
\textbf{3}&&&&\overline{\textbf{3}}\\
&\lusTo~{\sigma_0^3}&&\rusTo~{\sigma_2^3}&\\
\uiTo^{\delta_0^2}&&\textbf{4}&&\uiTo_{\delta_2^2}\\
&\ruiTo~{\iota_2}&&\luiTo~{\iota_1}&\\
\textbf{2}&&&&\overline{\textbf{2}}\\
\end{diagram}\H\H
\begin{diagram}[h=1.3em,w=3em,tight,inline,labelstyle=\ssst]
\textbf{3}&&&&\overline{\textbf{3}}\\
&\rdiTo~{\iota^{\text{bot}}}&&\ldiTo~{\iota^{\text{top}}}&\\
\dsTo^\delta&&\textbf{4}&&\dsTo_{\overline{\delta}}\\
&\ldsTo~{\varpi^{\text{bot}}}&&\rdsTo~{\varpi^{\text{top}}}&\\
\textbf{2}&&&&\overline{\textbf{2}}\\
\end{diagram}
\end{center}

\subsubsection{Some canonical endomorphisms of $\protect\Finv$. }
\label{canend}

As will be seen below (\ref{cancon}), 
dip-exact endo-functors of the canonical diptych $\Finv$ (and their natural transformations) 
define natural canonical operations on structured groupoids (and \emph{this will explain the following notations}: see \ref{cancon} below).

We shall denote by $\pmb{\sE}_0$ the \emph{monoid} of \emph{endo-functors} of $\,\pmb{\sF}$, and by $\pmb{\sE}$ the category (with base $\pmb{\sE}_0$) of its \emph{endomorphisms}, i.e. natural transformations $\gamma:\Gamma\stackrel{\centerdot}{\rightarrow}\Gamma'$ between such endo-functors (now regarded as objects), endowed with the composition law which is \emph{called vertical in} \cite{ML}.

However, to be coherent with our future conventions introduced below for diagrams (\ref{transpo}), \emph{we shall switch the names of horizontal and vertical compositions} in the present terminology.

By duality the \emph{identity map} from $\pmb{\Finv}$ to $\pmb{\sF}$ (``reversing the arrows'', which is a \emph{contra\-variant} functor) determines a \emph{bijection} between endo-functors $\Gamma,\,\Gamma'$ of the category $\pmb{\Finv}$ and natural transformations $\gamma:\Gamma\stackrel{\centerdot}{\rightarrow}\Gamma'$ on the one hand, and the corresponding objects denoted by $\ov{\ast}{\gamma}:\ov{\ast\ast}{\Gamma}\ov{\centerdot}{\leftarrow}\ov{\ast\ast}{\Gamma'}$ for the dual category $\pmb{\sF}$ on the other hand, and one notices that this bijection is \emph{horizontally contra\-variant} and \emph{vertically covariant} (\emph{with our present conventions}).

More precisely, denoting now by $\pmb{\exists}_0$ the monoid of endo-functors of $\pmb{\Finv}$ and by $\pmb{\exists}$ the category of its endomorphisms, endowed with what we prefer to call here the \emph{horizontal} law, the previous bijection determines two isomorphisms:
\begin{itemize}
	\item a \emph{covariant} isomorphism denoted by:\h\h $\ast\ast:\pmb{\exists}_0\longleftrightarrow\pmb{\sE}_0$;
	\item a \emph{contravariant} isomorphism:\H\H\H $\ast:\pmb{\exists}\ \longleftrightarrow\pmb{\sE}$.
\end{itemize}

The following three endo-functors of $\pmb{\sF}$ (and their dual associates), which turn out to be \emph{exact}, and which indeed \emph{extend} to the augmented category $\dot{\pmb{\sF}\,}$, will be of special importance.
\begin{itemize}
	\item To each object \textbf{n} of $\,\dot{\pmb{\sF}}$, there is attached an involutive arrow in $\,\dot{\pmb{\sF}}$, the ``\emph{mirror map}'' $\varsigma_\textbf{n}:\textbf{n}\rightarrow\textbf{n}$, which ``\emph{reverses the order}'' ($p\mapsto n-1-p$). This defines an involutive automorphism of $\pmb{\sF}\,$(or $\dot{\pmb{\sF}}$):
	\[\boldsymbol{\Sigma}:\pmb{\sF}\longleftrightarrow\pmb{\sF}\,,
	\tag{Mirror}
\]
preserving the objects and taking any arrow 
$f:\textbf{m}\rightarrow\textbf{n}$ to 
$\boldsymbol{\Sigma} f\un{\text{def}}{=}\varsigma_\textbf{n}\circ f \circ\varsigma_\textbf{m}\,$, 
and then an involutive natural isomorphism $\textbf{n}\mapsto\varsigma_\textbf{n}$ (between functors):
	\[\varsigma:\textbf{Id}\stackrel{\centerdot}{\longleftrightarrow}\boldsymbol{\Sigma}\,.
	\tag{mirror}
\]
We shall omit here the $\ast\ast$ for the dually associated objects.

(Note $\boldsymbol{\Sigma}$ induces an automorphism of $\boldsymbol{\Delta}\subset\pmb{\sF}$, but $\varsigma$ is no longer defined.)

	\item The (right\,/\,left) ``\emph{shift endo-functors}''\footnote{
	The use of symbol $\Delta$ seems rather traditional for denoting the simplicial category (augmented in \cite{ML}), and respected here. However this symbol and its dual $\nabla$ are now given, here and in the sequel, operator significances having a very strong \emph{geometric flavour} evoking \emph{commutative triangles} (which will be progressively revealed). The \emph{risk of confusion} seems negligible, according to the context.
	} 
	$\ov{\ast\ast}{\mathbf{\Delta}}$\,/\,$\ov{\ast\ast}{\pmb{\nabla}}$ of$\ {\pmb{\sF}}$ (or $\dot{\pmb{\sF}}$) are defined by:
	\[(\mathbf{n}\mapsto\mathbf{1+n},\ \lambda\mapsto\text{id}_{\textbf{1}}+\lambda)\h/\h
	(\mathbf{n}\mapsto\mathbf{n+1},\ \lambda\mapsto\lambda+\text{id}_{\textbf{1}})\,.
\]
They are endowed with natural ``shifting'' transformations (from the identical functor):
\begin{center}
\bdi[w=2.5em,inline,labelstyle=\sst]
\ov{\ast\ast}{\boldsymbol{\Delta}}&
\lTO^{\centerdot}_{\ov{\ast}{\delta}}&\textbf{Id}\\
\edi
\H/\H
\bdi[w=2.5em,inline,labelstyle=\sst]
\ov{\ast\ast}{\pmb{\nabla}}
&\lTO^{\centerdot}_{\ov{\ast}{\overline{\delta}}}&\textbf{Id}\\
\edi\,.
\end{center}
which are defined (using Mac Lane's notations \cite{ML} for the canonical generators of the simplicial category) by:

\[\textbf{n}\mapsto\delta_\textbf{n}\un{\text{def}}{=}\delta^n_0:\textbf{n}\rightarrow\textbf{1}+\textbf{n}\h/\h
\textbf{n}\mapsto\overline{\delta}_\textbf{n}\un{\text{def}}{=}\delta^n_n:\textbf{n}\rightarrow\textbf{n}+\textbf{1}\,.
\]
By duality this yields for $\dot{\Finv\ }$ exact endo-functors and natural transformations:
\begin{center}
\bdi[w=2.5em,inline,labelstyle=\sst]
\boldsymbol{\Delta}&\rTO^{\centerdot}_\delta&\textbf{Id}\\
\edi
\H/\H
\bdi[w=2.5em,inline,labelstyle=\sst]
\pmb{\nabla}
&\rTO^{\centerdot}_{\overline{\delta}}&\textbf{Id}\\
\edi\,,
\end{center}

One has: $\pmb{\nabla}=\boldsymbol{\Sigma}\circ\boldsymbol{\Delta}\circ\boldsymbol{\Sigma}$.

	\item The ``\emph{doubling}'' or ``\emph{co-squaring}'' endo-functor of$\ \dot{\pmb{\sF}}$ is defined by:
		\[\textbf{n}\mapsto\textbf{n}+\textbf{n},\H\H\lambda\mapsto\lambda+\lambda\,,
\]
and is endowed with natural transformations 
defined by the co-diagonal and the canonical injections into the co-product:
	\[\textbf{n}\leftarrow\textbf{n}+\textbf{n}\,(=\textbf{2n}),\H
	\textbf{n}\rightarrow\textbf{n}+\textbf{n}\leftarrow\textbf{n}\,.
\]
By duality ones gets for $\dot{\Finv\ }$ an exact endo-functor $\pmb{\square},$\footnote{
	This symbol has been used above (\ref{preoset}) with a quite different meaning to denote the (concrete) category $\mathbf{2}\times\mathbf{2}$. The context should avoid confusion.
	}
	 and natural transformations, which we denote by:
\begin{center}
\bdi[w=2.3em,inline,labelstyle=\sst]
\textbf{Id}&\rTO^{\centerdot}_\iota&\pmb{\square},&&\textbf{Id}&\lTO^{\centerdot}_{\varpi^{\ssst{\text{bot}}}}&
\pmb{\square}&\rTO^{\centerdot}_{\varpi^{\ssst{\text{top}}}}&\textbf{Id}
\edi
\,.
\end{center}
\end{itemize}
One can also write a canonical butterfly type diagram in the category $\pmb{\exists}$ (see below \ref{cancon} for the definition of the natural transformations $\iota^{\text{bot}}$ and $\iota^{\text{top}}$ and the significance of this diagram for structured groupoids):

\begin{center}
\begin{diagram}[h=1.8em,w=4.5em,tight,inline,labelstyle=\sst]
\pmb{\Delta}&&&&\pmb{\nabla}\\
&\rdiTo^{\centerdot}_{\iota^{\text{bot}}}&&\ldiTo^{\centerdot}_{\iota^{\text{top}}}&\\
\dsTo_{\centerdot}^{\delta}&&\pmb{\square}&&\dsTo_{\centerdot}^{\overline{\delta}}\\
&\ldsTo^{\centerdot}_{\varpi^{\text{bot}}}&&\rdsTo^{\centerdot}_{\varpi^{\text{top}}}&\\
\textbf{Id}&&&&\textbf{Id}\\
\end{diagram}\H.
\end{center}

{\footnotesize
\subsubsection{Canonical representation of $\protect\pmb{\Finv}$}
\label{repr}

It may be convenient to get some concrete representation of the abstract category $\pmb{\Finv}$.

Let us denote by \textbf{CLat} the category of [morphisms between] \emph{complete lattices}, which owns a standard (i.e. injections, surjections) (\ref{spe}) diptych structure $\mathsf{CLat}$. Then one defines a functor:
	\[\gP^\ast:\pmb{\Finv}\rightarrow\textbf{CLat}
\]
as follows (denoting by $f:\textbf{q}\rightarrow\textbf{p}$ an arrow of $\,\pmb{\sF}$ and by 
$\overline{f}:\textbf{p}\rightarrow\textbf{q}$ the dually associated arrow of $\pmb{\Finv}=\pmb{\sF}^{\text{op}}$, and using power sets and inverse images):
	\[\gP^\ast(\textbf{p})\un{\text{def}}{=}2^\textbf{p}=\gP(\textbf{p})\,;\h
	\gP^\ast(\overline{f})\un{\text{def}}{=}f^{-1}:\gP(\textbf{p})\rightarrow\gP(\textbf{q})\,.
\]
(Observe that identifying $A\in\gP(\textbf{p})=2^\textbf{p}$ with its \emph{characteristic map} $a:\textbf{p}\rightarrow\textbf{2}=\{0,1\}$, then $B=f^{-1}(A)$ is interpreted just as $b=a\circ f$).

Then it is a lengthy but elementary exercise in set theory to prove what follows:
\begin{itemize}
	\item $\gP^\ast:\Finv\rightarrow\mathsf{CLat}\rightarrow\mathsf{Set}$ is a \emph{faithful well exact diptych morphism} (\ref{dipmor}) (called the canonical embedding).
	\item $\gP^\ast:\Finv\rightarrow\mathsf{Set}$ may be identified with the \emph{symmetric nerve} (see below \ref{grouner}) of the \emph{banal groupoid} $\ov{\times}{\textbf{2}}\ (= 2\times2)$ of pairs of elements of the set $2=\{0,1\}$.
\end{itemize}

When $f:\textbf{q}\rsTo\textbf{p}$ is a \emph{surjection} in $\pmb{\sF}$ (defining an \emph{equivalence relation on} \textbf{q}), the good mono $\overline{f}:\textbf{p}\riTo\textbf{q}$ in $\pmb{\Finv}\,$ is represented (using $\gP^\ast$) by the embedding onto the subset of \emph{saturated} subsets of \textbf{q}.
When $f:\textbf{q}\riTo\textbf{p}$ is an \emph{injection}, the associated good epi $\overline{f}:\textbf{p}\rsTo\textbf{q}$ is represented by the surjection  ``\emph{trace}'' of a subset of \textbf{p} \emph{on the image} of $f$ in \textbf{p}.
}

\section{Squares in a diptych}
\label{dipsq}
A preliminary study of squares is needed for being able to study below more general diagram morphisms. (We modify and improve notations of \cite{P04}.)
\subsection{Three basic properties of squares in a (pre)diptych}

\subsubsection{Definitions and notations}
\label{tbp}

It turns out that, in a given prediptych $\mathsf{D}$, the following three types of commutative squares (introduced in \cite{P85, P89, P04} \emph{with different names and notations}) will play a \emph{basic role} in the present diptych framework\footnote{
The underlying purely algebraic properties are more or less classical in various contexts, with unsettled terminologies.
}. 
Here a (commutative)\footnote{
If the opposite is not expressly mentioned, \emph{all squares are commutative} (``\emph{quatuors}'', i.e. quartets, in \textsc{Ehresmann}'s terminology).
} square of $\cD$, with the notations of the diagrams just below, will be called:
\begin{itemize}
	\item an \textbf{ipb-square} if condition (ipb) of Rem. \ref{remax} (4) above is satisfied, which reads here:
	\[\boxed{A'\riTo^{(f',u)\ } B'\times A}\ ;\tag{ipb}
\]

	\item a \textbf{gpb-square} or \textbf{good pullback} (see Rem. \ref{remax} (5) above) if \emph{moreover} it is a pullback (or cartesian) square \cite{ML};
	\item an \textbf{spb-square} if the pair $(v,f)$ admits a \emph{good pullback}\footnote{
Observe this is, somewhat unsymmetrically, \emph{not} demanded for an ipb-square (in spite of the symmetrical terminology), but is notably automatically verified in the frequent cases when $f$ \emph{or} $v$ is in $\cD_s$ (``special'' ipb-squares).
} 
	and \emph{moreover} the canonical factorization is in $\cD_s$, i.e. one has: 
	\[\boxed{A'\rsTo^{(f',u)\ } B'\underset{B}{\times}A\riTo B'\times A}\,.\tag{spb}
\]
Note that: $\boxed{\text{(gpb)}\Longleftrightarrow\,\big(\,\text{(ipb)}\ \&\ \text{(spb)}\,\big)}$\,.
\end{itemize}

The ipb, gpb, spb -squares will be tagged respectively by $\ipb,\,\gpb,\,\spb$ as follows:
\begin{center}
\begin{diagram}[h=2em,w=2em,tight,labelstyle=\scriptstyle,inline]
A'&\rTO^{u}&A\\
\dTO^{f'}&{\ipb}&\dTO_{f}\\
B'&\rTO_{v}&B\\
\end{diagram}\h
\begin{diagram}[h=2em,w=2em,tight,labelstyle=\scriptstyle,inline]
A'&\rTO^{u}&A\\
\dTO^{f'}&{\gpb}&\dTO_{f}\\
B'&\rTO_{v}&B\\
\end{diagram}\h
\begin{diagram}[h=2em,w=2em,tight,labelstyle=\scriptstyle,inline]
A'&\rTO^{u}&A\\
\dTO^{f'}&{\spb}&\dTO_{f}\\
B'&\rTO_{v}&B\\
\end{diagram}\h.
\end{center}

It derives from axioms (iii)(b) and (vi) that spb-squares with $f$ \emph{and} $v$ (hence also, using $\cD_s\pitchfork\cD_s$, $f'$ and $u$) in $\cD_s$ are \emph{pushouts}.

\subsubsection{Some remarks}
\label{remsq}

Observe these properties are obviously stable by \emph{transposition} (i.e. diagonal symmetry); this enables \emph{doubling the following remarks}.

Note that, by axiom (iv)(a):
$\,(u\in \cD_i)\Longrightarrow\text{(ipb)}$.

One has also obviously: $\big((v\in\cD_\ast)\ \&\ (u\in\cD_s)\big)\Longrightarrow\text{(spb)}$.

There are also useful (horizontal or vertical) ``\emph{parallel transfer}'' properties:
	\[\bullet\ \text{in an \emph{ipb-square}: } (f\in\cD_i)\Rightarrow(f'\in\cD_i);\h
	\text{in an \emph{spb-square}: } (f\in\cD_s)\Rightarrow(f'\in\cD_s).
\]
The property on the left relies on the \emph{graph factorization} of any given square:

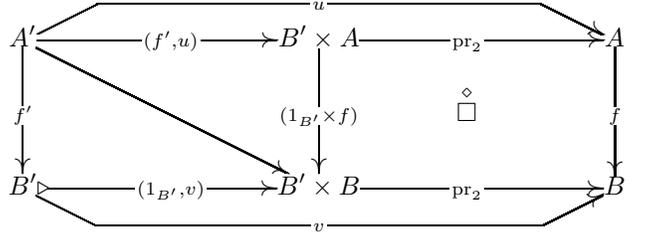
\begin{center}
\begin{diagram}[h=1.4em,w=2.8em,tight,labelstyle=\scriptstyle,inline,abut]
&{}&&&\rLine~u&&&{}&\\
A'\ruLine(1,1)&&\rTO~{(f',u)}&&B'\times A&&\rTO~{\text{pr}_2}&&A\rdTO(1,1)\\
&\rdTO(4,4)&&&&&&&\\
\dTO~{f'}&&&&\dTO~{(1_{B'}\times f)}&&\gpb&&\dTO~{f}\\
&&&&&&&&\\
B'&&\riTo~{(1_{B'},v)}&&B'\times B&&\rTO~{\text{pr}_2}&&B\\
&{}\rdLine(1,1)&&&\rLine~v&&&{}\ruTO(1,1)&\\
\end{diagram}\h,
\end{center}
using axioms (ii) and (iv)(a). The right property relies on $\cD_s\pitchfork\cD$ (\ref{remax}~ (4)).

The composition of squares gives rise to remarkable \emph{stability properties}, on which the constructions of the next section are relying. We cannot list here these properties, which may be proved in a purely diagrammatic way, using only diptych axioms. They are stated in \cite{P85} (Prop. A2) with a different terminology, and in a more restrictive context.

\subsection{Naïve, main and canonical diptych structures for squares in a diptych.}
Here $\mathsf{D}=(\cD;\cD_i,\cD_s)$ will be a \emph{diptych} (though the most obvious properties extend to pre\-diptychs).

We shall not give the (very geometric) proofs of the statements asserted in this section, since they are sometimes lengthy, but never very hard, once the suitable, sometimes (hyper)cubic, diagrams have been drawn
, chasing and juggling with the two compositions of squares and with diptych axioms.

\subsubsection{A piece of notation for squares}\label{notsq}

We denote by $\sq \cD$\footnote{
\textsc{Ehresmann}'s very geometric and useful notation. This is actually a basic example of \emph{double} category, with the horizontal and vertical composition laws satisfying the interchange law. The precise relation with the notation introduced for denoting the category $\mathbf{2}\times\mathbf{2}$ \eqref{diagex} will appear in formula (sq) of \ref{transpo} (2).
} 
 (or sometimes $\un{\ssst{\text{hor}}}{\square}\ \cD$, contrasting with $\ov{\ssst{\text{vert}}}{\square}\ \cD$, when more precision is required) 
the category of (commutative) squares with the \emph{horizontal} composition law (the \emph{objects} are \emph{vertical} arrows). The symbol $\sq$ is functorial; this means that any functor $f:\cD\rightarrow\cD'$ extends to a new functor:
	\[\sq f:\sq\cD\rightarrow\sq\cD'\,.\tag{sq.funct}
\]

Let $\cD_v$, $\cD_h$ be \emph{two subcategories} of $\cD$, which below may be $\cD$, $\cD_i$ or $\cD_s$:\emph{ without extra notification the symbols $h$ and $v$ are allowed to take independently the values $i$, $s$, or blank}. In what follows, transparent conventions will point to \emph{constraints on the vertical\,/\,horizontal arrows} of a square.
 We denote by:
	\[\lv\un{h}{\square}\ \cD
\]
the set of squares with \emph{vertical} arrows (= source and target objects) in $\cD_v$ and \emph{horizontal} arrows in $\cD_h$. This defines a subcategory of $\sq\cD$, and actually a \emph{double subcategory} (i.e. stable by both horizontal and vertical compositions).

Taking $h=i\,/\,s$ and $v=\text{blank}$, or $v=i\,/\,s$ and $h=\text{blank}$, the corresponding squares will be called $hi\,/hs-$, or $vi\,/vs-$squares respectively.

Adding now an extra constraint of type (ipb), (gbp), or (spb),
one can prove that the following additional constraints still define \emph{double subcategories} of $\sq\cD$:

	\[\lv\un{h}{\ipb}\ \cD,\h
	\lv\un{h}{\gpb}\ \cD,\h
	\ls\un{h}{\spb}\ \cD,\h
	\lv\un{s}{\spb}\ \cD\,.
\]

There are a lot of useful diptych structures on $\sq\cD$ and subcategories.

\subsubsection{Naïve diptych structure for squares}
\label{nasq}

First of all we get a rather obvious \emph{naïve diptych structure} on $\sq \cD$ (which will not be the most useful one in the sequel) by taking as good monos\,/\,epis the $hi\,/\,hs-$squares (\ref{notsq}), i.e.:
	\[\un{\ssst{\text{naïve}}}{\,\square}\sD\underset{\text{def}}{=}(\square\ \cD;\ \un{i}{\square}\ \cD,\ \un{s}{\square}\ \cD)\,.
\]

One has two obvious forgetful functors (projections onto the top\,/\,bottom arrow of the square):
	\[\pi^{\text{top}}\,/\,\pi^{\text{bot}}:\un{\ssst{\text{naïve}}}{\,\square}\sD\rightarrow\sD\,,\tag{top\,/\,bot}
\]
which are \emph{well exact} (\ref{dipmor}) dip-morphisms (but \emph{not} faithful).

We have also a well exact faithful dip-embedding $\sD\stackrel{\iota}{\rightarrow}\un{\ssst{\text{naïve}}}{\square}\sD$, identifying an object with its unit and an arrow with a vertically degenerate square (in which this arrow is horizontally doubled).

\subsubsection{Main and canonical diptych structures for squares}
\label{masq}

For the sequel, it will be now necessary to be able to \emph{add vertical constraints}, needed for studying what will be called dip-diagram morphisms in $\sD$.

Relying on axiom (iv)(a), we get first a \emph{full subdiptych}, corresponding to the vertical inclusion $\cD_i\subset\cD$ :
	\[\un{\ssst{\text{naïve}}}{\li\sq}\sD=(\li\square\ \cD;\ \li\un{i}{\square}\ \cD,\ \li\un{s}{\square}\ \cD).\tag{i-sq.naïve-dip}
\]

We might also consider the naïve PRE-diptych corresponding to the vertical inclusion $\cD_s\subset\cD$:
	\[\un{\ssst{\text{naïve}}}{\ls\sq}\sD=(\ls\square\ \cD;\ \ls\un{i}{\square}\ \cD,\ \ls\un{s}{\square}\ \cD).\tag{s-sq.naïve-predip}
\]

But it turns out (due to the important dissymmetry of axiom (iv)) that looking for a needed \emph{diptych} structure on the full subcategory $\ls\square\ \cD$
\emph{demands a stricter choice for good epis} 
in order to construct the required fibred products for the vertical $s-$arrows regarded as \emph{objects} of $\ls\square\ \cD$; for this purpose, it is suitable to strengthen the condition on $hs-$squares, \emph{using the} spb- or gpb-\emph{squares} just introduced for that purpose (relying on \ref{remsq}).

Actually, for the sequel, the \emph{most important diptych structure} on squares (called \emph{main}) will be:
	\[\square\,\hat{}\ \sD=(\sq\cD;\ \un{i}{\square}\ \cD,\ \un{s}{\spb}\ \cD),\tag{sq.spb-dip}
\]
(from now on the \emph{convention for superscripts} on the left hand side will be that \emph{left\,/\,right} super\-scripts describe \emph{constraints for defining} the good \emph{monos\,/\,epis}). Note that, since (after \ref{remsq}) one has $\un{i}{\square}\ \cD=\,\un{i}{\ipb}\ \cD$, one might write as well, \emph{equivalently but more symmetrically}, the \emph{main diptych structure for squares} as:
	\[\un{\ssst{\text{main}}}{\,\square}\sD
	\un{\text{def}}{=}{\check{}\,}\square\,\hat{}\ \,\sD=
	(\sq\cD;\ \un{i}{\ipb}\ \cD,\ \un{s}{\spb}\ \cD)
	\,.\tag{sq.main-dip}
\]
One proves this gives rise now to two full \emph{sub\-diptychs} (adding vertical constraints):
	\[\un{\ssst{\text{main}}}{\li\sq}\,\sD=(\li\square\ \cD;\ \li\un{i}{\square}\ \cD,\ \li\un{s}{\spb}\ \cD)
	\tag{i-sq.main-dip}
\]
	\[\un{\ssst{\text{main}}}{\ls\sq}\,\sD=(\ls\square\ \cD;\ \ls\un{i}{\square}\ \cD,\ \ls\un{s}{\spb}\ \cD)\,.\tag{s-sq.main-dip}
\]

One may also \emph{replace} the spb condition for \emph{good epis} by the \emph{stronger} condition gpb, thus getting:
	\[\un{\ssst{\text{maiN}}}{\,\square}\sD
	\un{\text{def}}{=}
	{\square\ov{\diamond}{\phantom a}\ }\sD
	\un{\text{def}}{=}
	(\sq\cD;\ \un{i}{\square}\ \cD,\ \un{s}{\gpb}\ \cD)
	\tag{sq.maiN-dip}
\]
and analogously the full subdiptychs:
	\[\un{\ssst{\text{maiN}}}{\li\sq}\sD
	\un{\text{def}}{=}
	{\li\square\ov{\diamond}{\phantom a}\ }\sD
	\un{\text{def}}{=}(\li\square\ \cD;\ \li\un{i}{\square}\ \cD,\ \li\un{s}{\gpb}\ \cD)\tag{i-sq.maiN-dip}
\]
	\[\un{\ssst{\text{maiN}}}{\ls\sq}\sD
	\un{\text{def}}{=}
	{\ls\square\ov{\diamond}{\phantom a}\ }\sD
	\un{\text{def}}{=}
	(\ls\square\ \cD;\ \ls\un{i}{\square}\ \cD,\ \ls\un{s}{\gpb}\ \cD)\tag{s-sq.maiN-dip}.
\]

One can \emph{strengthen also} the conditions for \emph{good monos}, demanding them to be spb, or equivalently 
gpb, and get various combinations for dip-pairs, such as, for instance (with transparent conventions about upper case letters!):
	\[\un{\ssst{\text{MaiN}}}{\,\square}\sD
	\un{\text{def}}{=}
	{\ov{\diamond}{\phantom a}\square\ov{\diamond}{\phantom a}\ }\sD=(\sq\cD;\ \un{i}{\gpb}\ \cD,\ \un{s}{\gpb}\ \cD)\tag{sq.MaiN-dip}
\]

Finally these diptychs \emph{induce} diptych structures on certain \emph{subcategories} of $\sq\cD$, and one gets notably the important subdiptychs:
	\[\spb\ \sD
	=(\,\spb\,\cD;\ \un{i}{\spb}\,\cD,\ \un{s}{\spb}\,\cD)
	=(\,\spb\,\cD;\ \un{i}{\gpb}\,\cD,\ \un{s}{\spb}\,\cD)
	\,,\tag{spb-dip}
\]
	\[\gpb\ \sD=(\,\gpb\,\cD;\ \un{i}{\gpb}\,\cD,\ \un{s}{\gpb}\,\cD)
	\,,\tag{gpb-dip}
\]
and the full subdiptychs (s-spb-dip) and (s-gpb-dip) below:
	\[\ls\spb\ \sD=(\,\ls\spb\,\cD;\ \ls\un{i}{\gpb}\,\cD,\ \ls\un{s}{\spb}\,\cD\,),
\H
\ls\gpb\ \sD=(\,\ls\gpb\,\cD;\ \ls\un{i}{\gpb}\,\cD,\ \ls\un{s}{\gpb}\,\cD).
\]
One observes that the above diptych structures on $\spb\,\cD$ and $\gpb\,\cD$ and their subcategories are \emph{induced both by the naïve and main} diptych structures: they will be called the \emph{canonical diptych structures}, and denoted by:
	\[\un{\ssst{\text{can}}}{\,\spb}\,\sD,\h
	\un{\ssst{\text{can}}}{\,\gpb}\,\sD
	\,.
\]

Note all these new square symbols define endo-functors of \textbf{Dip} (\ref{dipmor}). This allows to summarize the above constructions of diptych structures for categories of squares, using more condensed symbolic writings (omitting  $\sD$), such as:

\begin{eqnarray} \notag
	\boxed{
	\un{\ssst{\text{naïve}}}{\,\square\,}
	\un{\ssst\text{def}}{=}
	(\square\,;\ \un{i}{\square}\,,\ \un{s}{\square})
	}&&&
	\boxed{
	\un{\ssst{\text{main}}}{\,\square\,}
	\un{\ssst\text{def}}{=}
	\square\,\hat{}
	=(\square\,;\ \un{i}{\square}\,,\ \un{s}{\spb})\ 
	(={\check{}\,}\square\,\hat{}\,)}
	\\ \notag
	\boxed{
	\un{\ssst{\text{maiN}}}{\,\square\,}
	\un{\ssst\text{def}}{=}
	\square\ov{\diamond}{\phantom a}
	=(\square\,;\ \un{i}{\square}\,,\ \un{s}{\gpb})}
	&&&
	\boxed{
	\un{\ssst{\text{MaiN}}}{\,\square\,}
	\un{\ssst\text{def}}{=}
	\ov{\diamond}{\phantom a}\square\ov{\diamond}{\phantom a}=
	(\square\,;\ \un{i}{\gpb}\,,\ \un{s}{\gpb})
	}\\ \notag
	\boxed{
	\un{\ssst{\text{can}}}{\,\spb\,}
	\un{\ssst\text{def}}{=}
	(\spb\,;\ \un{i}{\spb}\,,\ \un{s}{\spb})
	}&&&
	\boxed{
	\un{\ssst{\text{can}}}{\,\gpb\,}
	\un{\ssst\text{def}}{=}
	(\gpb\,;\ \un{i}{\gpb}\,,\ \un{s}{\gpb})}\,.
\end{eqnarray}

\subsection{Second order diptych data: weak and strong exactness}
\label{seco}
Given a \emph{diptych} structure $\sD$ on a category $\cD$, the \emph{extra choice} of a diptych structure on $\sq\cD$ may be interpreted as ``\emph{second order diptych data}'' on $\cD$ and denoted by $\bbD$. The choices described in the previous section are then denoted shortly by:
	\[\boxed{
	\un{\ssst\text{naïve}}{\bbD},\ \un{\ssst\text{main}}{\bbD},\ \un{\ssst\text{maiN}}{\bbD},\ \un{\ssst\text{MaiN}}{\bbD}}
	\,.
\]

Given a functor $f:\cD\rightarrow\cD'$, one then has the notion of \emph{second order diptych morphism}: $\boxed{f:\bbD\rightarrow\bbD'}$\,, 
by which we mean that $f:\sD\rightarrow\sD'$ is a dip-morphism and moreover $\sq f:\sq\cD\rightarrow\sq\cD'$ is a dip-morphism too, with respect to these extra choices. This defines a new category $\pmb{\bbD}\textbf{ip}$.

Then the condition for a diptych morphism of being \emph{exact} in the sense defined in \ref{dipmor} (preservation of special pullbacks) may be expressed by saying 
$\boxed{f:\un{\ssst\text{maiN}}{\bbD}\rightarrow\un{\ssst\text{maiN}}{\bbD'}}$ 
is a second order diptych morphism, but this allows to introduce a variety of new notions, by combining various choices for $\bbD$ and $\bbD'$.

One observes that \emph{any} dip-morphism $f:\sD\rightarrow\sD'$ gives rise to a second order dip-morphism 
$\boxed{f:\un{\ssst\text{naïve}}{\bbD}\rightarrow\un{\ssst\text{naïve}}{\bbD'}}$ (and this remains valid 
for \emph{any} of the four previous choices for $\bbD$, keeping the naïve choice for $\bbD'$). 
This defines the (full and faithful) \emph{naïve embedding $\boxed{\mathbf{Dip}\subset\pmb{\bbD}\mathbf{ip}}$}\,.

One notes also that a second order dip-morphism 
$f:\un{\ssst\text{maiN}}{\bbD}\rightarrow\un{\ssst\text{maiN}}{\bbD'}$ 
\emph{implies} a second order dip-morphism 
$\boxed{f:\un{\ssst\text{main}}{\bbD}\rightarrow\un{\ssst\text{main}}{\bbD'}}$ (or equivalently $f:\un{\ssst\text{maiN}}{\bbD}\rightarrow\un{\ssst\text{main}}{\bbD'}$): this weaker latter notion might be called \emph{weak exactness}. There is also a stronger notion when $\boxed{f:\un{\ssst\text{MaiN}}{\bbD}\rightarrow\un{\ssst\text{MaiN}}{\bbD'}}$ is a dip-morphism, which may be called \emph{strong exactness}.

\subsection{Triangularly degenerate squares}

\subsubsection{Definition of s-objects}
\label{sob}

Let $\sov{}\cD$ denote the subcategory of $\ls\sq\cD$ for which the lower horizontal arrow is a unit, with its canonical embedding:
	\[\iota^{\text{top}}:\sov{}\cD\rightarrow\ls\sq\cD\,.
\]
Its objects, which are the arrows of $\cD_s$, will be called $s-$objects. It is endowed with a \emph{canonical diptych structure}, denoted by $\sov{}\sD$, which is \emph{induced both by the naïve and the main diptych structures} on $\ls\sq\cD$.

The (non-full) embedding $\iota^{\text{top}}$ is \emph{well exact} for both the naïve and the main diptych structures on $\ls\sq\cD$.

\subsubsection{s-localization, s-augmentation.}
\label{loc}

Now let us \emph{fix} an object $Q$ of $\cD$ and consider the full subcategory $\sov{Q}\cD$ of ``\emph{$s-$objects over $Q$}'', (i.e. with target $Q$), or briefly ``\,$s$-$Q$-objects\,''. 

It can be checked that this is an \emph{s-augmented diptych} in the sense of \ref{remax}~ (2), called the ``\,\emph{s-localization of} $\sD$ \emph{over} $Q\,$''. But the existence of products (which here are \emph{not} obtained from products of squares!) derives from the \emph{fibred} products over $Q$ (which becomes a \emph{terminal object} in the new category).

Precisely the canonical embedding into $\ls\sq\cD$ is dip-exact for both naïve and main diptych structures, but not well exact (it does \emph{not} preserve products).

We have an obvious faithful forgetful functor (forgetting the terminal object $Q$):
	\[\ptop{Q}:\sov{Q}\sD\rightarrow\sD\,,\tag{top}
\]
which is \emph{dip-exact}, but not (in general) well exact. There is now (in general) \emph{no embedding of} $\sD$ into $\sov{Q}\sD$.

However, in the \emph{special case} when the given diptych $\sD$ is already \emph{s-augmented} (\ref{remax}~ (2)), \emph{then} one gets \emph{canonical diptych isomorphisms} by means of:

	\[\boxed{
\bdi[w=3em,labelstyle=\ssst,tight,midshaft]
\sD&\rTO~{\approx}^{\itop{\bullet}\ }&\sov{\bullet}\sD&\rTO~{\approx}^{\ptop{\bullet}\ }&\sD\\
\edi}\ .\tag{s-augm.dip}
\]

More generally the objects of $\sov{\bullet}\sD$ may be identified with the s-condensed objects of \ref{remax}~ (2), and $\ptop{\bullet}$ defines a well exact full dip-embedding.

\subsubsection{The case of $\protect\mathbf{\Finv}$.}
\label{Finv}

Let us take $\sD=\Finv$ (\ref{sym}) (which is not augmented, but possesses an \emph{initial} object \textbf{1}).
Here the $s$-\textbf{1}-objects are just the \emph{pointed} (or marked) objects (with a distinguished element).

Moreover there are two \emph{canonical pointings}, taking as the distinguished element of $\textbf{n}=[n-1]$ the minimum one 0 (or alternatively the maximum one $n-1$). It is easy to check that this defines a \emph{well exact} (but non-full) \emph{diptych embedding} (``\emph{adding zero}'' denoted by~ $\odot\,$, or adding a maximum) of the $s$-augmented diptych $\dot{{\Finv\ }}$ (\ref{sym}) (which, as was just observed, may be \emph{identified} with $\sov{\textbf{0}}\dot{\Finv\ }$) into $\sov{\textbf{1}}\Finv$, while the (faithful non-full) forgetful functor (``\emph{forgetting}'' the zero pointing!) $\ptop{\textbf{1}}$ is not an isomorphism. Summarizing (and adding at each end well exact full diptych embeddings), one has:
\[\boxed{\Finv\subset\dot{\Finv\ }
(\approx\sov{\textbf{0}}\dot{\Finv\ })
\rTO^{\odot\ }\sov{\textbf{1}}\Finv\rTO^{\ptop{\textbf{1}}}\Finv
\subset\dot{\Finv\ }}\ .\tag{augm+shift}
\]
One can recover the \emph{shift} endo-functor $\Delta$ of $\Finv$ (\ref{canend}) or its extension to $\dot{\Finv\ }$ \emph{by composing three of these maps}.

The interest of the (non-full) shift embedding, exact though not well exact, derives from the fact that it enables to express the objects \textbf{n} of $\pmb{\Finv}$, not only as iterated \emph{products} of \textbf{1}, but, for $n\geq3$, as iterated \emph{fibred products} over \textbf{1} (deriving themselves from \emph{products} in $\sov{\textbf{1}}\pmb{\Finv}$). This will be basic for groupoids, since these will be defined as diagrams which preserve \emph{pullbacks}, but in general \emph{not products} (\ref{defgrou}).

\section{Diptych structures for diagram morphisms}
In the next section, our ``structured groupoids'' will be defined (using their nerves) as ``diagrams'' of a certain type in a ``structuring category'' $\cD$, more precisely a diptych $\sD$. We need some general preliminary study of such diagrams and morphisms between these diagrams. This study relies on our previous study of squares.
\subsection{Small diagrams in a (``structuring'') category: $\protect\overset{\bullet}{\cT}(\cD)$ and $\protect\overrightarrow{\cT}(\cD)$}

\subsubsection{Geometric notations and terminology for diagram morphisms.} \label{nat}

When we are given a \emph{small} category $\cT$ (``\emph{type of diagram}'') and a (possibly large) category $\cD$ (``\emph{structuring category}''), we denote\footnote{%
Here we give up the classical \cite{ML} notation $\textbf{Cat}(\cT,\cD)$ for $\overset{\bullet}{\cT}(\cD)$, and the \emph{exponential notation} $\cD^{\cT}=\text{Funct}(\cT,\cD)$ for $\overrightarrow{\cT}(\cD)$, and we slightly modify and simplify the notations of \cite{P04}, in order to make them more transparent and coherent.
} here by:
\begin{itemize}
	\item $\overset{\bullet}{\cT}(\cD)$ (\emph{objects}) the (large) set of \emph{functors} from $\cT$ to $\cD$, called here `` \emph{diagrams of type} $\cT\ $ in $\cD\,$'';
	\item $\overrightarrow{\cT}(\cD)$ (\emph{arrows}) the (large) category of \emph{natural transformations} \cite{ML} between these functors, (which means just \emph{commutative diagrams connecting two such diagrams}), called here ``\emph{diagram morphisms of type} $\cT$ in $\cD\,$''.
\end{itemize}

For picturing and interpreting these diagram  morphisms in various ways, we refer to the comments in \ref{transpo}, \ref{prod} and \ref{exnat} below, and to the figure hereafter (Fig. \ref{figex}), which displays an \emph{example of diagram morphism} (only generators of $\cT$ are figured, and the ``connecting arrows'' are dashed, for better legibility; a diagram morphism is a \emph{natural transformation}, but may also be regarded as a \emph{functor} between various categories, or as a \emph{double functor} between certain double categories).

\begin{figure}[ht]
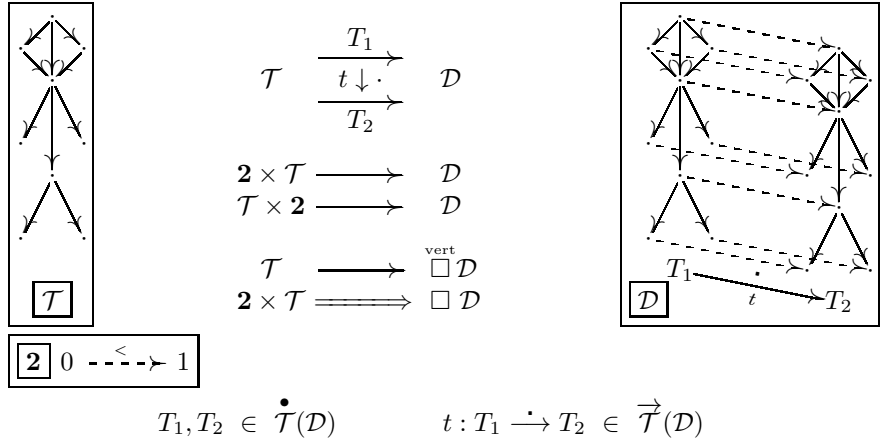

\begin{flushleft}
\H\boxed{
\bdi[h=1.2em,w=1.2em,tight,abut,inline]
&\cdot&\\
\cdot\ldTO(1,1)&\dTO&\cdot\rdTO(1,1)\\
&\cdot\rdTO(1,1)\ldTO(1,1)&\\
\ldTO(1,2)&\dTO&\rdTO(1,2)&\\
\cdot&&\cdot\\
&\cdot&\\
\ldTO(1,2)&&\rdTO(1,2)\\
\cdot&&\cdot\\
&&\\
&\boxed{\cT}&\\
\edi
}
\H\H\h
\bdi[w=1.7em,h=1.2em,tight,inline]
&&&&\\
&&&&\\
\h\cT\h&&
\pile
{\rTO^{T_1}\\
{t\downarrow\cdot}\\
\rTO_{T_2}}
&&\h\cD\h\\
&&&&\\
&&&&\\
\textbf{2}\times\cT&&\rTO&&\h\cD\h\\
\cT\times\textbf{2}&&\rTO&&\h\cD\h\\
&&&&\\
\h\cT\h&&\rTO&&\ \ov{\ssst\text{vert}}{\square}\cD\ \\
\textbf{2}\times\cT&&\rDoubleTo&&\ \sq\cD\\
\edi
\H\H\h
\boxed{
\bdi[h=1.2em,w=1.2em,tight,abut,inline,labelstyle=\sst]
&\cdot&&&\\
\cdot\ldTO(1,1)&\dTO&\cdot\rdTO(1,1)
\rdDashTO(5,1)&&&
&\cdot&\\
&\cdot\rdTO(1,1)\ldTO(1,1)\rdDashTO(5,1)&
&\rdDashTO(5,1)&&
\cdot\ldTO(1,1)&\dTO&\cdot\rdTO(1,1)\\
\ldTO(1,2)&\dTO&\rdTO(1,2)
\rdDashTO(5,1)&&&
&\cdot\rdTO(1,1)\ldTO(1,1)&\\
\cdot&&\cdot
&&&
\ldTO(1,2)&\dTO&\rdTO(1,2)\\
&\cdot\rdDashTO(5,1)&
&\rdDashTO(5,1)&&
\cdot&&\cdot\\
\ldTO(1,2)&&\rdTO(1,2)
\rdDashTO(5,1)&&&
&\cdot&\\
\cdot&&\cdot
&&&
\ldTO(1,2)&&\rdTO(1,2)\\
&T_1\rdDashTO(5,1)&
&\rdDashTO(5,1)&&
\cdot&&\cdot\\
\boxed{\cD}&&
\rdTO(5,1)^{\centerdot\ }_{t\ }&&&
&T_2&\\
\edi
}\\
\end{flushleft}
\begin{flushleft}
\H\boxed{
\bdi[h=1.1em,w=1.1em,tight,labelstyle=\ssst]
\boxed{\textbf{2}}\ &0&&\rDashTO^{<\ }&&1{}\\
\edi}
\end{flushleft}
\begin{center}
$T_1,T_2\ \in\ \overset{\bullet}{\cT}(\cD)$\H\H
$t:T_1\stackrel{\centerdot}{\longrightarrow}T_2\ \in\ \overrightarrow{\cT}(\cD)$
\end{center}
\caption{
Scheme of a ``diagram morphism'' $t$ 
between ``diagrams $T_1$, $T_2$ of type $\cT\,$'' in the ``structuring'' category $\cD$ (\ref{nat}, \ref{prod}, \ref{transpo}, \ref{exnat}).
}\label{figex}
\end{figure}

We note that the arrows of $\overrightarrow{\cT}(\cD)$ are defined by \emph{a collection of commutative squares} (indexed by $\cT$): 
this is expressed more precisely by the relation (Diag-switch) below (\ref{transpo}). This remark will be exploited in \ref{subcat} below in order to define various subcategories of $\overrightarrow{\cT}(\cD)\,$ deriving from the study made in \ref{dipsq}.

The notations are coined in order to emphasize the \emph{functorial} behaviour of the symbols $\overset{\bullet}{\cT}$ and $\overrightarrow{\cT}$ with respect to the objects $\cD$ (and to the arrows $\cD\rightarrow\cD'$)\footnote{
Being more precise and explicit about such functors would require some discussion of ``sizes'' of the universes involved, which would be here pedantic and irrelevant for our purpose, since we just wish to introduce \emph{convenient notations} for [morphisms between] diagrams of various types, in order to make transparent the twofold functorial behaviours.}.

\subsubsection{Iterated diagrams and interchange law}
\label{iter}
Now we can \emph{compose} such functorial symbols, supposing we are given a second small category $\cS$.

We rewrite, with the previous language and notations, some well known canonical isomorphisms, which yield a \emph{remarkable relation} between diagrams of types $\cS$, $\cT$, and $\cS\times\cT$, entailing a kind of commutativity (or \emph{interchange}) isomorphism $\circlearrowleft$. One has, for \emph{arrows}, i.e. diagram \emph{morphisms} (omitting $\cD$), the canonical isomorphisms:

\[\boxed{\overrightarrow{\cS}\circ\overrightarrow{\cT}=\overrightarrow{\cS\times\cT}
\ttt\overrightarrow{\cT\times\cS}=\overrightarrow{\cT}\circ\overrightarrow{\cS}}\,,\tag{Int.ch}
	\]
hence also, at the \emph{objects} (i.e. \emph{diagrams}) level: 
	\[\boxed{\overset{\bullet}{\cS}\circ\overrightarrow{\cT}=\Big(\cS\times\cT\Big)^{\bullet}
	\ttt\Big(\cT\times\cS\Big)^{\bullet}=\overset{\bullet}{\cT}\circ\overrightarrow{\cS}}\tag{int.ch}
\]
(related to \textsc{Ehresmann}'s \emph{interchange law} (``loi d'échange'') for double categories or groupoids).

Taking $\cS=\cT$ and iterating our functors $n$ times, we get:
	\[\left(\overrightarrow{\cT}\right)^{n}=\overrightarrow{\cT^n}\h\text{ and }\h
	\Big(\cT^{n+1}\Big)^{\bullet}=\overset{\bullet}{\cT}\circ\overrightarrow{\cT^{n}}
	\longleftrightarrow\Big(\cT^{n}\Big)^{\bullet}\circ\ovra{\cT}.
\]

\subsubsection{Picturing a product category}
\label{prod}
It is convenient to picture an arrow $w=(f,g)$ of the \emph{product category} $\cS\times\cT$ (with $f:a\rightarrow b,\ g:c\rightarrow d$), by the diagonal of a (commutative) \emph{square}, thought \emph{in a vertical plane}, with two directions, called \emph{longitude} (dashed) and \emph{latitude} (dotted)\footnote{
\textsc{Ehresmann}'s geometrical terminology.}
 (see \ref{transpo}):
\begin{equation}
	\bdi[h=3em,w=5em,tight,labelstyle=\ssst,objectstyle=\scriptstyle,inline,abut]
(a,c)&\rDashTO~{(f,1_c)}&(b,c)\\
\dTO~{(1_a,g)}[dotted]&\rdTO~{(f,g)} &\dTO~{(1_b,g)}[dotted]\\%
(a,d)&\rDashTO~{(f,1_d)}&(b,d)\\
\edi\tag{prod}
\H.
\end{equation}
This square owns the very special property of being \emph{uniquely determined by any pair of adjacent edges}, and also by its diagonal. In the absence of terminal objects, it may not be a pull back, though its two components are both trivial pullbacks.

Though this is a \emph{purely graphic convention}, it will be natural to call $\circlearrowleft$ the interchange or ``switching'' or \emph{transposition} isomorphism (in the matrix sense, i.e. a \emph{diagonal symmetry}). 

One can observe too that the product category $\cS\times\cT$ owns also a (very special) natural \emph{double category} structure (in \textsc{Ehresmann}'s sense\footnote{
This means two category composition laws satisfying the interchange law \cite{ML}.
}
), transparent through the above representation, defined by latitudinal and longitudinal compositions. Each arrow $w=(f,g)$ of $\cS\times\cT$ admits (in a twofold way, according to the orders of the factors) a unique factorization through the (dotted or dashed) subcategories $\cS_0\times\cT$ and $\cS\times\cT_0$ (which are the bases of the two category structures); such an arrow $w$ may be \emph{identified} with a \emph{double functor} (i.e. respecting both composition laws):
	\[\ddot{w}:\pmb{\square}=\textbf{2}\times\textbf{2}\Longrightarrow\cS\times\cT\,.
\]
When both $\cS\times\cT$ and $\sq\cD$ are regarded as \emph{double categories}, any (simple) functor $V:\cS\times\cT\rightarrow\cD$ may be regarded as defining a \emph{double functor}:
\begin{eqnarray}
	\ddot{V}:\cS\times\cT\Longrightarrow\sq\cD, \h\ddot{w}\longmapsto V\circ\ddot{w}:
	\pmb{\square}\Rightarrow\cD\tag{double}\ .\\
\notag
\end{eqnarray}

\subsubsection{Example: path functors.}\label{exnat}
When $\cT$ is a finite ordinal $\mathbf{n}=[n-1]$ (\ref{preoset})\footnote{
We remind that, with our notations, $n$ is the \emph{cardinal} of the \textbf{ordinal} $\textbf{n}=[n-1]$. This shift is traditional and unavoidable when dealing with simplicial objects \cite{ML}.
},
	\[\cD_{[n-1]}=\cD_{\textbf{n}}\underset{\text{def}}{=}\overset{\bullet}{\mathbf{n}}\,(\cD)
	\tag{path}
\]
is just the set of \emph{paths} (or strings) in $\cD$ 
with \emph{length} $n-1$ and $n$ \textbf{vertices}, while $\overrightarrow{\mathbf{n}}(\cD$) consists of \emph{morphisms} between two such paths. This may be written functorially:
	\[\boxed{(\,)_{[n-1]}=(\,)_{\textbf{n}}=\overset{\bullet}{\mathbf{n}}}\tag{nerve}
\]

\textbf{}
One has, by (int.ch) and (Int.ch):
\begin{eqnarray}
	\left(\overrightarrow{\textbf{q}}\right)_{[p-1]}=\left(\overrightarrow{\textbf{q}}\right)_{\textbf{p}}=
\ov{\bullet}{\textbf{p}}\circ\ovra{\textbf{q}}
	&=(\textbf{p}\times\textbf{q})^{\bullet}&&\ttt&
	(\textbf{q}\times\textbf{p})^{\bullet}
	&=\ov{\bullet}{\textbf{q}}\circ\ovra{\textbf{p}}=\left(\overrightarrow{\textbf{p}}\right)_{\textbf{q}}=
	\left(\overrightarrow{\textbf{p}}\right)_{[q-1]}\tag{double nerve}\\
	\overrightarrow{\textbf{p}}\circ\overrightarrow{\textbf{q}}&=\ 
\overrightarrow{\textbf{p}\times\textbf{q}}&&\ttt&
\overrightarrow{\textbf{q}\times\textbf{p}}\ \
&=\overrightarrow{\textbf{q}}\circ\overrightarrow{\textbf{p}}\notag.
\end{eqnarray}
It may be convenient to consider more generally \emph{paths of diagrams}, setting:
	\[\boxed{\Big(\overrightarrow{\cT}\Big)_{[n-1]}
	\underset{\text{def}}{=}\overset{\bullet}{\mathbf{n}}\circ\overrightarrow{\cT}=
	\Big(\mathbf{n}\times\cT\Big)^{\bullet}
	\,\ttt\,\Big(\cT\times\mathbf{n}\Big)^{\bullet}=
	\overset{\bullet}{\cT}\circ\overrightarrow{\mathbf{n}}}\,.\tag{Path-Diag}
\]

These paths define the nerve of the category $\ovra{\cT}(\cD)$, while the double nerve is associated to the double category $\sq\cD$. Disregarding the questions of ``sizes'', the previous formula might be interpreted as attaching to any small category $\cT$ (called here ``type of diagrams'') a \emph{simplicial object} \cite{ML} $\mathbf{n}\mapsto\Big(\overrightarrow{\cT}\Big)_{[n-1]}$ in the ``category of all categories'' \cite{ML}. One has by construction:
	\[\Big(\ovra{\cT}(\cD)\Big)_{[n-1]}=\Big(\overrightarrow{\cT}\Big)_{[n-1]}(\cD)
	=\Big(\mathbf{n}\times\cT\Big)^{\bullet}(\cD)
	\,\ttt\,\ov{\bullet}{\cT}\Big(\overrightarrow{\textbf{n}}(\cD)\Big).
\]

\emph{Special cases} ($\textbf{n}=\textbf{1},\textbf{2}$) of these ``\emph{path} functors'' are the \emph{forgetful} ``\emph{objects}'' and ``\emph{arrows}'' functors (which \emph{forget} a part of the categorical data $\cD\,$: the arrows for the former, the composition law for the latter). One has (still omitting $\cD$):
	\[\boxed{\overset{\bullet}{\mathbf{1}}=\overset{\bullet}{\pmb{\star}}=(\,)_{[0]},\h
	\overrightarrow{\mathbf{1}}=\overrightarrow{\pmb{\star}}=\text{Id},\h
	\overset{\bullet}{\mathbf{2}}=\,\overset{\bullet}{\pmb{\downarrow}}\,=(\,)_{[1]},\h
\overrightarrow{\mathbf{2}}=\overrightarrow{\pmb{\downarrow}}\underset{\text{def}}{=}
\un{\ssst{\text{hor}}}{\square}}\ ,
\tag{obj, arr, sq}
\]
where notations $\mathbf{1}=\pmb{\star}$ and $\mathbf{2}=\,\pmb{\downarrow}$ of \ref{preoset} are used, and $\un{\ssst{\text{hor}}}{\square}\ \cD$\footnote{
\textsc{Ehresmann}'s suggestive notation for what he called ``\emph{quatuors}'', i.e. quartets. The more classical (much less suggestive) exponential notation in \cite{ML} is $\cD^{\textbf{2}}$, called ``\emph{category of arrows of} $\cD$\,'' (which has not to be confused with $\cD^{2}=\cD\times\cD$ !).
} denotes, for more precision (\ref{notsq}), the category of \emph{commutative squares} 
in $\cD$ \emph{with the horizontal composition} (see Rem. \ref{transpo} below).
One has also more generally, specializing (int.ch), (Int.ch) and (Path-Diag):

	\[\boxed{\Big(\overrightarrow{\cT}\Big)_{[0]}=
	\overset{\bullet}{\cT},\H
\Big(\overrightarrow{\cT}\Big)_{[1]}=\,\Big(\mathbf{2}\times \cT\Big)^{\bullet}\,\ttt\,
	\Big(\cT\times \mathbf{2}\Big)^{\bullet}\,
=\overset{\bullet}{\cT}\circ\un{\ssst{\text{hor}}}{\square}}\ .\tag{diag, Diag}
\]
Using the remarks in \ref{prod}, this allows to regard also an arrow $t$ of $\ovra{\cT}(\cD)$ as a \emph{double functor} $\ddot{t}$, when regarding both $\textbf{2}\times\cT$ and $\sq\cD$ as \emph{double categories}:
	\[\ddot{t}:\textbf{2}\times\cT\rDoubleTo\ \sq\cD\,.
\]
\subsubsection{Remarks: diagram morphisms as diagrams of squares}
\label{transpo}\hfill
\begin{enumerate}
	\item Though the way of picturing the arrows of a category has \emph{no mathematical content}, it is unavoidable to make it precise for getting mental images, which are necessary for memorizing and phrasing intricate situations. When a single category is being considered, it will be often (not always) comfortable to think originally its arrows ``horizontally''. But as soon as a second category is involved, it is often necessary to be able to ``switch'' the directions of arrows, using pairs of terms such as horizontal\,/\,vertical or longitudinal\,/\,latitudinal (\textsc{Ehresmann}), \emph{which have to be interpreted in a relative rather than absolute sense}.

Note for instance that the composition law in $\ovra{\cT}(\cD)$ is called ``\emph{vertical}'' in \cite{ML} (though sometimes pictured horizontally!), but will be often better thought and written here \emph{horizontally}, the diagrams of type $\cT$ (which are the \emph{objects} of this category) having then to be pictured ``\emph{vertically}'' by means of their images in $\cD$ (indexed by~ $\cT$), and the connecting arrows ``\emph{horizontally}'' (\emph{see} Fig. \ref{figex} \emph{above}). However we stress that this implies ``\emph{switching}'' \emph{the direction of the diagrams when passing from} $\overset{\bullet}{\cT}(\cD)$ \emph{to} $\ovra{\cT}(\cD)$.

Notably the bijection (Diag): 
$\Big(\overrightarrow{\cT}\Big)_{[1]}\ttt
\overset{\bullet}{\cT}\circ\square$
implies a switching in the notations between horizontal and vertical diagrams. 

\begin{itemize}
	\item Precisely \emph{we shall agree} that, by ``\emph{abus de notation}'', $\,\overset{\bullet}{\cT}(\sq\cD)$ on the right hand side denotes functors from $\cT$ to the \emph{vertical} category of squares, and later these are \emph{composed horizontally}, using the alternative \emph{horizontal} composition of squares, in order to recover the composition law of the left hand side. 
\emph{These conventions will be implicit throughout whenever writing:}
	\[\boxed{\overrightarrow{\cT}\ \longleftrightarrow\ 
	\overset{\bullet}{\cT}\circ\,\sq}\ .\tag{Diag-switch}
\]
This will lead below (\ref{subcat}) to the \emph{convention}: $\ovra{\cT}=\ov{\ssst{\square}}{\cT}$.
\end{itemize}

\item 
The above formula (Diag) (\ref{exnat}) means the useful fact that (when \emph{disregarding the composition} of diagram morphisms) any \emph{morphism} (arrow) between diagrams of type $\cT$ in $\cD$ may be interpreted as a \emph{diagram} (object), \emph{in a twofold way}: either as a \emph{diagram} of new type $\cT\times\mathbf{2}\ (\approx\mathbf{2}\times\cT)$ \emph{in the same category} $\cD$, or alternatively as a \emph{diagram} of the \emph{same type} $\cT$ in the new\emph{ vertical category} $\ov{\ssst{\text{vert}}}{\square}\cD$ \emph{of squares of} $\cD$ (see above Fig.\ref{figex}).

The latter interpretation enables then to \emph{recover the composition law} by using the \emph{alternative} (horizontal) composition law of squares. This will be used systematically in the next section for \emph{transferring properties of squares to general diagrams of type} $\cT$.
\end{enumerate}

Taking $\cT=\mathbf{2}$, one recovers the description of commutative squares, when forgetting their two composition laws:
\[\big(\square\ \cD\big)_{[1]}\approx\Big(\mathbf{2}\times\mathbf{2}\Big)^{\bullet}(\cD)\,=\,\overset{\bullet}{\pmb{\square}}\ (\cD),
\text{ or briefly }\boxed{\big(\square\big)_{[1]}\approx\,\overset{\bullet}{\pmb{\square}}}\ ,\tag{sq}
\]
where $\square$ is interpreted in the left hand sides as the operator $\overrightarrow{\mathbf{2}}$ (\ref{exnat}), while the boldface $\pmb{\square}$ in the right hand sides is picturing the product category $\mathbf{2}\times\mathbf{2}$ (notation of \ref{preoset}). All these geometric notations are at once very precise and coherent, while yielding an easy handling of more intricate diagrams.

Note that we shall be led below (see also \ref{canend}, above), to assign to the square symbols some \emph{special meanings for groupoids}; the context should avoid confusion.

\subsubsection{Defining subcategories of the category $\protect\ovra{\cT}(\cD)$ of diagram morphisms}
\label{subcat}

The correspondence (Diag-switch) of the previous paragraph \ref{transpo} (1) allows to \emph{transfer} many definitions, introduced in section \ref{dipsq} about squares, \emph{to morphisms} between diagrams of general type $\cT$ (the case of squares corresponding to $\cT=\mathbf{2}$). We remind that this correspondence uses successively the vertical and horizontal composition laws of $\sq D$.

The following \emph{general scheme} may be used to define \emph{various subcategories} of $\ovra{\cT}(\cD$). 
We consider various symbols $\Xi$ assigning to $\cD$ (in a functorial way) \emph{double subcategories} $\Xi\,\cD$ of $\sq\cD$ (i.e. stable by both vertical and horizontal laws). This may be denoted symbolically by: 
$\boxed{\Xi\subset\sq}$\,.

One has to agree that, for $\Xi=\bullet\subset\sq$, one has $\bullet\,\cD=\cD$ ``written vertically'' (\emph{horizontally degenerate squares}, which become \emph{units for the horizontal law}).

Then the definition of the right hand side of (Diag-switch) (\ref{transpo}) (1) assigns, to each such symbol $\Xi$ (in a functorial way) a (special) \emph{subcategory} of the left hand side, denoted by:%
	\[\ovra{\cT}(\cD)\supset\ov{\ssst{\Xi}}{\cT}(\cD)\longleftrightarrow\ov{\bullet}{\cT}(\Xi\,\cD)\subset\ov{\bullet}{\cT}(\sq\cD)
\]
or symbolically (having to agree that $\boxed{\ov{\ssst{\square}}{\cT}=\ovra{\cT}}$\,, as announced in \ref{transpo}):
	\[\boxed{\ovra{\cT}\supset\ov{\ssst{\Xi}}{\cT}
	\un{\text{def}}{\longleftrightarrow}
	\ov{\bullet}{\cT}\circ\Xi\subset\ov{\bullet}{\cT}\circ\sq}\,.\tag{subsq-switch}
\]

\subsubsection{Application to prediptychs}
\label{appred}
Particularly, in the case when $\sD=(\cD;\cD_i,\cD_s)$ is a pre\-diptych, then, using notations of \ref{notsq}, the previous considerations apply notably 
 to the symbols:
	\[\Xi\,=\,\un{h}{\square}\,,\ {\un{h}{\ipb}}\,,\ {\un{h}{\gpb}}\,,\ {\un{h}{\spb}}\,,\lv\square,\,\lv\un{h}{\square},\cdots,
\]
where the symbols $h$, $v$ may take the values $i$, $s$, or blank.

Then we introduce the following \emph{shortened notations} (extending $\ov{\ssst{\square}}{\cT}=\ovra{\cT}$):
	\[\boxed{\un{h}{\ovra{\cT}}\un{\text{def}}{=}\ov{\un{h}{\square}}{\cT}}\,,\h
	\boxed{\check{\un{h}{\ovra{\cT}}}\un{\text{def}}{=}\ov{\un{h}\ipb}{\cT}}\,,\h
	\text{and analogously}\h
	\boxed{\ov{\diamond}{{\un{h}{\ovra{\cT}}}},\h\hat{\un{h}{\ovra{\cT}}}}\,,\tag{spec.sq}
\]
which allow to describe, in a functorial and transparent way, \emph{various styles of constraints and ways of combining them} (bearing either on the \emph{horizontal connecting arrows} or on some \emph{global} properties), \emph{imposed to the squares} connecting two diagrams of the same type $\cT$.

Symbols such as $\boxed{\check{\ovra{\cT}},\ \hat{\ovra{\cT}},\ \ov{\diamond}{{\ovra{\cT}}},\cdots}$
 have then to be regarded as defining \emph{functors 
with domain} \textbf{Predip} (\ref{dipmor}).

\subsection{(Pre)dip-diagrams in a (pre)diptych}
\label{diadip}

\subsubsection{Definition of $\protect\overset{\bullet}{\mathsf{T}}(\sD)$ and $\protect\overrightarrow{\mathsf{T}}(\sD)$}
\label{notsub}

$\sT$ is now a fixed small (\emph{pre\-})\emph{diptych} and $\sD$ a (possibly large) (\emph{pre\-})\emph{diptych}. We are going to describe some constructions attached to $\sT$ depending functorially on $\sD$. We remind that any categories $\cT$ and $\cD$ may (and will) be regarded as \emph{trivial} pre\-diptychs (\ref{defdip}).

First we introduce the \emph{subclasses}:
\[\boxed{\overset{\bullet}{\mathsf{T}}{_{\ssst\mathsf{ex}}}(\mathsf{D})\subset\overset{\bullet}{\mathsf{T}}(\mathsf{D})\subset\overset{\bullet}{\cT}(\cD)}
\tag{predipdiag, dipex}
\]
and \emph{subcategories}:
	\[\boxed{
	\ovra{\mathsf{T}}_{\ssst\mathsf{ex}}(\mathsf{D})\subset\ovra{\mathsf{T}}(\mathsf{D})\subset\ovra{\cT}(\cD)}\,,
	\tag{predipDiag, dipEx}
\]
where the symbols $\overset{\bullet}{\mathsf{T}}(\mathsf{D})$ and 
$\overset{\bullet}{\mathsf{T}}{_{\ssst\mathsf{ex}}}(\mathsf{D})$ 
mean respectively (\emph{pre})\-\emph{dip-morphisms} from $\sT$ to $\sD$ and \emph{dip-exact} morphisms (\ref{dipmor}) (in other words (\emph{pre\-})\emph{dip-diagrams} and \emph{exact diagrams of type} $\sT$ in $\sD$), and the over arrows denote the \emph{full subcategories generated} by such objects ((\emph{pre\-})\emph{dip-diagram morphisms}).

More generally, following the ideas of \ref{seco} about the interpretation of dip-exact\-ness, one might introduce also \emph{second order diptych data} $\bbT$ and $\bbD$ on $\cT$ and $\cD$ (which would have  been either the ``naïve'' or the ``maiN'' ones just above), then define analogously 
	\[\boxed{\ov{\bullet}{\bbT}(\bbD)\subset\overset{\bullet}{\mathsf{T}}(\mathsf{D})
	\subset\overset{\bullet}{\cT}(\cD)}
\]
as the subclasses of first and second order diptych morphisms, and finally consider the full subcategories thus generated:
	\[\boxed{\ovra{\bbT}(\bbD)\subset\ovra{\mathsf{T}}(\mathsf{D})\subset\ovra{\cT}(\cD)}\,.
\]
The symbols $\ovra{\mathsf{T}}$ or $\ovra{\bbT}$ now give rise to functors defined on the categories 
\textbf{Predip} (\ref{dipmor}) or $\bbD\textbf{ip}$ (\ref{seco}) as their domains.

Next we are looking for defining \emph{endo-functors} of (\textbf{Pre})\textbf{dip}, which requires \emph{defining (pre)\-diptych structures} on the categories of diagrams arising in (predipDiag, dipEx) just above.

\subsection{Naïve and main diptych structures on categories of diagrams}
\subsubsection{General method}
\label{gemeth}
The \emph{outline of the general process} for constructing new (pre)\-diptychs, attached to the small (pre)diptych $\sT$, and to the (pre)diptych $\sD$, will be the following (using the scheme \ref{subcat} for notations).

We consider a double set of symbols $\boxed{[\Xi]=\ov{\ssst{\text{vert}}}{\un{\ssst{\text{hor}}}{\Xi}}=\Big(\!\lv{\un{h}{\Xi}}\Big)}$ 
(where the arguments h,\,v may take independently the values i,\,s,\,blank, and $\lv\un{h}{\Xi}=\lv\Xi\,\cap\,\un{h}{\Xi}\,$). It is understood that each such symbol attaches functorially to the (pre)\-diptych $\sD$ a \emph{double subcategory} $\lv{\un{h}{\Xi}}\,\sD\subset\sq\cD$ as in \ref{subcat}, and that this yields a ``\emph{double pre\-diptych}'' structure on $\Xi\,\sD$: this means that $\un{h}{\Xi}\,\sD$\,/\,$\lv\Xi\,\sD$ defines a horizontal\,/\,vertical pre\-diptych structure 
$\un{\ssst{\text{hor}}}{\Xi}\,\sD$\,/\,$\ov{\ssst{\text{vert}}}{\Xi}\,\sD$, on $\Xi\,\sD$ when h\,/\,v takes the values (blank;\ i,\,s). One may write symbolically (omitting $\sD$):
	\[\un{\ssst{\text{hor}}}{\Xi}
	=\left(\Xi;
	\,\un{i}{\Xi},
	\un{s}{\Xi}
	\right),\h
	\ov{\ssst{\text{vert}}}{\Xi}
	=\Big(\Xi;
	\,\li{\Xi},
	\ls{\Xi}
	\Big)\,.
\]

Then (taking into account successively first the \emph{vertical} diptych structures and the definitions of \ref{notsub}, second the \emph{horizontal law} for defining the composition, as explained in \ref{transpo}), we can define the triple:
\[\Big(\ov{\bullet}{\sT}\,(\,\ov{\ssst{\text{vert}}}{\Xi}\,\sD);
\ov{\bullet}{\sT}\,(\,\un{i}{\ov{\ssst{\text{vert}}}{\Xi}}\,\sD),
\ov{\bullet}{\sT}\,(\,\un{s}{\ov{\ssst{\text{vert}}}{\Xi}}\,\sD)
\Big)
\]
and this defines, by applying (subsq-switch) of \ref{subcat} , a triple of \emph{subcategories} of $\ovra{\sT}(\sD)\subset\ovra{\cT}(\cD)$, which is the announced pre\-diptych. 
With shortened notations analogous to (spec.sq) in \ref{appred},
 we shall denote this pre\-diptych by:
	\[\ov{[\ssst\Xi]}{\sT}(\sD)
	\un{\text{def}}{=}
	\Big(\ov{\ssst{\Xi}}{\sT}(\sD);
	\un{i}{\ov{\ssst{\Xi}}{\sT}}(\sD),
	\un{s}{\ov{\ssst{\Xi}}{\sT}}(\sD)
	\Big)\,.
\]
or symbolically:
	\[\ov{\ssst{[\Xi]}}{\sT}
	\,\un{\text{def}}{=}\,
	\left(\ov{\ssst{\Xi}}{\sT};
	\un{i}{\ov{\ssst{\Xi}}{\sT}},
	\un{s}{\ov{\ssst{\Xi}}{\sT}}\right)
	\un{\text{def}}{\longleftrightarrow}
	\left(\ov{\bullet}{\sT}\circ\ov{\ssst{\text{vert}}}{\Xi};
	\ov{\bullet}{\sT}\circ\un{i}{\ov{\ssst{\text{vert}}}{\Xi}},
	\ov{\bullet}{\sT}\circ\un{s}{\ov{\ssst{\text{vert}}}{\Xi}}\right)\,.
	\tag{predip-switch}
\]
This triple defines a symbol $\ov{\ssst{[\Xi]}}{\sT}$ (depending on the data
$[\Xi]=\ov{\ssst{\text{vert}}}{\un{\ssst{\text{hor}}}{\Xi}}$\,), which determines an \emph{endo-functor of} \textbf{Predip}, 
with the symbolic correspondence explained above:
	\[\boxed{\ov{\ssst{[\Xi]}}{\sT}
	\,\un{\text{def}}{\longleftrightarrow}\,
	\ov{\bullet}{\sT}\circ\ov{\ssst{\text{vert}}}{\un{\ssst{\text{hor}}}{\Xi}}}\,.
\]

Using this general scheme, we are now looking for various basic (\emph{pre})\-\emph{diptych} structures on certain \emph{subcategories} of $\ovra{\sT}(\sD)$ and $\ovra{\mathsf{T}}_{\ssst\mathsf{ex}}(\mathsf{D})$ (introduced in \ref{subcat}), canonically and functorially associated to the (\emph{pre})\-\emph{diptych data} on the categories $\cD$, $\cT$.

\subsubsection{Naïve \emph{(}pre\emph{)}diptych structure for diagram morphisms}
\label{naive}
Here $\cT$ denotes just a small category regarded as a \emph{trivial} pre\-diptych (\ref{defdip}), and $\ov{\ssst{\text{vert}}}{\un{\ssst{\text{hor}}}{\Xi}}=
\ov{\ssst{\text{vtriv}}}{\un{\ssst{\text{hnaïve}}}{\square}}$, by which we mean the \emph{horizontal naïve} diptych of \ref{nasq} and the \emph{vertical trivial} pre\-diptych.

Then the previous scheme defines the \emph{naïve pre\-diptych} structure on $\ovra{\cT}(\sD)$, which actually is still a \emph{diptych} when $\sD$ is. Explicitly:
	\[\un{\ssst{\text{naïve}}}{\ovra{\cT}}(\mathsf{D})\underset{\text{def}}{=}\Big(\ovra{\cT}(\cD);\un{i}{\ovra{\cT}}(\cD),\un{s}{\ovra{\cT}}(\cD)\Big),
\]
or more briefly and symbolically:
\[\boxed{\un{\ssst{\text{naïve}}}{\ovra{\cT}}\,\underset{\text{def}}{=}\,\Big(\ovra{\cT};\un{i}{\ovra{\cT}},\un{s}{\ovra{\cT}}\Big)}\,,\tag{Diag naïve (pre)dip}
\]
where (as defined in \ref{subcat} (spec.sq)) $\un{i}{\ovra{\cT}}(\cD)\,/\,\un{s}{\ovra{\cT}}(\cD)$ denotes the subcategory of those diagram morphisms for which all the connecting (``horizontal'') arrows are taken from $\cD_i\,/\,\cD_s$.

Taking $\cT=\mathbf{2}$, this is coherent with the notations (sq.naïve-dip) of \ref{nasq}:
	\[\un{\ssst{\text{naïve}}}{\square}\sD\underset{\text{def}}{=}\Big(\square\ \cD;\ \un{i}{\square}\ \cD,\ \un{s}{\square}\ \cD\Big).\tag{sq.naïve (pre)dip}
\]
The switch correspondence reads here:

	\[\boxed{\un{\ssst{\text{naïve}}}{\ovra{\cT}}
	\un{\text{def}}{\longleftrightarrow}
{\ov{\bullet}{\cT}}\circ
\ov{\ssst{\text{vtriv}}}{\un{\ssst{\text{hnaïve}}}{\square}}}\,.\tag{Diag naïve (pre)dip}
\]

This is an \emph{endo\-functor of the category} \textbf{Dip}, attached to any type of diagram~ $\cT$.

\subsubsection{Square description for dip-exact diagram morphisms}
\label{sqdes}

Here $\mathsf{T}=(\cT;\cT_i,\cT_s)$ is again \emph{also} a small \emph{pre\-diptych}, and $\mathsf{D}=(\cD;\cD_i,\cD_s)$ still a (possibly large) \emph{pre\-diptych}, and we consider, as in \ref{notsub}, predip-diagrams and morphisms:

	\[\overset{\bullet}{\mathsf{T}}(\mathsf{D})\subset\overset{\bullet}{\cT}(\cD)
	\h\text{and}\h
	\ovra{\mathsf{T}}(\mathsf{D})\subset\ovra{\cT}(\cD)\,.
\]

These symbols $\overset{\bullet}{\mathsf{T}}$ and $\ovra{\mathsf{T}}$ (attached to the \emph{pre\-diptych} $\sT$) are functorial with respect to the \emph{pre\-diptychs} $\sD$ [and predip-morphisms].

Still interpreting a diagram morphism $f$ as a functor from $\cT$ to the (\emph{vertical}!) category of squares (\ref{transpo}), the conditions for a diagram morphism being a \emph{predip-diagram} morphism may be written, with our notations:
	\[f(\cT_i)\subset\li\sq\cD,\H f(\cT_s)\subset\ls\sq\cD.
\]
This may be expressed by saying $f$ is still a \emph{pre\-diptych morphism }%
from the pre\-diptych $\sT$ to the (here simple, not double) \emph{naïve pre\-diptych} $\ov{\ssst{\text{vnaïve}}}{\square}\sD$ of \ref{nasq} (but written here \emph{vertically}; it was the trivial structure in the previous paragraph!). So we have now the symbolic correspondence (for functors on \emph{pre\-diptychs}):
	\[\boxed{\ovra{\sT}\longleftrightarrow
	\overset{\bullet}{\mathsf{T}}\circ\ov{\ssst{\text{vnaïve}}}{\sq}}\,.\tag{predip-Diag switch}
\]

Now if \emph{moreover} $\sT$ and $\sD$ are both \emph{diptychs}, then it is possible (and necessary for the applications of the next section) to consider the \emph{dip-exact diagrams of type~ $\sT$}, i.e. the \emph{exact}\footnote{
We remind this means (in the dip-sense) preserving special pullbacks(\ref{remax} (4)) (but in general not products, nor other limits or co-limits).
} 
diptych morphisms from $\sT$ to $\sD$ (\ref{dipmor}), and \emph{morphisms} between such exact diagrams. (Observe that this new constraint is of \emph{global} nature, meaning thereby that it involves several arrows of $\cT$ together). We remind the notations of \ref{notsub}:\[\overset{\bullet}{\mathsf{T}}{_{\ssst\mathsf{ex}}}(\mathsf{D})\subset\overset{\bullet}{\mathsf{T}}(\mathsf{D})\subset\overset{\bullet}{\cT}(\cD)
	\h\text{and}\h
	\ovra{\mathsf{T}}_{\ssst\mathsf{ex}}(\mathsf{D})\subset\ovra{\mathsf{T}}(\mathsf{D})\subset\ovra{\cT}(\cD)\ .\tag{dipex, dipEx}
\]

In this way we have now attached, to each \emph{small diptych} $\sT$, new sym\-bols $\overset{\bullet}{\mathsf{T}}{_{\ssst\mathsf{ex}}}$ and $\ovra{\mathsf{T}}_{\ssst\mathsf{ex}}$, still functorial with respect to \emph{diptychs} $\sD$ (more precisely associating to \emph{diptychs} certain \emph{classes} or \emph{categories}, in a functorial way). It may be proved (using cubical diagrams and diptych axioms) that the basic correspondence between diagram morphisms and diagrams of squares is \emph{still valid for exact diagrams}:
	\[\boxed{\ovra{\sT}_{\ssst\mathsf{ex}}\longleftrightarrow
	\ov{\bullet}{\sT}_{\ssst\mathsf{ex}}\circ\ov{\ssst{\text{vnaïve}}}{\sq}}\,.\tag{dipEx switch}
\]

\subsubsection{Main diptych structures for morphisms between exact diagrams}
\label{maindip}

We suppose here that \emph{both} $\sT$ and $\sD$ are still \emph{diptychs} (not only pre\-diptychs!) and are now looking for the announced \emph{diptych} structures on the categories $\ovra{\sT}(\sD)$ and $\ovra{\sT}_{\ssst\mathsf{ex}}(\sD)$.

The \emph{basic point} is that the \emph{naïve} diptych structure on $\ovra{\cT}(\cD)$, defined in \ref{naive}, \emph{induces} on $\ovra{\sT}(\sD)$ and $\ovra{\sT}_{\ssst\mathsf{ex}}(\sD)$ only \emph{pre}diptych structures, denoted by 
$\un{\ssst{\text{naïve}}}{\ovra{\sT}}(\sD)$ and $\un{\ssst{\text{naïve}}}{\ovra{\sT}_{\ssst\mathsf{ex}}}(\sD)$, 
which are \emph{not} in general 
\emph{diptych} structures\footnote{
Essentially because of the \emph{parallel transfer requirements for pullbacks}. 
This has \emph{important consequences for the existence of pullbacks for Lie groupoids}, which demands some ``spb'' type conditions on morphisms (\emph{unnecessary in the purely algebraic case}). These conditions, which are happily often satisfied in most applications, \emph{are widely ignored in the literature}, where it is very often \emph{implicitely admitted} that, since the pullbacks are always \emph{algebraically} defined for groupoid morphisms, it should be enough to add existence conditions for the \emph{pullbacks of manifolds} in order to get pullbacks for Lie groupoids; but the point is that the sub\-mersion requirements for the source maps demand more care, and our stronger spb type requirements (which indeed are \emph{ubiquitous in various theories}) yield the suitable conditions!
}, 
hence are rather irrelevant for our purpose.

But it turns out that this drawback may be removed when replacing on squares the horizontal \emph{naïve} diptych structure by one of the \emph{main} ones, introduced for this purpose in \ref{nasq}, \emph{following the general process} explained above for defining pre\-diptych structures (\ref{gemeth}), which now may be proved to be actually \emph{diptychs}.

This enables to attach to any small \emph{diptych} $\sT$ (diagram type) several symbols, which now associate functorially a \emph{diptych} to a \emph{diptych} by means of the switch correspondence:

	\[	\boxed{\un{\ssst{\text{main}}}{\ovra{\sT}}
	\ \un{\text{def}}{\longleftrightarrow}
	\ov{\bullet}{\sT}\ \circ\ov{\ssst{\text{vnaïve}}}{\un{\ssst{\text{hmain}}}{\square}},\h
	\un{\ssst{\text{maiN}}}{\ovra{\sT}}
	\ \un{\text{def}}{\longleftrightarrow}
	\,\ov{\bullet}{\sT}\ \circ\ov{\ssst{\text{vnaïve}}}{\un{\ssst{\text{hmaiN}}}{\square}},\h
	\un{\ssst{\text{MaiN}}}{\ovra{\sT}}
	\ \un{\text{def}}{\longleftrightarrow}\,
	\ov{\bullet}{\sT}\ \circ\ov{\ssst{\text{vnaïve}}}{\un{\ssst{\text{hMaiN}}}{\square}}}\,,
	\tag{maindipDia}
	\]
inducing on exact diagrams:
	\[\boxed{\un{\ssst{\text{main}}}{\ovra{\sT}_{\ssst\mathsf{ex}}}
	\un{\text{def}}{\longleftrightarrow}
\ov{\bullet}{\sT}_{\ssst\mathsf{ex}}\circ\ov{\ssst{\text{vnaïve}}}{\un{\ssst{\text{hmain}}}{\square}},\h\,
	\un{\ssst{\text{maiN}}}{\ovra{\sT}_{\ssst\mathsf{ex}}}
	\un{\text{def}}{\longleftrightarrow}
\ov{\bullet}{\sT}_{\ssst\mathsf{ex}}\circ\ov{\ssst{\text{vnaïve}}}{\un{\ssst{\text{hmaiN}}}{\square}},\h\,
	\un{\ssst{\text{MaiN}}}{\ovra{\sT}_{\ssst\mathsf{ex}}}
	\un{\text{def}}{\longleftrightarrow}
	\ov{\bullet}{\sT}_{\ssst\mathsf{ex}}\circ\ov{\ssst{\text{vnaïve}}}{\un{\ssst{\text{hMaiN}}}{\square}}}\,.
	\tag{maindipEx}
\]

If for instance one makes the second definition explicit, this very condensed notation means, using (dipex-Diag switch) of \ref{sqdes} (in order to express the conditions that the source and target have to be vertically written exact dip-diagrams) and our usual switching conventions of \ref{notsub} for horizontal\,/\,vertical laws:
	\[\un{\ssst{\text{maiN}}}{\ovra{\sT}_{\ssst\mathsf{ex}}}(\sD)\,{\longleftrightarrow}\,
	\left(\ov{\bullet}{\sT}_{\ssst\mathsf{ex}}\left(\ov{\ssst{\text{vnaïve}}}{\square} \sD\right);\,
	\ov{\bullet}{\sT}_{\ssst\mathsf{ex}}\left(\ov{\ssst{\text{vnaïve}}}{\un{i}{\square}}\sD\right),
	\ov{\bullet}{\sT}_{\ssst\mathsf{ex}}\left(\ov{\ssst{\text{vnaïve}}}{\un{s}{\gpb}}\sD\right)\right)\,,
\]
where the subscripts i and s mean that the horizontal arrows have to be taken from $\cD_i$ or from $\cD_s$.

These symbols obey certain rules of composition and interchange analogous to (Int.ch) in \ref{iter}, but much more caution is needed for precise statements. For instance, though one can still define (among others) correspondences such as:
	\[\ov{\bullet}{\sS}_{\ssst\mathsf{ex}}
	\circ\un{\ssst{\text{main}}}{\ovra{\sT}_{\ssst\mathsf{ex}}}
	\ttt
	\ov{\bullet}{\sT}_{\ssst\ssst\mathsf{ex}}
	\circ\un{\ssst{\text{main}}}{\ovra{\sS}_{\ssst\mathsf{ex}}}
	\H\text{and}\H
	\un{\ssst{\text{main}}}{\ovra{\sS}_{\ssst\mathsf{ex}}}
	\circ
	\un{\ssst{\text{main}}}{\ovra{\sT}_{\ssst\mathsf{ex}}}
	\ttt
	\un{\ssst{\text{main}}}{\ovra{\sT}_{\ssst\mathsf{ex}}}
	\circ
	\un{\ssst{\text{main}}}{\ovra{\sS}_{\ssst\mathsf{ex}}}\,,
\]
it is no longer possible to identify:
	\[{\ov{\bullet}{\sS}_{\ssst\mathsf{ex}}}\circ\un{\ssst{\text{main}}}{\ovra{\sT}_{\ssst\mathsf{ex}}}
	\H\text{with}\H
\left(\sS\times\sT\right)_{\ssst\mathsf{ex}}^\bullet\ .
\]

\subsubsection{Some important subdiptychs}
\label{sis}

It turns out that most of the previous diptychs \emph{induce} useful \emph{sub\-diptychs} on the subcategories defined by the symbols
	\[\check{\un{h}{\ovra{\cT}}},\h
	\hat{\un{h}{\ovra{\cT}}},\h
	\ov{\diamond}{\un{h}{\ovra{\cT}}}
\]
$(h=i,s)$ introduced in \ref{notsub}. These may be denoted for instance by:
$\un{\ssst{\ \text{main}}}{\Big(\un{s}{\ov{\diamond}{\overrightarrow{\mathsf{T}}}}\Big)_{\ssst\mathsf{ex}}}$, etc.

\subsubsection{Scholium}
\label{scho}

The statements that, starting with a given \emph{diptych}, the previous constructions yield actually \emph{new diptychs} require indeed a very large number of checkings, which cannot be detailed here; just a part of them are straightforward, though none of them is really very hard\footnote{
Though (hyper)cubic diagrams are often basically needed: a \emph{cube} in $\cD$ encapsulates \emph{three commutative squares in} $\sq\cD$ (chasing with the two composition laws and the interchange law), and that is very powerful for deriving properties of the faces, since each face is involved in \emph{two} such commutative squares.
}. %
But they are made \emph{once and for all}, and all these checkings are encapsulated and summarized in the statements. These statements point out highly \emph{remarkable stability properties of diptych structures}.

The point is that, thanks to these stability properties, these constructions may be iterated indefinitely and lead to more and more intricate types of diagrams which it would be very hard to handle and study directly, and to which all the general results established using only diptych axioms apply immediately. This will be used below for various kinds of \emph{structured groupoids}.

\section{Structured groupoids and structured actions}
\label{sgsa}

\subsection{(Symmetric) nerves of small groupoids (Grothendieck)}
\label{grouner}

The (``well known'') characterization of categories and groupoids among simplicial objects, i.e. \emph{contravariant} functors\footnote{
Sometimes called presheaves.
} from $\mathbf{\Delta}\subset\pmb{\mathsf{F}}$ (notations of \ref{simp}) towards \textbf{Set}, by some pushout\,/\,pullback properties of their \emph{nerves} seems to be generally attributed to \textsc{Gro\-then\-dieck} \cite{GRO}, though it does not seem to be very explicit in this reference. (In what follows we often work with the \emph{dual} categories $\pmb{\nabla\subset\pmb{\Finv}}$ in order to avoid contra\-variant functors and co-diptychs whenever possible.)

This algebraic characterization is particularly well adapted for an application of the general ``meta-principle of internalization'' alluded to in \ref{internal}, and enables to transfer ``automatically'' a large variety of algebraic constructions for groupoids to a huge number of structured contexts (including basically Lie groupoids).

Recall that the \emph{nerve} of a (let us say here \emph{small}) \emph{category} $\textbf{G}$ associates (in a functorial way) to each object $[n-1]=\textbf{n}=\{0<\ldots<n-1\}\ (\,\subset\ov{\times}{\textbf{n}}=\left\langle n-1\right\rangle\,)$ of $\pmb{\nabla}\subset\pmb{\Finv}$ (\ref{preoset}, \ref{simp}) the set $G_{[n-1]}$ of \emph{paths} of length $n-1$ (hence $n$ vertices) in $\textbf{G}$, and this actually allows to \emph{identify} $\textbf{G}$ with some (covariant!) functor:
	\[\textbf{G}:\pmb{\nabla}\rightarrow\mathbf{Set},\h\textbf{n}\longmapsto G_{[n-1]}\,,\tag{nerve}
\]
in other words, using the terminology of \ref{nat}, some \emph{diagram of type} $\pmb{\nabla}$ in \textbf{Set}.

We remind also (\ref{exnat}) that each path $u\in G_{[n-1]}$ may be viewed in turn as defining a functor:
	\[[\textbf{u}]:\mathbf{n}\rightarrow \textbf{G}.
\]

Now when the category $\textbf{G}$ is a \emph{groupoid} (i.e. all its arrows are invertible), the path $u$ generates also a \emph{commutative diagram} $\langle u\rangle$ in $\textbf{G}$ with $n$ vertices, which allows to \emph{extend} the previous functor as follows (notations of \ref{sym}):
	\[\langle \textbf{u}\rangle:\ov{\times}{\mathbf{n}}\rightarrow \textbf{G}.
\]
and then the nerve functor \emph{extends} to a functor (which we call here the \emph{symmetric nerve}) defined on the \emph{dual} $\pmb{\Finv}$ of what we called above the \emph{symmetric simplicial category} $\pmb{\mathsf{F}}$ (\ref{sym}):
	\[\textbf{G}:\pmb{\Finv}\rightarrow\textbf{Set},\h\ov{\times}{\textbf{n}}\longmapsto G_{\langle n-1\rangle}\,\tag{sym.nerve}.
\]
In other terms this identifies any groupoid \textbf{G} with some \emph{diagram of type} $\pmb{\Finv}$ in \textbf{Set}.

Now the remarkable \textsc{Grothendieck}'s \emph{characterization of groupoids} among such ``\emph{symmetric simplicial objects}'' is that they \emph{preserve pullbacks} of $\pmb{\Finv}$. Of course these pullbacks come dually from \emph{pushouts} of the symmetric simplicial category $\pmb{\mathsf{F}}$\footnote{
We remind again that many authors use the notation $\Phi$ for our category $\pmb{\sF}$ of finite cardinals.
}.

One may also observe that all the pushouts of $\,\pmb{\mathsf{F}}$ are \emph{generated} (when composing squares) by the ``\emph{elementary squares}'' expressing the classical \emph{canonical relations} between the \emph{canonical generators} $\delta^m_i$ and $\sigma^n_j$ of the \emph{simplicial} subcategory $\mathbf{\Delta}$ (see \cite{ML} for instance) and it is enough to demand preserving the corresponding dual \emph{elementary pullbacks}\footnote{%
Looking more precisely at the structure of $\pmb{\Finv}$ through its dual $\pmb{\mathsf{F}}$ (using for instance its description in \cite{ML}, and some general stability properties of pullback or push\-out squares), one can see that it is enough to preserve the \emph{special} pullbacks (\ref{remax} (4)), and this allows to apply our precise definition of \emph{dip-exactness} introduced in \ref{dipmor}.
} %
of $\pmb{\Finv}$.

Using now the canonical \emph{diptych structure} $\Finv$ on $\pmb{\Finv}$ defined in \ref{sym} (the \emph{dual} of which, $\mathsf{F}$, is the \emph{co}-diptych structure of $\pmb{\sF}$)\footnote{
We insist on the fact that the products in the diptych $\Finv$ we are considering comes from set-theoretic \emph{co}-products of $\pmb{\mathsf{F}}$.
}, this can be expressed quickly in our diptych language by saying that (\ref{dipmor}):
\begin{itemize}
	\item \fbox{groupoids may be identified with \emph{dip-exact diagrams of type} $\Finv\ $ in {\sffamily Set}}.
\end{itemize}

Thus, with the notations of \ref{notsub}, \textsc{Grothendieck}'s characterization allows to \emph{identify} the category of \emph{groupoid morphisms} (\emph{alias} functors between groupoids) with the category of \emph{dip-exact diagram morphisms of type} $\mathsf{\Finv}\ $ in the \emph{standard} diptych $\mathsf{Set\,}$, i.e. with:
	\[\boxed{
	\ovra{\Finv}_{\ssst\mathsf{ex}}(\mathsf{Set})
	}.\tag{Set-Groupoids}
\]

\subsection{Structured groupoids as groupoids in a diptych}
\label{groudip}

\subsubsection{Definition of $\,\protect\sD$-groupoids \emph{[}morphisms\emph{]}}
\label{defgrou}
By application of our meta-principle of internalization, we are now in a position to define \emph{groupoids} [\emph{morphisms}] in any, possibly large and possibly non-concrete (Rem. \ref{remax}), given \emph{diptych} $\sD$ (called the \emph{structuring diptych}), just replacing $\mathsf{Set}$ by $\sD$, as being \emph{dip-exact diagram} [\emph{morphisms}] \emph{of type} $\Finv\ $ in $\sD$. So we shall set, using the definitions (dipex, dipEx) of \ref{notsub}:
	\[\textbf{gpd}(\mathsf{D})
	\un{\text{def}}{=}
	{\overset{\bullet}{\Finv\,}_{\ssst\mathsf{ex}}} (\mathsf{D})
	,\H
	\textbf{Gpd}(\mathsf{D})
	\un{\text{def}}{=}
	{\ovra{\Finv}_{\ssst\mathsf{ex}}}(\mathsf{D}),\tag{D-Groupoids}
\]
or, more shortly and more functorially (omitting the diptych $\sD$):
	\[\boxed{
	\textbf{gpd}
	\un{\text{def}}{=}
	{\overset{\bullet}{\Finv\,}_{\ssst\mathsf{ex}}}}\ 
	,\H
	\boxed{
	\textbf{Gpd}
	\un{\text{def}}{=}
	{\ovra{\Finv}_{\ssst\mathsf{ex}}}}\ .\tag{gpd, Gpd}
\]
The elements of the class $\textbf{gpd}(\mathsf{D})$ are called $\sD-$\emph{groupoids}, and the arrows of the (large) category $\textbf{Gpd}(\mathsf{D})$ are called $\sD-$\emph{groupoid morphisms} (or functors).

\subsubsection{Elementary description.}
\label{edescr}
The exactness property implies that the nerve functor is \emph{uniquely} determined by various subcategories\footnote{
These can be interpreted as sketches (\emph{esquisses}) in \textsc{Ehresmann}'s sense. The comparison of various possible choices with our choice $\Upsilon$ would lead to interesting developments, which cannot be made here. They rely on the remarks in \ref{Finv}, which enable to ``generate'' (in various ways) the objects and arrows of $\pmb{\Finv}$ by means of pullbacks over \textbf{1}.
} %
(more precisely sub-pre\-diptychs) of $\,\mathsf{\Finv}$, notably $\mathsf{\nabla}^{[2]}$ and $\mathsf{\Upsilon}$ (\ref{grouda}, \ref{simp}, \ref{sym}). This allows notably to \emph{identify} a groupoid \textbf{G} with a certain \emph{ diagram of type} $\mathsf{\Upsilon}$:
\begin{center}
\boxed{
\begin{diagram}[w=3em,tight,labelstyle=\ssst]
\Delta G&\rsTo~{\delta_\textbf{G}}&G&\pile{\rsTo~{\alpha_\textbf{G}}\\ \liTo~{\omega_\textbf{G}}\\}&B\\
\end{diagram}}
\end{center}
(the \emph{axioms} of groupoids being then expressed, shortly though abstractly, by the condition that this diagram \emph{extends} to an \emph{dip-exact} diagram of type $\mathsf{\Finv}$, the ``\emph{symmetric nerve}''); in other words $\textbf{G}$ is \emph{uniquely} defined by the sole data:
	\[\boxed{
	(G,B,\,\omega_\textbf{G},\,\alpha_\textbf{G},\,\Delta G,\,\delta_\textbf{G})}\,,
\]
where $B=\textbf{G}_{\langle  0\rangle  }$ (base, objects), $G=\textbf{G}_{\langle  1\rangle  }$ (arrows), $\Delta G=\textbf{G}_{\langle  2\rangle  }$ (\emph{commutative triangles})\footnote{
Remind that, for groupoids, we use what is called here the (extended) \emph{symmetric} nerve $\textbf{G}_{\left\langle n\right\rangle}$ rather than the customary nerve $\textbf{G}_{[n]}$, for which $\textbf{G}_{[2]}$ would mean ``paths of length 2''.
}, and $\alpha_\textbf{G},\,\omega_\textbf{G},\,\delta_\textbf{G}$ (\emph{source} map\footnote{
Of course here the term ``map'' is used loosely for naming an \emph{arrow of} $\cD$, possibly not a set-theoretic map.}, %
\emph{unit} map, \emph{division} map) are determined by the embedding $\Upsilon\subset\pmb{\Finv}$ described in \ref{sym}.
One can check (mimicking the set-theoretic case, in which case one can write $\delta_\textbf{G}(y,x)=yx^{-1}$) that it is then possible to recover, in a purely diagrammatic way, the \emph{projections} $\pi_1,\pi_2:\Delta G\rsTo G$ (or equivalently the embedding $\Delta G\riTo G\times G$), the \emph{inverse} map $\varsigma_\textbf{G}$, the \emph{target} map $\beta_\textbf{G}$ and the \emph{multiplication} (or composition) map $\gamma_\textbf{G}$ (which make up the classical data), as other parts of the \emph{symmetric} nerve.

The arrow $\tau_{\textbf{G}}=(\beta_\textbf{G},\alpha_\textbf{G}):G\rightarrow B\times B$ is generally known as the \emph{anchor map}, following \textsc{K. Mackenzie}, but will be called here, more significantly, the \emph{transitor}, since it encapsulates the (in)\emph{transit}ivity properties of \textbf{G}, and, at least when $B$ is ``$s-$basic'' in the sense defined above in \ref{remax} or below in \ref{ban}, defines a groupoid morphism or funct\emph{or}.

 The squares below are both \emph{s-exact squares} (\ref{Dipdax}) in $\cD$:
\begin{center}
\bdi[size=2em,inline]
\Delta G&\rsTo^{\pi_1}&G\\
\dsTo^{\pi_2} &\gpb&\dsTo_{\alpha_\textbf{G}}\\
G&\rsTo^{\alpha_\textbf{G}}&B\\
\edi\H
\bdi[size=2em,inline]
\Delta G&\rsTo^{\delta_\textbf{G}}&G\\
\dsTo^{\pi_2} &\gpb&\dsTo_{\alpha_\textbf{G}}\\
G&\rsTo^{\beta_\textbf{G}}&B\\
\edi\h.
\end{center}
The left one enables to \emph{identify} $\Delta G$ with $\pmb{\wedge} G$, the \emph{pullback of the pair} $(\alpha_\textbf{G},\alpha_\textbf{G})$.

Note that now $\textbf{G}_{\left\langle 3\right\rangle}=\sq \textbf{G}$ (\emph{commutative squares} in \textbf{G}), which might be regarded also as the \emph{pullback} ${\bigwedge\mspace{-7mu}\!\!\pmb{|}}\ G$ \emph{of three copies} of $\alpha_\textbf{G}$ . Note also that, in the case of groupoids, the projections $\pi^{\text{top}}\,/\,\pi^{\text{bot}}$ of \ref{nasq} become now \emph{equivalences} (to be developed below).

\subsubsection{Co-groupoids}
\label{cogr}
Given a prediptych $\sD$, if $\sD$ happens to be a co-diptych (Rem. \ref{remax} (8)), a $\sD^{\text{op}}-$groupoid will be called a $\sD-$\emph{co-groupoid}.

\subsection{Degenerate D--groupoids.}
One of the interests of the notion of groupoid is to create a link between different notions which appear as degeneracies in various directions.

\subsubsection{Pluri-groups}
\label{plg}
A $\sD-$groupoid \textbf{G} is called a $\sD-$\emph{group} if $\textbf{G}_{\left\langle 0\right\rangle}=\bullet$ (\emph{terminal object}, supposed to exist). This \emph{implies}:
	\[\alpha_\textbf{G}=\beta_\textbf{G}\tag{plg}.
\]

More generally \textbf{G} will be called a $\sD-$\emph{pluri-group}\footnote{
In the case of $\mathsf{Set}-$groupoids, this just means a \emph{sum of groups}. Some authors, such as \textsc{Douady-Lazard} \cite{DOULAZ}, have extended the term \emph{group} to this case, but this may be confusing. The alternative term ``multi-group'' might also create confusion with a special case of \textsc{Ehresmann}'s multiple categories, a quite different notion.
} 
when condition (plg) is satisfied. This is one basic degeneracy of groupoids.

\subsubsection{Null groupoids and null morphisms.}
\label{nul}
For \emph{any} object $B$, the \emph{constant} diagram is a $\sD-$groupoid, denoted by $\mathring{\textbf{B}}$ (and sometimes loosely $B$), called the \emph{null groupoid} associated to $B$: one has $\big(\mathring{\textbf{B}}\big)_{\langle  n-1\rangle  }=B$ for all $n\in\bbN^{\ast}$.

This determines a full and faithful embedding:
	\[\cD\hookrightarrow\textbf{Gpd}(\cD)\,.
\]

A \emph{null morphism} will be defined as a morphism with null source \emph{or} target.

\subsubsection{s-basic objects and banal D-groupoids.}
\label{ban}
Note that, given an object $B$, our axioms do \emph{not} demand the projections $\text{pr}_i:B\times B\rightarrow B$ to be good epis, though this is often true (see Rem. \ref{remax} (2) for sufficient conditions). When this condition is satisfied, the object $B$ will be called \emph{s-basic}. \emph{Then} it defines a $\sD-$groupoid with base $B$, denoted by $\ov{\times}{\textbf{B}}$ (or loosely $B\times B$), called the \emph{banal groupoid}\footnote{
Of course, for above-mentioned reasons, we have to proscribe terminologies such as discrete\,/\,coarse instead of null\,/\,banal, though these are widespread in the literature.
} 
associated to the s-basic object $B$, for which $\Big(\ov{\times}{\textbf{B}}\Big)_{\left\langle n-1\right\rangle}\un{\text{def}}{=}B^n$ (powers of $B$), and, if \textbf{G} is a $\sD-$groupoid with s-basic base $B$, the transitor defines a $\sD-$groupoid morphism:
	\[\tau_{\textbf{G}}:\textbf{G}\rightarrow\ov{\times}{\textbf{B}}\,.
\]

A banal groupoid is a \emph{well exact} (\ref{dipmor}) diagram of type $\Finv$ in $\sD$. When $B$
 is s-condensed (\ref{remax} (2)), this diagram may be ``augmented'', i.e. extended to a diagram of type $\dot{\Finv\ }$.

\subsubsection{Principal D-groupoids and principal morphisms.}
\label{ppal}
In our diptych framework the notion of \emph{regular equivalence} on a manifold $B$ \cite{VAR, SLG} extends immediately, by our general meta-principle, as being the datum of a good epi $B\rsTo^q Q$, i.e. (\ref{loc}) an object of the s-augmented diptych $\vQD$.

The previous construction of \ref{ban} then applies to the s-condensed object $q$ of $\vQD$. Composing with the forgetful functor $\pi_Q^{\text{top}}$, this yields a (degenerate) $\sD-$groupoid \textbf{R}, defined by 
$\textbf{R}_{\left\langle n-1\right\rangle}=\un{\ssst{Q}}{\ov{\ssst{n}}{\times}}B$ (iterated \emph{fibred} product over $Q$), loosely written $R=B\un{Q}{\times}B$. This generalizes the \emph{graph of a regular equivalence relation}, and this is the second basic degeneracy for $\sD-$groupoids.


This groupoid has the property that the (symmetric) nerve may be \emph{augmented}, which means that the exact diagram of type $\Finv$ \emph{extends} to a \emph{still dip-exact} (not well exact in general!) diagram of type $\dot{{\Finv\ }}$: \textbf{0} goes to $Q$ (cf. \ref{sym}).

Note that $q:\textbf{R}\rightarrow\mathring{\textbf{Q}}$ is a null $\sD-$groupoid morphism.

Such an ``\emph{augmented} $\sD-$\emph{groupoid}'' will be called \emph{principal} (over $Q$).

\emph{Principal $\sD-$groupoids over $Q$ derive from banal $\sov{Q}\sD-$groupoids} by the (exact) forgetful functor $\pi_Q^{\text{top}}$ of \ref{loc}.

The \emph{null} groupoid $\mathring{\textbf{B}}$ may be viewed as the special case $q=\text{id}_B$.

The \emph{banal} groupoid of an s\emph{-condensed} (which implies s-basic) object (\ref{remax}) may be viewed as the special case $Q=\bullet$ (terminal object, which is then supposed to exist).

A $\sD-$groupoid \emph{morphism} $\textbf{R}\rightarrow\textbf{G}$ will be called \emph{principal} when its \emph{source} \textbf{R} is a principal $\sD-$groupoid. Basic examples are the right and left ``division morphisms'':
	\[\delta_{\textbf{G}}:\boldsymbol{\Delta}\textbf{G}\longrightarrow\textbf{G}
	\,,\H
{\overline{\delta}}_{\textbf{G}}:\pmb{\nabla}\textbf{G}\longrightarrow\textbf{G}\,.
\]

Given a groupoid $\textbf{G}$, the previous construction applies notably to the source\,/\,tar\-get map $\alpha_\textbf{G} \,/\,\beta_\textbf{G}:\textbf{G}_{\langle  1\rangle  }\rsTo \textbf{G}_{\langle  0\rangle  }$, and this endows $\textbf{G}_{\langle  2\rangle  }$ with \emph{two} canonical \emph{principal} $\sD-$groupoid structures (with base $G$) over $\textbf{G}_{\langle  0\rangle  }$, denoted by $\mathbf{\Delta G}\,/\,\pmb{\nabla} \textbf{G}$.

\subsubsection{ Godement groupoids, Godement's axiom, Godement diptychs}
\label{god}
By its very definition, a principal $\sD-$group\-oid \textbf{R} over $Q$ gives rise to an \emph{s}-exact square (\ref{defdip}):

\begin{center}
\bdi[size=2em,inline]
R&\rsTo^{\beta_{\textbf{R}}}&B\\
\dsTo^{\alpha_{\textbf{R}}} &\gpb&\dsTo_{q}\\
B&\rsTo^{q}&Q\\
\edi\h,
\end{center}
which \emph{implies} the transitor to be a good mono in $\sD$:
	\[\tau_\textbf{R}:R\riTo B\times B\,.\tag{God}
\]

We shall say that a $\sD-$groupoid \textbf{R} is a \emph{Godement $\sD-$groupoid} when its transitor owns that property $\tau_\textbf{R}\in\cD_i$, and that \emph{Godement's axiom} is satisfied in $\sD$ or that $\sD$ is a \emph{Godement diptych} when \emph{conversely} any Godement $\sD-$groupoid is principal.

Following \cite{SLG} (and \cite{VAR}), \emph{Godement's theorem} may be expressed in the present language by the property that $\mathsf{Dif}$ is a \emph{Godement} diptych.

Now we have observed the remarkable fact that most of the diptychs used by ``working mathematicians'' and quoted in the above examples (including toposes and abelian categories) are Godement. 

Moreover it is also remarkable that most of the above constructions for getting diptych structures on categories of diagrams (particularly of groupoids) in a diptych yield new Godement diptychs when starting with a Godement diptych. This cannot be developed and stated more precisely here, and we refer to \cite{P03, P04} for more explanations, and to future papers.

\subsection{Canonical (or universal) groupoids}
\label{cancon}
The definitions may be applied with $\sD=\Finv$, and one has, with the notations of \ref{canend}:
	\[\textbf{gpd}(\Finv)\subset\pmb{\exists}_0
	\H\text {and}\H
	\textbf{Gpd}(\Finv)\subset\pmb{\exists}\,.
\]

$\mathsf{\Finv}-$groupoids [morphisms] are [natural transformations between] exact endo-functors of $\mathsf{\Finv}\,$. They correspond bi\-jectively \emph{by duality } (see \ref{canend})to [natural transformations between] endo\-functors $\Gamma$ of $\,\pmb{\mathsf{F}}$ preserving injections, sur\-jections and push\-outs.

\emph{By right composition} with the symmetric nerves, 
\emph{they operate functorially on $\sD-$groupoids}. (More generally any diagram in the category \textbf{Gpd($\Finv$)} gives rise to a diagram of the same type in \textbf{Gpd}($\sD$)).

For instance (referring to the examples studied in \ref{canend}) one has:
\begin{itemize}
	\item The ``\emph{tautological}'' $\mathsf{\Finv}-$groupoid $\pmb{\mathsf{I}}$, taking for $\Gamma$ the \emph{identical functor} \textbf{Id} of $\pmb{\mathsf{F}}$. It operates identically on $\sD-$groupoids. One has:
	$\pmb{\mathsf{I}}_{\langle  n-1\rangle  }=\mathbf{n}$ (the base is \textbf{1}, the arrow set is \textbf{2}).
	
	It should be carefully noted that, though this groupoid would be banal as a $\dot{\Finv\ }-$groupoid (defined by the \emph{embedding} of $\pmb{\Finv}$ in $\pmb{\dot{\Finv\ }}$), however (in spite of its very naïve appearance) it is \emph{neither banal nor principal} as a $\Finv-$groupoid (which would imply the degeneracy of \emph{all} $\sD-$groupoids!).

	\item The \emph{null} groupoid \textbf{O}, associated to the constant endo\-functor $\boldsymbol{1}$, endowed with the morphism :
	
	\[\omega:\textbf{O}\longrightarrow\pmb{\mathsf{I}}\,,
\]
associated by duality to the natural transformation (null constant map):
	\[\ov{\ast}{\omega}:\textbf{Id}\longrightarrow\textbf{1} ,\textbf{n}\mapsto(\textbf{n}\stackrel{0}{\rightarrow}\textbf{1})\,.
\]
	Its effect on a $\sD-$groupoid \textbf{G} yields the \emph{unit law}:
	\[\omega_\textbf{G}:\textbf{OG}\rightarrow\textbf{G}\,,
\]
where $\textbf{OG}=\mathring{\textbf{B}}$ is the null groupoid attached to the base $B=\textbf{G}_{\langle  0\rangle  }$.

	\item The ``\emph{shift}'' (or ``\emph{triangle}'') $\mathsf{\Finv}-$groupoids $\boldsymbol{\Delta},\pmb{\nabla}$, associated to the ``\emph{shift functors}'' of $\pmb{\mathsf{F}}$, 
endowed with canonical groupoid morphisms:

	\[\delta:\boldsymbol{\Delta}\longrightarrow\pmb{\mathsf{I}}
	\,,\H
{\overline{\delta}}:\pmb{\nabla}\longrightarrow\pmb{\mathsf{I}}\,.\\
\]
These groupoids are \emph{principal} (over \textbf{1}), but \emph{not banal} (see \ref{loc}, \ref{Finv}, and \ref{ppal})!
One has $\boldsymbol{\Delta}_{\langle  n-1\rangle  }=\mathbf{1+n}=\mathbf{n+1}=\pmb{\nabla}_{\langle  n-1\rangle  }$: the base is \textbf{2}, the ``arrows object'' is \textbf{3}.

The action on a $\sD-$groupoid $\textbf{G}$ yields (this explains the notations: see \ref{ppal}) the right and left ``\emph{division}'' morphisms:

	\[\delta_{\textbf{G}}:\boldsymbol{\Delta}\textbf{G}\longrightarrow\textbf{G}
	\,,\H
{\overline{\delta}}_{\textbf{G}}:\pmb{\nabla}\textbf{G}\longrightarrow\textbf{G}\,,
\tag{div}
\]
where $\mathbf \Delta\textbf{G}$ and $\pmb{\nabla} \textbf{G}$ are \emph{principal} (over $G$). When $\sD=\mathsf{Set}$, these are the maps defined by:
	\[\delta_\textbf{G}(y,x)=yx^{-1}\h(\alpha y=\alpha x)\,,\H
	\overline{\delta}_\textbf{G}(y,x)=y^{-1}x\h(\beta y=\beta x)\,.
\]

\item The``\emph{square}'' $\mathsf{\Finv}-$groupoid $\pmb{\square}$\footnote{
Not to be confused with the category $\mathbf{2}\times\mathbf{2}$, for which we used the same symbol in (\ref{preoset}).
}, 
endowed with the following canonical group\-oid morphisms:
	\[\iota:\pmb{\mathsf{I}}\longrightarrow\pmb{\sq},\h
	\pmb{\mathsf{I}}\xleftarrow{\pi^{\text{top}}}\pmb{\sq}
	\xrightarrow{\pi^{\text{bot}}}\pmb{\mathsf{I}}\,.
\]
One has $\pmb{\square}_{\langle  n-1\rangle  }=\textbf{2n}$: the base is \textbf{2} and the arrows object is \textbf{4}.

The action on a $\sD-$groupoid $\textbf{G}$ yields (whence the notations):
	\[\textbf{G}\xrightarrow{\iota_\textbf{G}}{\sq} \textbf{G},\h
	\textbf{G}\xleftarrow{\pi_\textbf{G}^{\text{top}}}{\sq} \textbf{G}\xrightarrow{\pi_\textbf{G}^{\text{bot}}}\textbf{G}\,.
\]
Just as the tautological groupoid, it is neither banal nor principal (the augmented functor is no longer $\pmb{\Finv}\,$-valued, since $\mathbf{0+0}=\textbf{0}$\,!).

\item The``\emph{mirror}'' $\mathsf{\Finv}-$groupoid $\boldsymbol{\Sigma}$, neither banal nor principal, isomorphic to the tautological groupoid by the ``mirror'' isomorphism $\varsigma:\pmb{\mathsf{I}}\rightarrow\boldsymbol{\Sigma}$. The action on a $\sD-$groupoid $\textbf{G}$ yields the \emph{inverse map} (a $\sD-$groupoid isomorphism):
	\[\mathbf{\varsigma_\textbf{G}}: \textbf{G}\rTO^{\approx\ }\textbf{G}^{\text{op}}.
\]
\end{itemize}

In the category $\textbf{Gpd}(\Finv)\subset\pmb{\exists}$, there is a basic diagram of type $\pmb{\bowtie}$:

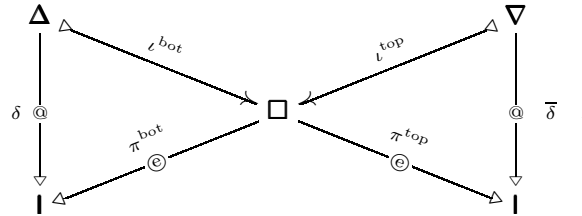
\begin{center}
\begin{diagram}[h=1.8em,w=4.5em,tight,inline,labelstyle=\sst]
\pmb{\Delta}&&&&\pmb{\nabla}\\
&\rdiTo^{\iota^{\text{bot}}}&&\ldiTo^{\iota^{\text{top}}}&\\
\dsTo~{\text{@}}^{\delta}&&\pmb{\square}&&\dsTo~{\text{@}}_{\overline{\delta}}\\
&\ldsTo~{\eq}^{\pi^{\text{bot}}}&&\rdsTo~{\eq}^{\pi^{\text{top}}}&\\
\pmb{\mathsf{I}}&&&&\pmb{\mathsf{I}}\\
\end{diagram}\H.
\end{center}
The tags $\text{@}$ and $\eq$ (which mean actor and equivalence) will be explained below.

The morphisms $\iota^{\text{bot}}$ and $\iota^{\text{top}}$ derive by duality (as explained in \ref{canend}) from the natural transformations associating to the object \textbf{n} the maps $\textbf{\textbf{n}+\textbf{n}}\rightarrow\textbf{1}+\textbf{n}$ and $\textbf{\textbf{n}+\textbf{n}}\rightarrow\textbf{n}+\textbf{1}$ defined respectively by $0_\textbf{n}+\text{id}_{\textbf{n}}$ and $\text{id}_{\textbf{n}}+0_\textbf{n}$, 
where $0_\textbf{n}$ means here the 0 constant map from \textbf{n} to the terminal object $\textbf{1}=\{0\}$ of $\pmb{\sF}$.

Composing on the left with any $\sD-$groupoid \textbf{G} (still interpreted as an exact morphism from $\Finv$ to $\sD$) this diagram  gives rise to the following ``\emph{canonical butterfly diagram}'' (functorial with respect to \textbf{G} and to $\sD-$groupoid morphisms):
	\[\begin{diagram}[h=1.8em,w=4.5em,tight,inline,labelstyle=\sst]
\boldsymbol{\Delta}\textbf{G}&&&&\pmb{\nabla}\textbf{G}\\
&\rdiTo^{\iota_\textbf{G}^{\text{bot}}}&&\ldiTo^{\iota_\textbf{G}^{\text{top}}}&\\
\dsTo~{\text{@}}^{\delta_\textbf{G}}&&\pmb{\sq}\textbf{G}&&\dsTo~{\text{@}}_{\overline{\delta}_\textbf{G}}\\
&\ldsTo~{\eq}^{\pi_\textbf{G}^{\text{bot}}}&&\rdsTo~{\eq}^{\pi_\textbf{G}^{\text{top}}}&\\
\pmb{\textbf{G}}&&&&\pmb{\textbf{G}}\\
\end{diagram}\H.
\tag{can.but}
\]
\begin{itemize}
	\item The transparent geometric notations are of course \emph{suggestive of the set-theoretic case}, where the abstract objects $G_{\langle  n-1\rangle  }$ are sets, the elements of which are pictured by commutative diagrams with $n$ vertices (successively vertices, edges, triangles, squares$,\ldots$), while the maps between these objects are interpreted by iterating operations consisting in dropping some vertex (and the corresponding adjacent arrows), or on the contrary doubling some vertex and the corresponding adjacent arrows, and adding a degenerate unit arrow. The oblique sequences may be regarded as very special cases of short exact sequences of groupoids (cf \cite {P86}).
\end{itemize}

It should be noted that any exact endo\-functor of $\Finv$, notably $\boldsymbol{\Delta}, \pmb{\nabla},\pmb{\square}$, operates also on the former butterfly diagram by right composition, as well as the natural transformations between them (using here the alternative monoidal composition law of $\pmb{\exists}$, called horizontal in \cite{ML} and vertical here), giving rise to a lot of interesting objects. In particular one gets in this way a ``\emph{double butterfly diagram}'' (a diagram of type $\pmb{\bowtie}\times\pmb{\bowtie}$, which it would be difficult to draw explicitly and study directly), which will not be studied here. All of this deserves further studies which are out of our present scope.

\subsection{Various examples} We refer to the examples of diptychs with notations given in section~ \ref{Exdip}, and particularly \ref{large}.
\label{varex}

\begin{enumerate}
	\item Taking successively(see \ref{mainex}): $\sD=\mathsf{Dif},\,\mathsf{DifH},\,\mathsf{Dif}^{(\mathsf{\acute{e}t})},\mathsf{Dif}^{(\mathsf{r.c})},\,\mathsf{Dif}^{(\mathsf{r.1-c})}$, we get notably, among many other variants:
\begin{enumerate}
	\item \textsc{Ehresmann}'s differentiable (or briefly smooth) groupoids, nowadays called \emph{Lie groupoids} (Rem. \ref{Lie pseudogroup});
	\item the \emph{Hausdorff} Lie groupoids (the $\omega_\textbf{G}$ embedding of the base is closed);
	\item the \emph{étale} Lie groupoids (called S-atlases by van Est in \cite{vE}); any (classical) atlas of a manifold $B$ determines a (very special) principal étale groupoid over $B$, which leads to regard general $\sD-$groupoids as \emph{generalized atlases} (see \cite {P89, P04} for a development of this point of view); note that here s-basic means discrete, and the banal groupoids have to be (topologically!) discrete;
	\item the $\alpha-$\emph{connected} Lie groupoids (i.e. the fibres of the source map $\alpha$ are connected: this is the right generalization of connected Lie groups), including the holonomy groupoids of regular foliations (\cite{P66}); here s-basic means connected;
	\item the $\alpha-$\emph{simply connected} Lie groupoids (this is the right generalization of simply connected Lie groups: \cite{P66}).
\end{enumerate}
	\item For $\sD= \mathsf{Top}$ (i.e. using the \emph{open} surjections as good epis, see \ref{large}), one gets those \textsc{Ehresmann}'s \emph{topological groupoids} for which the source (and then target) maps are \emph{open}, which are certainly the most useful ones.
	\item Starting with $\sD=\mathsf{Dif}$ (or variants), let $\mathsf{VB}(\sD)$ be the diptych structure on the category of (morphisms) between smooth \emph{vector bundles}, defined by using spb-squares for defining the good epis: then the $\mathsf{VB}(\sD)-$groupoids are the ``\emph{vector bundle groupoids}'' introduced in \cite{P86} (where a dual object is constructed inside this category). It seems that these groupoids should be the natural framework for studying \emph{linear representations}.
	\item For $\sD=\mathsf{Poset}$, ones gets \textsc{Ehresmann}'s \emph{ordered groupoids}.
	\item One can take for $\sD$ \emph{any topos or abelian category} with its \emph{standard} diptych structure (i.e. \emph{all} monos\,/\,epis are ``good''). This gives rise to an unlimited and very rich range of examples, thus handled in a unified way.
	\item As above-mentioned we hope that a slight modification of our diptych axioms, without weakening their efficiency (introducing for $\cD$ ``good isomorphisms'' and an ``involution''), will enable us to include in the same unified treatment new basic examples, such as Riemannian groupoids, Poisson groupoids, and perhaps measured groupoids.
\end{enumerate}

\subsection{New examples derived from functoriality}%
\label{nex}
(See \cite{P04} for more details and developments).
Thanks to the functorial way we defined $\sD-$groupoids, it is obvious that:
\begin{itemize}
	\item Any \emph{dip-exact} morphism $\Phi:\sD\rightarrow\sD'$ defines (by \emph{left} composition) a functor:
	\[\textbf{Gpd}(\Phi):\textbf{Gpd}(\sD)\rightarrow\textbf{Gpd}(\sD')\,,
\]
and any natural transformation $\varphi:\Phi\stackrel{\centerdot}{\rightarrow}\Phi'$ defines a natural transformation:
	\[\textbf{Gpd}(\varphi):\textbf{Gpd}(\Phi)\stackrel{\centerdot}{\rightarrow}\textbf{Gpd}(\Phi')
\]
This may be applied for instance to many forgetful functors, whenever they are dip-exact.
\item Any [natural transformation between] endo-functors of $\mathbf{DipEx}$ defines a [natural transformation between] endo\-functors of \textbf{Gpd}.

For instance, with the notations of the above example (3), one can prove that the ``\emph{tangent functor}'' (which associates to any manifold $B$ its tangent bundle $\text{T}B$) defines indeed such an endo\-functor of $\mathbf{DipEx}$ (the choice of our good epis using (spb) condition was made for this purpose):
	\[\textbf{T}:\sD\mapsto \mathsf{VB}(\sD),
\]
where $\sD$ may be one of the variants for $\mathsf{Dif}$, 
and we have natural transformations (``zero section'' and ``canonical projection to the base'') connecting \textbf{T} with the ``zero functor'' \textbf{0} (by which any manifold $B$ may be \emph{identified} with a \emph{null} vector bundle $\textbf{0}B$):
	\[\textbf{0}\stackrel{0}{\rightarrow}\textbf{T}\stackrel{\text{t}}{\rightarrow}\textbf{0}.
\]
Then we know immediately that, for any $\sD-$groupoid \textbf{G}, its tangent space has a canonical structure of $\mathsf{VB}(\sD)-$groupoid\footnote{
On the contrary it is by no means obvious (as proved by \textsc{A. Weinstein} \emph{et alii} \cite{CDW}, using symplectic structures) that the \emph{co\-tangent} space owns such a structure. In \cite{P88}, this is derived from a more general quite unexpected and remarkable general duality property of 
 $\textbf{Gpd}(\mathsf{VB(D)})$, which does not use any symplectic structure.
}, 
with canonical morphisms:
\[\textbf{G}\stackrel{\textbf{0}_{\textbf{G}}}{\rightarrow}\textbf{TG}\stackrel{\textbf{t}_{\textbf{G}}}{\rightarrow}\textbf{G}\,,
\]
	and that this construction is functorial in \textbf{G} (identified with $\textbf{0}\textbf{G}$).
\end{itemize}

\subsection{D-groupoids morphisms; D-actors, D-action laws}
\label{groumor}

\subsubsection{Description of D-groupoid morphisms and notations.}
\label{desmor}

We are in the general situation described in \ref{sqdes}, \ref{maindip}, with many special properties deriving from $\sT=\Finv$. We cannot make an exhaustive study here but just gather some needed facts.

We fix some \emph{shortened notations}. The diptych $\sD$ is fixed. Let $\textbf{f}:\textbf{H}\rightarrow\textbf{G}$ be a $\sD-$groupoid morphism. 
We define (for $n\in\bbN$) $H_n=\textbf{H}_{\langle  n\rangle  },\ G_n=\textbf{G}_{\langle  n\rangle  },\ E=H_0,\ H=H_1,\ \Delta H=H_2,\ B=G_0,\ G=G_1,\ \Delta G=G_2$. 
We have a family of horizontal ``connecting'' arrows of $\cD$, $\textbf{f}_{\left\langle n\right\rangle}=f_n:H_n\rightarrow G_n$, which turns out to be \emph{uniquely} determined by $f=f_1$. They determine ``connecting'' squares in $\cD$, which will be called $hi$/$vi$, or $hs$/$vs$ when their horizontal\,/\,vertical arrows belong to $\cD_i$ or $\cD_s$. All these squares are \emph{generated} by vertical composition of $vi$ or $vs$ squares (more precisely the vertical arrows may be taken as being the opposite (or dual) of the canonical generators of the simplicial category).

A $\sD-$\emph{natural transformation} may be defined as $\sD-$morphism from \textbf{H} to $\sq\textbf{G}$.

\subsubsection{ Special morphisms: hyper- and hypo-actors, (partial) action laws.}
\label{spemor}
It turns out that $f\in\cD_i$ implies $f_n\in\cD_i$ for all $n$ (then all the connecting squares are ipb-squares, and $\textbf{f}$ will be called an \emph{i-morphism}), and $f\in\cD_s$ implies $f_0\in\cD_s$ (\textbf{f} will be called an \emph{s-morphim} if \emph{all} the connecting arrows $f_n$ belong to $\cD_s$). These define \emph{subcategories} of \textbf{Gpd}, denoted by \textbf{iGpd} and \textbf{sGpd}.

In order all the $vs$ connecting squares (\ref{desmor}) to be respectively gpb, ipb, spb (\ref{tbp}), it is enough that this property is satisfied by the square $\textbf{a}(\textbf{f})$ below:
\begin{center}
\bdi[w=2.8em,h=2.8em,labelstyle=\ssst,inline]
H&\rTO~{f}&G\\
\dsTo~{\alpha_\textbf{H}}&\textbf{a}(\textbf{f})&\dsTo~{\alpha_\textbf{G}}\\
E&\rTO~{f_0}&B\\
\edi\h.
\end{center}
(In the gpb and ipb cases, one can prove that the property is then actually satisfied by \emph{all} the connecting squares).
Then the $\sD-$groupoid morphism \textbf{f} is called respectively an @-morphism (or $\sD-$\emph{actor}), $\check{@}-$morphism (or $\sD-$\emph{hypo-actor}), $\widehat{@}-$mor\-phism (or $\sD-$\emph{hyper-actor})\footnote{
We slightly modify the terminology (in-actors, ex-actors) of our previous papers. \textsc{Ehresmann}'s terminology for our hypo-actors was ``well faithful'', while our actors were called ``foncteurs d'hypermorphismes'' in his book ``Catégories et structures'' (see our remark \ref{universe}).
}. 
A (hyper)-actor is an s-morphism iff $f_0\in\cD_s$, and is then called an \emph{s-hyper-actor}, tagged by $\rsTo^{\widehat{@}\ }$.
\begin{itemize}
	\item The (seemingly very natural) term \emph{actor} (funct\emph{or} describing an \emph{act}ion law, substituted here to the proscribed ``classical'' terminology ``discrete (op)fibration''!)\footnote{
While our \emph{hyper-actors} are known classically as ``(op)fibrations'', in the categorical sense.
} 
stresses in a straightforward way that, when $\sD=$ {\sffamily Set}, such a morphism or \emph{functor} $f:H\rightarrow G$ is associated (\emph{in an essentially bijective way}) to an \emph{action law} $\lambda$ of the groupoid \textbf{G} on ($f_0:E\rightarrow B$), obtained by composing $\lambda=\beta_\textbf{H}\circ\varphi^{-1}$, as shown below:

\begin{center}
\bdi[w=6em,h=2em,tight,labelstyle=\sst,inline,midshaft]
&{}&\hLine~{\lambda}[abut]&{}&\\
\ldLine(1,1)[abut] G\un{B}{\times}E&\lTO~{\approx}^{\varphi=(f,\,\alpha_\textbf{H})}&H&\rsTo^{\beta_\textbf{H}}&E\rdTO(1,1)[abut]\\
\edi
\h\h.
\end{center}
Conversely the groupoid structure on the domain of $\lambda$ is derived from its graph embedding\footnote{
We modify the order of the factors in order to \emph{systematically write the source on the right and the target on the left}.
}:
	\[\widehat {\lambda}:G\un{B}{\times}E\rTO G\times(E\times E),\h
	(g,x)\longmapsto(g,(\lambda(g,x),x))\,.
\]
(More precisely $\widehat {\lambda}$ takes its values in the fibre product $G\un{B\times B}{\times}(E\times E)$).
The definition of a groupoid action law is precisely expressed by the property that the image of $\widehat {\lambda}$ be a sub\-groupoid of the groupoid 
$\textbf{G}\un{\ov{\times}{\textbf{B}}}{\times}\ov{\times}{\textbf{E}}$.
In the general case, according to our general meta-principle of internalization, we shall \emph{define} $\lambda=\beta_\textbf{H}\circ\varphi^{-1}$ as the $\sD-$\emph{action law} associated to the $\sD-$\emph{actor} $f:H\to G$. The most important case is when one has $(E\rsTo B)$ in $\cD_s$ ($\sD-s$-actors, $\sD-s$-action laws), but the definitions remain meaningful in the general case.

This generalization of group action laws was introduced by \textsc{Ehresmann}.
\item The terminology of \emph{hypo-actor} is also very natural, since the definition means that the previous map $\varphi$, instead of being an isomorphism, is now a good mono. Then the map $\beta_{\textbf{G}}= \lambda\circ\varphi$ may be interpreted as a (special) \emph{partial action law}: see \ref{univac} below for the problem of globalizing such a partial action law.
\end{itemize}

Such special morphisms will be tagged by the symbols @, $\check{@}$, $\widehat{@}$; they define subcategories of \textbf{Gpd} and \textbf{sGpd}, denoted by @\textbf{Gpd}, $\check{@}$\textbf{Gpd}, $\widehat{@}$\textbf{Gpd}, and \textbf{s}@\textbf{Gpd}, $\textbf{s}\check{@}$\textbf{Gpd}, $\textbf{s}\widehat{@}$\textbf{Gpd}.

We refer to \cite{P89} for various statements which can be proved diagrammatically in the diptych context, notably concerning the existence of pullbacks; in particular the pullbacks along any s-hyper-actor \textbf{f} \emph{always do exist}.

\subsubsection{ Special morphisms: inductors, equivalences, and extensors}
\label{eqext}

We wish now to ``internalize'' the algebraic notions of fullness and faithfulness.

Analagously to the square \textbf{a}(\textbf{f}) of \ref{spemor}, we can consider the square \textbf{t}(\textbf{f}) built from the transitors (keeping the notations of \ref{desmor}, and setting $p=f_0:E\rightarrow B$):

\begin{center}
\bdi[w=5em,h=2.8em,labelstyle=\sst,inline,midshaft]
H&\rTO~{f}&G\\
\dTO~{\tau_{\textbf{H}}}&\textbf{t}(\textbf{f})&\dTO~{\tau_{\textbf{G}}}\\
E\times E&\rTo~{p\times p}&B\times B\\
\edi\h.
\end{center}

Then $\textbf{f}:\textbf{H}\rightarrow\textbf{G}$ will be called: \emph{i-faithful}, \emph{s-full}, a $\sD-$\emph{inductor}, depending on whether the square \textbf{t}(\textbf{f}) is ipb, spb, gpb (\ref{tbp}). Note that any hypo-actor is i-faithful.

In the third case we can say \textbf{H} is \emph{induced} by \textbf{G} along $f_0=p:E\rightarrow B$. When moreover one has $p\in\cD_s$, \textbf{f} will be called an \emph{s-equivalence} (or \emph{s-inductor}).

To ``internalize'' the general algebraic notion of equivalence, we need a diagrammatic description of \emph{essential} (or \emph{generic}) \emph{sur\-jectivity}. This is achieved by means of the following diagram, which always does exist:

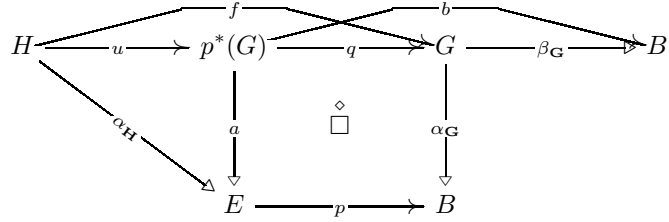
\begin{center}
\begin{diagram}[w=2em,h=1.5em,tight,labelstyle=\scriptstyle,textflow,midshaft,inline]
&&&{}&\rLine[abut]~ f&{}&&  {}&\rLine[abut]~ b&{}&&&\\
H&&\ruLine(3,1)[abut]\rTO~ u&&p^{\ast}(G)&&\rdTO(3,1)[abut]\rTO~ q\ruLine(3,1)[abut]&&G  &&\rsTo~{\beta_{\textbf{G}}}\rdTO(3,1)[abut]&&B\\
&\rdsTo(4,4)~{\alpha_{\textbf{H}}}&&&&&&&  &&&&\\
&&&&\dsTo~ a&&\gpb&&\dsTo~{\alpha_{\textbf{G}}}  &&&&\\
&&&&&&&&  &&&&\\
&&&&E&&\rTO~ p&&B  &&&&\\
\end{diagram}\H.
\end{center}
Then \textbf{f} will be called \emph{essentially-s} whenever one has $b\in\cD_s$. One can prove that this \emph{implies} the existence of a $\sD-$groupoid $\textbf{H'}=p^{\ast\ast}(\textbf{G})$ \emph{induced} by \textbf{G} along $p$, with a canonical factorization $\textbf{h}:\textbf{H}\rightarrow\textbf{H'}$. Then (mimicking the algebraic definition) \textbf{f} is called a $\sD-$\emph{equivalence} whenever \textbf{h} is an isomorphism and \textbf{f} is essentially-s. When moreover \textbf{f} is an i-morphism, one speaks of an \emph{i-equivalence}.

A $\sD-$equivalence will be tagged by $\rTO^{\eq\ }[midshaft]$.

Finally \textbf{f} will be called an \emph{s-extensor} whenever it is s-full and moreover one has $p\in\cD_s$. This \emph{implies} that \textbf{f} is an s-hyper-actor (\ref{spemor}) (hence an s-morphism).

An s-extensor will be tagged by $\rsTo^{\widehat{\eq}\ }[midshaft]$.

We explain now why this (non-classical) terminology seems very natural.

For an s-extensor \textbf{f}, the kernel \textbf{N} exists as a $\sD-$groupoid, and is defined by the following pullback (which turns out to be a pushout too):

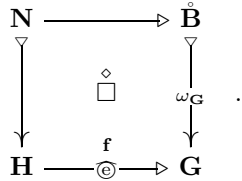
\begin{center}
\bdi[w=3.2em,h=2.8em,labelstyle=\sst,inline,midshaft]
\textbf{N}&\rsTo&\mathring{\textbf{B}}\\
\diTo&\gpb&\diTo~{\omega_{\textbf{G}}}\\
\textbf{H}&\rsTo^{\textbf{f}}~{\widehat{\eq}}&\textbf{G}\\
\edi\h.
\end{center}

It is a \emph{normal}\footnote{
It is worth noting that, for a subgroupoid \textbf{N}, the algebraic condition of normality bears only on the isotropy groups of \textbf{N}.
} sub\-groupoid of \textbf{H} (and principal (\ref{ppal}) in the case of an s-equivalence). Then 
	\[\mathring{\textbf{E}}\riTo \textbf{N}\riTo\textbf{H}\rsTo\textbf{G}
\]
may be called a \emph{short exact sequence} of $\sD-$groupoids, since, from the underlying algebraic point of view, this is the \emph{obvious exact generalization}\footnote{
Which however seems strangely widely despised or underestimated in the literature, in spite of its basic simplicity and importance among the various kinds of more sophisticated quotients (systematically studied notably in \textsc{Ehresmann}'s book, in the more general context of categories, as solutions of various universal problems, and, in the smooth framework, in \cite{MK}, where the study of our s-hyperactors, named fibrations, is privileged).
} of the quotient of a group by a normal subgroup (the only slight difference is that one has to consider the \emph{two-sided} cosets $NxN$ instead of the left or right ones in the case of groups, or more generally plurigroups (\ref{plg})). The purely diagrammatic proof of these facts given in \cite{P86} remains valid in any Godement diptych. Then the funct\emph{or} \textbf{f} describes the \emph{extens}ion of \textbf{G} by \textbf{N}, whence our terminology.
\subsubsection{$\protect\sD-$groupoid morphisms as groupoids in $\protect\sq\protect\sD$}
\label{gpdsq}

The general discussion of \ref{sqdes} applies with $\sT=\Finv$, though some of the statements below require extra verifications.

The basic correspondence (dipEx switch) of \ref{sqdes} reads now:
	\[\textbf{Gpd}(\sD)
	\longleftrightarrow
	\textbf{gpd}\Big(\ov{\ssst{\text{vnaïve}}}{\square}\sD\Big)\,,
\]
and means that giving a $\sD-$groupoid \emph{morphism} is equivalent to giving a \emph{groupoid} structured by the category of \emph{squares} of $\sD$, endowed with its \emph{naïve} diptych structure. (Note that even for $\sD=\mathsf{Set}$, this gives an interesting algebraic interpretation of groupoid morphisms).

Using \ref{nasq} and \ref{masq}, we can get also interpretations for the \emph{special morphisms} introduced in \ref{spemor} ($\sD$ is omitted):

\begin{eqnarray}
	@\textbf{Gpd}&\leftrightarrow
	\textbf{gpd}\Bigg(\ov{\ssst{\text{vcan}}}{\gpb}\Bigg)\,,&
	\check{@}\textbf{Gpd}&\leftrightarrow
	\textbf{gpd}\bigg(\ov{\ssst{\text{vnaïve}}}{\ipb}\bigg)\,,&
	\widehat{@}\textbf{Gpd}&\leftrightarrow
	\textbf{gpd}\Big(\ov{\ssst{\text{vcan}}}{\square}\Big)\,,\\
	\textbf{s}@\textbf{Gpd}&\leftrightarrow
	\textbf{gpd}\,\Bigg(\un{s}{\ov{\ssst{\text{vcan}}}{\gpb}}\Bigg)\,\,,&
\textbf{i}@\textbf{Gpd}&\leftrightarrow
	\textbf{gpd}\,\Bigg(\un{i}{\ov{\ssst{\text{vcan}}}{\gpb}}\Bigg),&
	\textbf{s}\widehat{@}\textbf{Gpd}&\leftrightarrow
	\textbf{gpd}\,\Big(\un{s}{\ov{\ssst{\text{vcan}}}{\square}}\Big)\,.\notag
\end{eqnarray}

\subsubsection{Canonical diptych structures on the category of \emph{D}-groupoids.}
\label{dipgpd}
Applying now the discussion of \ref{maindip} (maindipEx) with $\sT=\Finv$, we can get various \emph{diptych structures} on the category of groupoid morphisms and on various subcategories, for which we introduce special notations :

\begin{align}
	\mathsf{Gpd}&
	\un{\text{def}}{=}
	\un{\ssst{\text{main}}}{\textbf{Gpd}}=
	\big(\textbf{Gpd};\,\textbf{iGpd},\,\textbf{s}\widehat{@}\textbf{Gpd}\big)&&
	\longleftrightarrow\H
	\mathbf{gpd}\,\circ\ov{\ssst{\text{vnaïve}}}{\un{\ssst{\text{hmain}}}{\sq}}\label{Gmain}\\
	\mathsf{Gpd}_{(@)}&
	\un{\text{def}}{=}
	\un{\ssst{\text{maiN}}}{\textbf{Gpd}}=
	\big(\textbf{Gpd};\,\textbf{iGpd},\,\textbf{s}@\textbf{Gpd}\big)&&
	\longleftrightarrow\H
	\mathbf{gpd}\,\circ\ov{\ssst{\text{vnaïve}}}{\un{\ssst{\text{hmaiN}}}{\sq}}\,.\label{GmaiN}
\end{align}
For instance these formulas mean that, for the \emph{main diptych structure} on the category of $\sD-$groupoid morphisms, one should take as \emph{good epis} the \emph{s-hyper-actors} \eqref{Gmain}, while it is possible also to take the \emph{s-actors} \eqref{GmaiN}. One may also observe, for more symmetry, that the good monos will be automatically hypo-actors, i.e. one may replace \textbf{iGpd} by $\textbf{i}\check{@}\textbf{Gpd}$ and write 
$\mathsf{Gpd}={_{i\check{@}}}\textbf{Gpd}_{\,s\widehat{@}}$.

There are also several interesting \emph{subdiptychs}. Among others:

\begin{align}
	\widehat{@}\mathsf{Gpd}&
	\un{\text{def}}{=}
	\big(\,\widehat{@}\textbf{Gpd};\,\textbf{i}@\textbf{Gpd},\,\textbf{s}\widehat{@}\textbf{Gpd}\,\big)&&
	\longleftrightarrow\H
	\mathbf{gpd}\,\circ\un{\ssst{\text{hcan}}}{\ov{\ssst{\text{vcan}}}{\widehat{\square}}}\label{subhyper}\\
		\widehat{@}\mathsf{Gpd}_{(@)}&
	\un{\text{def}}{=}
	\big(\,\widehat{@}\textbf{Gpd};\,\textbf{i}@\textbf{Gpd},\,\textbf{s}@\textbf{Gpd}\,\big)&&
	\longleftrightarrow\H
	\mathbf{gpd}\,\circ\un{\ssst{\text{hMaiN}}}{\ov{\ssst{\text{vcan}}}{\widehat{\square}}}\label{SuBhyper}\\
	@\mathsf{Gpd}&
	\un{\text{def}}{=}
	\big(\,@\textbf{Gpd};\,\textbf{i}@\textbf{Gpd},\,\textbf{s}@\textbf{Gpd}\,\big)&&
	\longleftrightarrow\H
	\mathbf{gpd}\,\circ\ov{\ssst{\text{vcan}}}{\un{\ssst{\text{hcan}}}{\gpb}}\label{suba}\,.
\end{align}
This means notably that the \emph{canonical diptych structure} (which is just induced by the \emph{naïve} one \emph{in those special cases}) is suitable for the special \emph{subcategories of actors}, and \emph{hyper-actors}. (Actually the extra condition demanded for good monos, reduces indeed to just being an i-morphism. Remind also, for more symmetry: $\textbf{i}\widehat{@}\textbf{Gpd}=\textbf{i}@\textbf{Gpd}$).

It is also highly remarkable that starting with a \emph{Godement} diptych $\sD$, the previous constructions yield a variety of new Godement diptychs, but this cannot be developed here.

\subsection{Double groupoids}
\label{doub}
We cannot discuss here thoroughly the general problem of ``internalizing'' the notion of double groupoids, which certainly deserves future research and further developments. We just mention briefly a few tracks; these lead to a huge variety of unequivalent more or less restrictive notions, which are of interest even in the purely algebraic (i.e. set-theoretic) context.

Using the ``double nerve'', one may consider pre\-dip-diagrams $\Finv\times\Finv\rightarrow\sD$, preserving ``horizontal'' and ``vertical'' pullbacks. This defines two $\sD-$groupoid structures satisfying \textsc{Ehresmann}'s interchange law.

Alternatively (but here in a more restrictive way), following the basic observation made by \textsc{Ehresmann} in the purely algebraic context, one may also iterate or compose the functors $\mathsf{Gpd}$ or $\mathsf{Gpd}_@$ introduced in \ref{dipgpd}. But even when taking $\sD=\mathsf{Set}$, this defines interesting classes of \emph{special} double groupoids (considered by various authors).

The objects of $(\mathsf{Gpd})^2(\mathsf{Set})$ satisfy certain special ``filling conditions'' considered notably by Ronnie \textsc{Brown}.

The more special objects of $(\mathsf{Gpd}_@)^2(\mathsf{Set})$ 
(for which double arrows are \emph{uniquely defined} by their two sources, or alternatively by the source for one law and the target for the other law) are particularly important and interesting in our present context. In the special case when moreover the two groupoid laws are both principal, they have been called affinoids by A. \textsc{Weinstein} and pregroupoids by A. \textsc{Kock}'s, who presents them in his lecture \cite{K} in the underlying algebraic context (see also \cite{P85} and \cite{P03}). These objects own a \emph{third category structure}, the composition law of which (\emph{mixed} or \emph{oblique law}) is suggested by the following diagram:
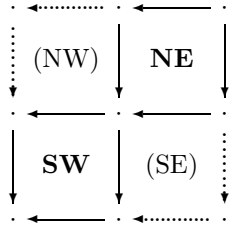
\begin{center}
\begin{diagram}[size=2em,tight,labelstyle=\scriptstyle,inline]
\cdot&\lDotsto&\cdot&\lTo&\cdot\\
\dDotsto&(\text{NW})&\dTo&\textbf{NE}&\dTo\\
\cdot&\lTo&\cdot&\lTo&\cdot\\
\dTo&\textbf{SW}&\dTo&(\text{SE})&\dDotsto\\
\cdot&\lTo&\cdot&\lDotsto&\cdot\\
\end{diagram}\ ,\emph{}
\end{center}
where the dotted arrows and the corresponding (NW )and (SE) ``squares'' (picturing ``double arrows'') are \emph{built} by the basic property; the \emph{oblique composite} of the \emph{given} \textbf{NE} and \textbf{SW} ``squares'' is then the big ``square'' defined by applying the \emph{exchange law} to the four small ``squares''.

In the present context these notions become meaningful (and may be studied in a purely diagrammatic way) when replacing $\mathsf{Set}$ by any diptych $\sD$.

\section{Principal actions; conjugate and associate actions}
\label{sec:pact}
It was initially planned to give longer development to the present section, which illustrates the general methods described above. Since the preliminary sections are more meaty than expected, we have to postpone further details to future papers, and will be content with presenting the main leading ideas.
\subsection{Diagrammatic description of Ehresmann's constructions}
\label{diagdesc}

\subsubsection{Classical principal fibrations}
\label{pfib}

According to the description, as given in \cite{VAR}, of a (smooth) \emph{principal fibration} in \textsc{Ehresmann}'s original sense, let $p_0:P\rightarrow B$ be a principal fibration with base $B$ and \emph{structural group} $\textbf{H}$ (a Lie group with underlying manifold $H$) acting\footnote{
Contrary to the tradition, we consider only \emph{left} action laws, using the fact that any right action law $(h,x)\mapsto x\cdot h$ may be described by the left action law $h\cdot x=x\cdot h^{-1}$. A right action law of $\textbf{H}$ may be also considered as a left action law of $\textbf{H}^{\text{op}}$, using the canonical isomorphism $\varsigma_{\textbf{H}}:\textbf{H}\stackrel{\approx}{\rightarrow}\textbf{H}^{\text{op}},\ h\mapsto h^{-1}$.
} 
smoothly and \emph{freely}\footnote{
And properly when one wants to deal only with Hausdorff manifolds, but the definitions (and constructions) remain meaningful when these assumptions are dropped simultaneously.
} on $P$.

The base $B$ is the orbit space of this action, which defines on $P$ a regular equivalence relation; the \emph{graph} of this equivalence relation is the manifold $S=P\un{B}{\times}P$, and we have a (smooth) \emph{isomorphism}:
	\[\varphi^{-1}:H\times P\stackrel{\approx}{\rightarrow}S,\ (h,x)\mapsto(h\cdot x,x)\,.
\]
Let us set $q=\text{pr}_1\circ\varphi:S\rightarrow H$.

Now, relying on \emph{Godement's theorem} (\cite{SLG, VAR}), and using the language introduced in \ref{spemor} and \ref{ppal}, it is readily seen that we can summarize, in a highly condensed and very simple way, all the data, definitions and axioms (including the local triviality conditions), as given in \cite{VAR}, by saying that $S$ has a \emph{principal} $\sD-$groupoid structure \textbf{S} (taking $\sD=\mathsf{Dif}$, or $\mathsf{DifH}$) (\ref{mainex}) with base $S_0=P$, and:
	\[\boxed{
	\textbf{q}:\textbf{S}\rsTo^{@\ }\textbf{H}
	}%
	\tag{ppl.act}
\]
is a \emph{principal} $\sD-$\emph{s-actor}\footnote{
One should note that, in our terminology, ``\emph{principal}'' includes simultaneously an \emph{algebraic} condition (the action is \emph{free}), and several \emph{topological} (more precisely smooth) conditions encapsulated via Godement's theorem and the precise definitions of structured groupoids and actors. Briefly said, we introduce the notion of $\sD-$\emph{principal action} as the ``diptych internalization'' of the purely algebraic notion of \emph{free action}, via a diagrammatic description of such an algebraic object. The goal of this section is to show that this encapsulates a very rich extra information!
} (\ref{spemor}, \ref{ppal}).
One observes that the canonical factorization of $q:S\stackrel{j}{\rightarrow}K=H\times(P\times P)\stackrel{v}{\rightarrow}H$, where $j=(q,\tau_\textbf{S})$ and $v=\text{pr}_1$, may be read in the category of smooth groupoids:
	\[	\boxed{
	\bdi[labelstyle=\sst,w=4em,midshaft,inline]
\textbf{S}&\riTo^{\textbf{j}\ }&\textbf{K}=\textbf{H}\times\ov{\times}{\textbf{P}}&\rsTo~{\eq}^{\textbf{v}}&\textbf{H}
\edi}\ .\tag{left wing}
\]
\textbf{S} may be viewed as a \emph{normal} sub\-groupoid embedded by \textbf{j} in the (trivial and trivialized) groupoid \textbf{K} (since the condition of normality bears only on the isotropy groups, which vanish here, \textbf{S} being principal), and \textbf{v} as a (trivial) \emph{s-equivalence} of groupoids (\ref{eqext}). We note also that one has $K\approx S\times P$ (as manifolds).

\subsubsection{Diagrammatic description of Ehresmann's structural groupoid}
\label{stgpd}
Remind that \textsc{Ehres\-mann}'s structural groupoid \textbf{G} may be defined as the quotient of $P\times P$ by the \emph{diagonal action} of \textbf{H}; more precisely, as a $\sD-$groupoid, it may be regarded as a quotient of the \emph{banal} (\ref{ban}) groupoid $\textbf{R}=\ov{\times}{\textbf{P}}$. Its base is $G_0=B$ (quotient of $P$ by the action of \textbf{H}).
Moreover the canonical projection defines a \emph{principal} $\sD-$\emph{s-actor}:
	\[\boxed{
	\textbf{p}:\textbf{R}\rsTo^{@\ }\textbf{G}
	}%
	\tag{conj.act}
\]
and one checks immediately that the associated \emph{action law} (\ref{spemor}) is precisely the (transitive) action law of \textbf{G} on $(p_0:P\rightarrow B)$ introduced by \textsc{Ehresmann}. One has $R\approx G\un{B}{\times}P$.

The canonical factorization of \textbf{p} is now:

	\[	\boxed{
	\bdi[labelstyle=\sst,w=4em,midshaft,inline]
\textbf{R}&\riTo^{\textbf{i}\ }&\textbf{K'}=\textbf{G}\un{\ov{\times}{\textbf{B}}}{\times}\ov{\times}{\textbf{P}}&\rsTo~{\eq}^{\textbf{u}}&\textbf{G}
\edi}\ .\tag{right wing}
\]
with $\textbf{j}=(\textbf{p},\tau_{\textbf{R}})$, and $\textbf{u}=\text{pr}_1$.

Now we let the reader check (using associativity of fibre products) that one has a \emph{canonical isomorphism} (by means of which we \emph{identify} \textbf{K} and \textbf{K'}):
	\[\boxed{
	\textbf{K}\approx\textbf{K'}
	}\,.
\]
This allows to gather the previous left and right wings into a beautiful \emph{butterfly diagram} (precisely a diagram of type $\pmb{\bowtie}$ (\ref{but}) in $\textbf{Gpd}(\mathsf{Dif})$):
	\[
\begin{diagram}[h=1.8em,w=4.5em,tight,inline,labelstyle=\sst]
\textbf{S}&&&&\textbf{R}\\
&\rdiTo^{\textbf{j}}&&\ldiTo^{\textbf{i}}&\\
\dsTo~{\text{@}}^{\textbf{q}}&&\textbf{K}&&\dsTo~{\text{@}}_{\textbf{p}}\\
&\ldsTo~{\eq}^{\textbf{v}}&&\rdsTo~{\eq}^{\textbf{u}}&\\
\textbf{H}&&&&\textbf{G}\\
\end{diagram}\H.
\tag{conj}
\]

The symmetry of the diagram may be achieved by setting $H_0=C=\bullet$ (the base of a group is a singleton!) and noting that $q_0$ is the canonical map $P\rsTo\bullet\,$. This allows to write products $\times$ as fibred products $\un{\bullet}{\times}=\un{C}{\times}\,$. As explained in the introduction, the fact that the corresponding equivalence relation is the coarse one makes it ``invisible'', and hides the symmetry.

The main \emph{source of dissymmetry} lies here in the fact that $q_0$ is \emph{split} (though \emph{non-canonically}!), while $p_0$ is not in general: when it is, the principal fibration is called trivial\emph{izable} (a trivial\emph{ization} being a \emph{choice} of a section of $p_0$). One may say that one wing is trivial while the other one is only locally trivial. In a certain sense the local triviality of this original situation leads to wrong tracks that hide the fundamental simplicity of the situation, which will be more apparent in the general situation as described below.

Moreover one observes that $\textbf{S}=\text{Ker}\,\textbf{u}$ and $\textbf{R}=\text{Ker}\,\textbf{v}$, so that the two diagonals may be viewed as (degenerate) \emph{short exact sequences of $\sD-$groupoids} (\ref{eqext}).

This means that the structural groupoid \textbf{G} might be also described (up to a uniquely defined isomorphism) as the quotient of the heart \textbf{K} of the butterfly by the normal subgroupoid \textbf{S}.

We shall express this relation between the principal actions of the \emph{structural group} and the \emph{structural groupoid} onto the principal bundle $P$ by saying that they are \emph{conjugate}. This is indeed a very special case of a general symmetric relation described below.

We shall not go further without restoring the symmetry, and forgetting the (local) transitivity and triviality conditions. Even the notion of associate bundles will become much clearer in the general context.

Moreover the above presentation, suitable for an immediate application of our meta-principle of internalization, allows now to forget and bypass the set-theoretic constructions, and leads to very wide-ranging notions in the diptych framework.

\subsection{General conjugation of principal morphisms}
\label{cpa}
From now on $\sD$ will be a \emph{Godement diptych} (\ref{god}). For the sake of simplicity, we shall also add to the definition of $\sD-$groupoids  the condition that the base is an s-basic object (\ref{remax} (2), \ref{ban}) (when $\sD=\mathsf{Dif}$ this would eliminate just the empty groupoid). This allows to use freely the banal groupoid associated to any object of $\cD$, and to consider the transitors $\tau$ (or anchor maps) as $\sD-$groupoid morphisms.

For general notations and conventions concerning a $\sD-$groupoid morphism $\textbf{f}:\textbf{H}\rightarrow\textbf{G}$, we refer to \ref{desmor}. In particular \textbf{f} determines a map between the bases $f_0:E\rightarrow B$, and a map between the arrows objects $f=f_1:H\rightarrow G$.

Most of the situations presented here are indeed special instances of more general ones, introduced in \cite{P89} in order to interpret the generalized morphisms in the sense of Skandalis-Haefliger \cite{H}, as \emph{fractions} obtained by \emph{inverting the s-equivalences} (and actually all the $\sD-$equivalences). More precisely we encounter here only the \emph{generalized isomorphisms}, also known as \emph{Morita equivalences}. Though written in the differentiable framework, the presentation given in \cite{P89} was designed in order to extend to the Godement diptychs framework.

\subsubsection{Butterfly diagram of a principal morphism}
\label{cdpm}
We assume here \textbf{H} is \emph{principal} (\ref{ppal}) and this is stressed by writing now $\textbf{H}=\textbf{R}$, $R_0=P$, and $\textbf{f}=\textbf{r}:\textbf{R}\rightarrow \textbf{G}$. This means we have 
$(\tau_{\textbf{R}}:R\riTo P\times P)\in\cD_i$
. We suppose moreover one has: 
	\[(r_0:P\rsTo B)\in\cD_s\,.
\]
We have not to assume \textbf{r} to be an actor for being able to build a butterfly diagram (more precisely a diagram of type $\pmb{\bowtie}$ in the category $\textbf{Gpd}(\sD)$) as in the previous subsection:

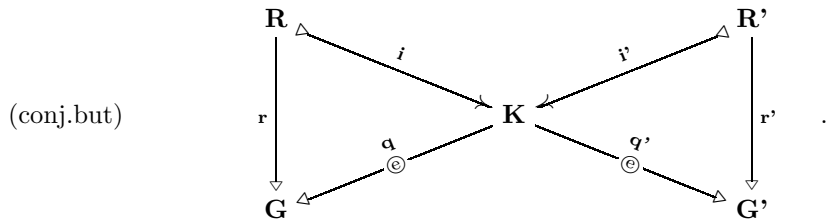
\[
\begin{diagram}[h=1.8em,w=4.5em,tight,inline,labelstyle=\sst]
\textbf{R}&&&&\textbf{R'}\\
&\rdiTo^{\textbf{i}}&&\ldiTo^{\textbf{i'}}&\\
\dsTo^{\textbf{r}}&&\textbf{K}&&\dsTo_{\textbf{r'}}\\
&\ldsTo~{\eq}^{\textbf{q}}&&\rdsTo~{\eq}^{\textbf{q'}}&\\
\textbf{G}&&&&\textbf{G'}\\
\end{diagram}\tag{conj.but}\H.
\]
It is defined as follows:
\begin{itemize}
	\item $\textbf{i}=(\textbf{r},\tau_{\textbf{R}}):\textbf{R}\rightarrow\textbf{K}=\textbf{G}\un{\ov{\times}{\textbf{B}}}{\times}\ov{\times}{\textbf{P}}$\,;
	\item $\textbf{q}=\text{pr}_1$ (it is an s-equivalence);
	\item $\textbf{G'}=\textbf{K}//\textbf{R}$ (which denotes the \emph{two-sided quotient} of $\textbf{K}$ by the normal subgroupoid \textbf{R} (\ref{eqext})), and the projection $\textbf{q'}:\textbf{K}\rightarrow \textbf{G'}$ is an s-equivalence;
	\item $\textbf{R'}= \text{Ker}\,\textbf{q}$, with the embedding \textbf{i'}.
\end{itemize}
The base of \textbf{K} has two projections: $B\lsTo^{\ q_0\ }P\rsTo^{q'_0\ }B'$ (with $q_0=r_0$, $q'_0=r'_0$).
We have again oblique short exact sequences of $\sD-$groupoids, and each wing is determined by the other one up to a unique isomorphism.

\subsubsection{Conjugation butterflies}
\label{conbut}
More generally a diagram of type $\pmb{\bowtie}$ as above will be called a \emph{conjugation butterfly} whenever:
\begin{itemize}
	\item \textbf{q} and \textbf{q'} are s-equivalences (\ref{eqext}), with $\textbf{R}=\text{Ker}\,\textbf{q'}$, $\textbf{R'}=\text{Ker}\,\textbf{q}$ (\textbf{i} and \textbf{i'} being the canonical embeddings).
\end{itemize}
The (principal) $\sD-$morphisms \textbf{r} and \textbf{r'} are then called \emph{conjugate} (by means of this diagram).
By definition of s-equivalences we have two \emph{isomorphisms}:
	\[\bdi[w=5em,tight]
\textbf{G}\un{\ov{\times}{\textbf{B}}}{\times}\ov{\times}{\textbf{P}}&\lTO~{\ \approx}^{(\textbf{q},\tau_{\textbf{K}})}&\textbf{K}&\rTO~{\approx}^{(\textbf{q'},\tau_{\textbf{K}})\ }&\textbf{G'}\un{\ov{\times}{\textbf{B'}}}{\times}\ov{\times}{\textbf{P}}\\
	\edi
\]

\subsubsection{Conjugation functor and its properties}
\label{cf}
The previous construction is functorial: to a square in $\sq\textbf{Gpd}(\sD)$ connecting two s-principal morphisms, there corresponds functorially a diagram morphism of type $\pmb{\bowtie}$ in $\textbf{Gpd}(\sD)$ (which would be pictured by a ``thickened'' butterfly, adding a third dimension).

It is readily seen that such a conjugation butterfly diagram determines a \emph{self-adjoint functor} in the full sub\-category of $\sq\textbf{Gpd}(\sD)$ generated by the s-principal $\sD-$morphisms (viewed as special objects of this category).

So, to any s-principal morphism (or to any square connecting two such s-principal morphisms) this functor associates an s-principal morphism (or a square), which will be called its \emph{mirror image}.

This functor has very remarkable properties, which may be proved diagrammatically (as interesting exercises!), and might be checked more directly in the set-theoretic case, i.e. when $\sD=\mathsf{Set}$. These properties lead to very interesting and powerful constructions, which certainly have to be exploited more extensively; we are going to sketch just a few of them.
\begin{itemize}
	\item When \textbf{G} is \emph{s-transitive} (i.e. $\tau_{\textbf{G}}\in\cD_s$), so is \textbf{G'}. When $\sD=\textbf{Dif}$, this characterizes the transitive locally trivial smooth groupoids (anciently called Lie groupoids till 1987: cf. \ref{Lie pseudogroup})\footnote{
It should be noted that, in the smooth case, the properties of \emph{local triviality} are reflected (relying on Godement's theorem) by a very simple \emph{property of the transitor} $\tau$, which allows an immediate ``internalization'' in our framework of Godement diptychs. In the purely topological context, there is no such simple characterization, and the local triviality becomes a cumbersome condition, with very poor stability properties.}. 
The correspondence between the structural group and the structural groupoid, recalled in \ref{pfib}, is a very special instance of this general property. However in our present framework, this appears as an \emph{extra property}, which is not basic for most constructions.
\item When \textbf{r} is a principal \emph{s-actor}, so is \textbf{r'}; this is the situation alluded to in the introduction and on which we focus now.
\end{itemize}

\subsubsection{Conjugation of principal actors}
\label{cpac}
We shall not develop further the study of conjugation for the general principal morphisms, though it may be of interest.

Indeed the situation becomes geometrically much richer (and in direct connection with our present subject) in the case of principal \emph{actors}. Actually the conjugate of an actor is an actor, as mentioned just above.

As stated in \cite{P89}, this situation is \emph{characterized} by the (symmetric) condition that $\delta_{\textbf{K}}$ \emph{induces} an isomorphism:
	\[\delta_{\textbf{K}}:\Delta\textbf{K}\supset\textbf{R}\un{\alpha}{\times}\textbf{R'}\stackrel{\approx}{\rightarrow}\textbf{K}\,.
	\tag{transv}
\]
Then \textbf{R} and \textbf{R'} are called in \cite{P89} \emph{transverse} subgroupids of \textbf{K}.

This gives rise to a very rich geometrical situation which can be described in various equivalent ways. Indeed the underlying \emph{algebraic} situation is exactly the situation described in \textsc{A.Kock}'s lecture \cite{K} (using the language of torsors and pregroupoids), but our diagrammatic description allows to transfer it into any Godement diptych, while enlightening the computations by a geometric vision.

Actually the data of the two \emph{principal} subgroupoids \textbf{R} and \textbf{R'} of \textbf{K}, satisfying such a \emph{transversality} condition, define on the underlying object $K$ of the groupoid \textbf{K} an \emph{extra structure}, which may be described equivalently as a rule of three on its base $P$, or alternatively as this special double groupoid structure arising in \ref{doub} as the structure of an object of $(\mathsf{Gpd}_@)^2(\mathsf{Set}))$, admitting the groupoid structure of \textbf{K} as its \emph{mixed law}.

This cannot be detailed in the present account, which aims at stressing the leading principles arising from \textsc{Ehresmann}'s concepts. We shall be content with reprinting below the following figure 
(Fig. \ref{fpac}), extracted from \cite{P03}, which pictures the underlying set-theoretic situation, with its strong geometric perfume, while suggesting the diagrammatic description. 

The remarkable fact is that the whole structure arises from the data of principal (i.e. very degenerate) groupoids. (In the same way, one notes that any groupoid structure \textbf{G} arises from a degenerate groupoid $\Delta\textbf{G}$, and an equivalence relation on $\Delta \textbf{G}$, to know $\delta_{\textbf{G}}:\Delta\textbf{G}\rsTo\textbf{G}$, i.e. again a degenerate groupoid, thus illustrating the slogan that all the structured groupoids arise from principal ones and from Godement's axiom).

\begin{figure}[htb]
\def\KK{{
\begin{diagram}[size=1.2em,p=.5em,tight,thin]
\textbf{K}\\ \dsTo\dsTo\\P
\end{diagram}}}
\def\gg{{
\begin{diagram}[size=1.2em,p=.5em,tight,thin]
\textbf{G'}\\ \dsTo\dsTo\\B'
\end{diagram}}}
\def\gh{{\begin{diagram}[w=1.3em,p=0.5em,tight,labelstyle=\scriptstyle,thin]
\textbf{G'}&\pile{\rsTo[abut]\\ \rsTo[abut]}&B'\\
\end{diagram}}}
\def\GG{{
\begin{diagram}[size=1.2em,p=.5em,tight,thin]
\textbf{G}\\ \dsTo\dsTo\\B
\end{diagram}}}
\fbox{
\begin{diagram}[w=1.6em,h=1.6em,p=.3em,tight,labelstyle=\scriptstyle,abut,midshaft,inline,thick]
&&&&&&&&&&&&&&&&&&&&\\
{}&\rLine&{}&\rLine&{}&\rLine&{}&&&&&{}&&&&&&&&\\
\dLine&&\dLine&&\dLine&&\dLine&&&\gh&&\dLine&&&&&&&&&\\
{}&\rLine&\sst{y}\bullet\phantom{\sst{y}}&\rLine&\phantom{\sst{x}}\bullet\sst{x}&\rLine&{}&{}&&&&\phantom{\scriptstyle{\alpha(g')}}\bullet\scriptstyle{\alpha(g')}&&&&&&&&&\\
\dLine&&\dLine&\ldTo~{k}&\dLine&&\dLine&\dTo~{g'}&&&{}&\dLine&&\dLine\dLine&\left(\KK\right)&&&\rsTo_{q'}~{\eq\,}&&&\left(\gg\right)\\
{}&\rLine&\sst{y'}\bullet\phantom{\sst{y'}}&\rLine&\phantom{\sst{x'}}\bullet\sst{x'}&\rLine&{}&{}&&&&\phantom{\scriptstyle{\beta(g')}}\bullet\scriptstyle{\beta(g')}&&&&&&&&&\\
\dLine&&\dLine&&\dLine&&\dLine&&&&&\dLine&&&\dsTo~{\eq}_q&&&&&&\\
{}&\rLine&{}&\rLine&{}&\rLine&{P}&{}&&\rsTo~{q'_0}&{}&B'&&&&&&&&&\\
&&{}&\lTo~{g}&{}&&{}&&&&&&&&&&&&&&\\
&\GG&{}&&&&\dsTo~{q_0}&&&&&&&&\left(\GG\right)&&&&&&\\
&&&&&&{}&&&&&&&&&&&&&&\\
{}&\rLine&{\overset{{\beta(g)}}\bullet}&\rLine&\overset{{\alpha(g)}}\bullet&\rLine&B&&&&&&&&&&&&&&\\
&&&&&&&&&&&&&&&&&&&&\\
\end{diagram}}
\begin{tabular}{||c||}
\hline
$y=g\lhd x,\,x'=g'\triangledown x$\\
$y'=g\lhd x'=g'\triangledown y$\\
$\alpha(k)=x,\,\beta(k)=y'$\\
$q_0^{\ast\ast}(\textbf{G})\approx \textbf{K}\approx{q'_0}^{\ast\ast}(\textbf{G'})$\\
$g=q(k),\,g'=q'(k)$\\
\hline
\end{tabular}
\caption{Conjugation of two principal actors.}\label{fpac}
\end{figure}

Now we display some remarkable examples of the corresponding properties of a square and its mirror (i.e. conjugate) image (the vertical arrows are s-principal morphisms, connected by the horizontal arrows); the reflection is horizontal:

	\[\bdi[s=3em,inline,tight,midshaft,labelstyle=\sst]
\textbf{S}&\rTO&\textbf{R}\\
\dsTo~{@}^{\textbf{s}}&\gpb&\dsTo~{@}_{\textbf{r}}\\
\textbf{H}&\rTO^{\textbf{f}}&\textbf{G}\\
\edi
\HH\leftarrow\Biggr|\!\rightarrow\HH
\bdi[s=3em,inline,tight,midshaft,labelstyle=\sst]
\textbf{S'}&\rTO~{@}&\textbf{R'}\\
\dsTo~{@}^{\textbf{s'}}&&\dsTo~{@}_{\textbf{r'}}\\
\textbf{H'}&\rTO~{@}^{\textbf{f'}}&\textbf{G'}\\
\edi\tag{a}
\]

	\[\bdi[s=3em,inline,tight,midshaft,labelstyle=\sst]
\textbf{S}&\rTO~@&\textbf{R}\\
\dsTo~{@}^{\textbf{s}}&\gpb&\dsTo~{@}_{\textbf{r}}\\
\textbf{H}&\rTO~@^{\textbf{f}}&\textbf{G}\\
\edi
\HH\leftarrow\Biggr|\!\rightarrow\HH
\bdi[s=3em,inline,tight,midshaft,labelstyle=\sst]
\textbf{S'}&\rTO~@&\textbf{R'}\\
\dsTo~{@}^{\textbf{s'}}&\gpb&\dsTo~{@}_{\textbf{r'}}\\
\textbf{H'}&\rTO~@^{\textbf{f'}}&\textbf{G'}\\
\edi\tag{b}
\]

	\[\bdi[s=3em,inline,tight,midshaft,labelstyle=\sst]
\textbf{S}&\rTO~{\eq}&\textbf{R}\\
\dsTo~{@}^{\textbf{s}}&&\dsTo~{@}_{\textbf{r}}\\
\textbf{H}&\rTO~{\eq}^{\textbf{f}}&\textbf{G}\\
\edi
\HH\leftarrow\Biggr|\!\rightarrow\HH
\bdi[s=3em,inline,tight,midshaft,labelstyle=\sst]
\textbf{S'}&\rTO~{@}&\textbf{R'}\\
\dsTo~{@}^{\textbf{s'}}&&\dsTo~{@}_{\textbf{r'}}\\
\textbf{H'}&\rTO^{\textbf{f'}}~{\approx}&\textbf{G'}\\
\edi\tag{c}
\]

\subsection{A few examples and applications}
\label{exap}

\subsubsection{The canonical butterfly}
\label{canbut}
The canonical butterfly (canbut), encountered in \ref{cancon}, and which we reproduce here for convenience, 
\begin{diagram}[h=1.8em,w=4.5em,tight,inline,labelstyle=\sst,textflow]
\boldsymbol{\Delta}\textbf{G}&&&&\pmb{\nabla}\textbf{G}\\
&\rdiTo^{\iota_\textbf{G}^{\text{bot}}}&&\ldiTo^{\iota_\textbf{G}^{\text{top}}}&\\
\dsTo~{\text{@}}^{\delta_\textbf{G}}&&\pmb{\sq}\textbf{G}&&\dsTo~{\text{@}}_{\overline{\delta}_\textbf{G}}\\
&\ldsTo~{\eq}^{\pi_\textbf{G}^{\text{bot}}}&&\rdsTo~{\eq}^{\pi_\textbf{G}^{\text{top}}}&\\
\pmb{\textbf{G}}&&&&\pmb{\textbf{G}}\\
\end{diagram}
may be obtained by applying the previous construction to the canonical principal actor $\delta_{\textbf{G}}:\Delta\textbf{G}\rightarrow\textbf{G}$ attached to any $\sD-$groupoid \textbf{G}. Its conjugate is $\overline{\delta}_{\textbf{G}}$. The corresponding action laws are the left and right\footnote{
More precisely the left action law associated to the right translation by the process explained above.
}
 translations of \textbf{G} on $B\stackrel{\alpha_{\textbf{G}}}{\leftarrow}G\stackrel{\beta_{\textbf{G}}}{\rightarrow}B$ (the symmetry is fully apparent in this example, which, though trivial, is basic).

\subsubsection{Associate bundles}
\label{ass}
Let us look at the second mirror correspondence (b) between squares (which is a consequence of property (a)). Considering \textbf{r} and \textbf{r'} as fixed, we get (by building pullbacks) a correspondence (conjugation) between (possibly non-principal) \textbf{G}-actors \textbf{f} (or action laws) and \textbf{G'}-actors \textbf{f'} (or action laws).

In the classical case recalled in \ref{pfib}, one can check that we recover in this way the correspondence between a given action of the structural group on a ``\emph{fibre type}'' $E$ (base of \textbf{H}), and the \emph{associate action}, introduced by \textsc{Ehresmann}, of the structural groupoid on the so-called ``\emph{associate bundle}'' (see also \cite{VAR} for instance, for the classical presentation), which here is constructed as \emph{the base $E'$ of} \textbf{H'}. 

In the general case, these associate action laws play perfectly \emph{symmetric roles}, and are determined by one another.

Note that, when this second action, corresponding to the second actor \textbf{f}, is itself also principal, the properties of the squares are preserved by transposition. This gives rise to a second mirror image (obtained by ``vertical'' reflections). This ``kaleidoscopic'' situation leads to very interesting questions, which we intend to examine elsewhere (its seems that a triple periodicity arises).

\subsubsection{Conjugating horizontal arrows}
\label{conjh}
Given any (``horizontal'') arrow $\textbf{f}:\textbf{H}\rightarrow\textbf{G}$ of $\textbf{Gpd}(\sD)$, we can pullback the (``vertical'') \emph{canonical} principal actor $\delta_{\textbf{G}}:\Delta\textbf{G}\rsTo~{@}\textbf{G}$ along \textbf{f}, and consider the mirror image of this pullback square (not a pullback in general when \textbf{f} is no longer an actor!) and this defines a correspondence $\textbf{f}\mapsto\textbf{f}^{\text{c}}$, which has very interesting properties.

First we note that when \textbf{f} is any \textbf{G}-\emph{actor}, i.e. an actor with target \textbf{G} (not supposed to be principal here), this correspondence defines a new \textbf{G}-actor by property (b), and this determines an interesting involution which will not be studied here.

We give now another important application to the ``activation'' of a groupoid morphism, relying on properties (a) and (c).

\subsubsection{``Activation'' of a \emph{D}-groupoid morphism}
\label{univac}
We start with a definition. 
Given a $\sD-$groupoid morphism $\textbf{f}:\textbf{H}\rightarrow\textbf{G}$ we define an \emph{activation} of \textbf{f} as a (commutative) diagram:

	\[
\bdi[w=10em,h=2em,tight,midshaft,inline,labelstyle=\sst]
&\textbf{H}_{\textbf{1}}&\\
\textbf{H}\ruTO~{\textbf{h}_{\textbf{1}}}(1,1)&\rTO~{\textbf{f}}&\textbf{G}\rdTO~{@}^{\textbf{f}_\textbf{1}}(1,1)\\
\edi\H,
\tag{activ}
\]
where (as indicated by the tagging) $\textbf{f}_{\textbf{1}}$ is an actor.

Such an activation $\widehat{\textbf{H}}$ is called \emph{universal} if it factorizes any other given activation \textbf{H'}, as shown by the dotted arrow in the following diagram:

	\[\bdi[w=10em,h=2em,tight,midshaft,inline,labelstyle=\sst]
&\textbf{H'}&\\
\ruTO~{\textbf{h'}}(1,3)&\dTO[dotted]&\rdTO~{@}^{\textbf{f'}}(1,3)\\
&\widehat{\textbf{H}}&\\
\textbf{H}\ruTO~{\widehat{\textbf{h}}}(1,1)&\rTO~{\textbf{f}}&\textbf{G}\rdTO~{@}^{\widehat{\textbf{f}}}(1,1)\\
\edi\H.
\tag{univac}
\]
The following wide-ranging result (where the above notations are used) extends an algebraic one to the structured context:

\begin{thm} Let \emph{$\textbf{f}:\textbf{H}\rightarrow\textbf{G}$} be an \emph{i-faithful} $\sD-$morphism. Then:
\begin{enumerate}
	\item \emph{\textbf{f}} admits an (essentially unique) \emph{universal activation} \emph{$\textbf{f}=\widehat{\textbf{f}}\circ\widehat{h}$}, for which moreover \emph{$\widehat{\textbf{h}}$} is an \emph{equivalence}.
\item \emph{\textbf{f}} is a \emph{hypo-acto}r iff \emph{$\widehat{\textbf{h}}$} is an \emph{i-equivalence}.
\item\emph{ \textbf{f}} is a \emph{hyper-actor} iff \emph{$\widehat{\textbf{h}}$} is an \emph{s-equivalence}.
\end{enumerate}
\end{thm}

We just indicate briefly how to build the diagram on which the proof of (1) relies.

One starts with the big gpb square (written on the left) obtained by pulling back $\delta_{\textbf{G}}$ along \textbf{f}, and one considers its right mirror image (to which property (a) applies). It is not in general a pullback, but may be factorized into a gpb square and a triangularly degenerate (on the bottom edge) square. Then one considers the left mirror image of this factorization: this yields a horizontal factorization of the big left square, and property (c) applies. Then the decomposition of the bottom arrow is the desired one. Look at the following scheme below. The universal property is left to the reader, as well as the other statements.

	\[\bdi[s=2.9em,inline,tight,midshaft,labelstyle=\sst]
&{}&\hLine[abut]&{}&\\
\ruLine(1,1)[abut]\textbf{S}&\rTO~{\eq}&\textbf{R}&\rTO~{@}&\Delta\textbf{G}\rdTO[abut](1,1)\\
\dsTo~{@}_{\textbf{s}}&\gpb&\dsTo~{@}^{\textbf{r}}&\gpb&\dsTo~{@}^{\delta_{\textbf{G}}}\\
\textbf{H}&\rTO~{\eq}^{\widehat{\textbf{h}}}&\widehat{\textbf{H}}&\rTO~{@}^{\widehat{\textbf{f}}}&\textbf{G}\\
&{}\rdLine(1,1)[abut]&\hLine^{\textbf{f}}[abut]&{}\ruTO(1,1)[abut]&\\
\edi
\h\leftarrow\Biggr|\!\rightarrow\h
\bdi[s=3em,inline,tight,midshaft,labelstyle=\sst]
&{}&\hLine~{@}[abut]&{}&\\
\ruLine(1,1)[abut]\textbf{S'}&\rTO~{@}&\textbf{R'}&\rTO~{@}&\nabla\textbf{G}\rdTO[abut](1,1)\\
\dsTo~{@}_{\textbf{s'}}&&\dsTo~{@}^{\textbf{r'}}&\gpb&\dsTo~{@}^{\overline{\delta}_\textbf{G}}\\
\textbf{H'}&\rDoubleLine&\textbf{H'}&\rTO~{@}^{\textbf{f'}}&\textbf{G}\\
&{}\rdLine(1,1)[abut]&\hLine~{@}^{\textbf{f'}}[abut]&\ruTO(1,1)[abut]&\\
\edi
\]

\subsection{A few examples of universal activations}
\label{exactiv}

\subsubsection{Homogeneous spaces}
\label{homsp}

When applying the previous construction to the inclusion of an embedded subgroupoid $\textbf{G'}\subset\textbf{G}$, one just gets the canonical action of \textbf{G} on the (\emph{one-sided}!) quotient of \textbf{G} by \textbf{G'}, which does exist (using Godement's axiom) as an object of $\cD$ (not a groupoid), and which extends obviously the theory of homogeneous spaces in the case of groups. In particular, $\delta_{\textbf{G}}$ may be regarded as the universal activation of $\omega_{\textbf{G}}$.

\subsubsection{Palais' globalization of partial action laws}
\label{pal}

A much more interesting example is a generalization of \textsc{R. Palais'} theory of globalization of partial action laws \cite{Pal}. It is important to note that this theory holds in $\mathsf{Dif}$ (not $\mathsf{DifH}$!) (\ref{mainex}). Of course this causes no problem in our framework. This might be one reason of extending the definitions of principal fibrations and the basic constructions to non-Hausdorff manifolds. (Classically these are rejected, by demanding the action of the structural group to be proper \cite{VAR}, as leading to too exotic objects; however these objects do exist, and are nowadays recognised to be important, notably in non-commutative geometry. Of course, with our presentation, all the constructions which rely only on our Godement diptych properties remain valid).

It turns out that the condition introduced by \textsc{Palais} for his so-called \emph{globalizable partial actions}, may be expressed by the fact that the ``graph embedding'' $\hat{\lambda}$ we described in \ref{spemor} for (global) action laws still defines a $\sD-$subgroupoid \textbf{H}. Now $\varphi:H\riTo G\un{B}{\times}E$ is no longer an isomorphism in general, but an arrow of $\cD_i$ [even more, in \textsc{Palais}' context, it is an open (but \emph{non-closed}!) embedding (but that is not necessary in our framework)]; this is precisely, in this $\mathsf{Dif}$ context, the definition of our \emph{hypo-actors} (which is the special case (2) of the previous theorem), and one can check that \textsc{Palais}' globalization (though it may look very specific and \emph{ad hoc}) may be actually regarded as a special instance of the previous very general result.

\subsubsection{Cocycles and Haefliger cocycles}
\label{cocy}
We go on working in $\sD=\mathsf{Dif}$, and not $\mathsf{DifH}$ (since we need open embeddings).
An open covering $\mathcal{U}=(U_i)_{i\in I}$ of a (possibly non-Hausdorff) manifold $B$ may be described diagrammatically in the following way. We set $U=\coprod_{i\in I}U_i$, which has an \'etale surjective projection on $B$ and we regard $\textbf{R}=U\un{B}{\times}U$ as a principal (\'etale) $\sD-$groupoid (with base $U$) (which is equivalent to the null groupoid $\mathring{\textbf{B}}$).

Following the terminology of \cite{VAR}, a (non-abelian) cocycle, attached to $\mathcal{U}$, with values in \textbf{G} (a Lie group), is nothing else than a principal $\sD-$morphism of $\sD-$groupoids
	\[\boxed{
	\textbf{g}:\textbf{R}\rightarrow \textbf{G}}\,.
\]

Moreover two cocycles are \emph{cohomologous} \cite{VAR} iff the corresponding $\sD-$mor\-phisms are \emph{isomorphic} (in the sense of $\sD-$natural transformations between functors, see \ref{desmor}).

Now such a formulation makes sense for any kind of $\sD-$groupoid \textbf{G}. So, taking for \textbf{G} a \emph{pseudogroup} $\Gamma$, which, as explained in \ref{remextra}, may always be completed (by considering its \emph{germs}) and then identified with an (\'etale!) $\sD-$groupoid, this general definition contains also, by another specialization, the notion of the so-called \emph{Haefliger cocycles}.

Coming back to the case of groups, we let the reader check that the universal activation of \textbf{g} yields precisely ``the'' (the quotation marks are to remind that one has to work up to isomorphisms) principal fibration associated to the isomorphy class of \textbf{g}.

The construction applies also to Haefliger cocycles, and is related to the construction of the associated foliation and its holonomy pseudogroup, but precise formulations require some more care, and that will be another story.

We would hope to tease some reluctant geometers with showing that ``abstract nonsense'' may have unexpected connections with very concrete Geometry, which of course \textsc{Ehresmann} would knew perfectly, though sometimes taking pleasure in keeping these connections secrete, and kept for Zeus.
\section{Conclusion}
\label{concl}

The main messages we would like to stress, stemming from \textsc{Ehresmann}'s legacy, would be the following:

\begin{itemize}
	\item Graphs of equivalence relations and groups are dissymetrical degeneracies of groupoids in which a part of deep reality becomes invisible, and the continuity of the link between various incarnations of this reality is broken.
	\item Life is possible in ``bad'' categories, and it may even be comfortable and pleasant.
	\item If you are able to write diagrams, then you can make Geometry while making Algebra, and, in this way, live two lives simultaneously: this is the essence of \textsc{Ehresmann}'s meta-principle of ``internalization''.
\end{itemize}

\providecommand{\bysame}{\leavevmode\hbox to3em{\hrulefill}\thinspace}

\end{document}